\documentclass[3p, 12pt, preprint]{elsarticle}

\makeatletter
\def\ps@pprintTitle{%
    \let\@oddhead\@empty
    \let\@evenhead\@empty
    \let\@evenfoot\@oddfoot
    }
\usepackage{tabularx}
\usepackage{amsmath}
\usepackage{amssymb}
\usepackage{amsthm}
\usepackage{mathrsfs}
\usepackage{dsfont}
\usepackage{bm}
\usepackage[usenames,dvipsnames]{color}
\usepackage{graphicx}
\usepackage{float}
\usepackage{enumitem}
\usepackage{array}
\usepackage{multirow}
\usepackage{makecell}
\usepackage{empheq}
\usepackage{cases}
\usepackage{textcomp}
\usepackage{gensymb}
\usepackage{subcaption}
\usepackage{hyperref}
\hypersetup{
	colorlinks=true,       %
	linkcolor=blue,        %
	citecolor=blue,        %
	filecolor=magenta,     %
	urlcolor=blue    
}
\usepackage[capitalise,noabbrev]{cleveref}
\usepackage{etoolbox}

\AtBeginDocument{\setlength{\bibsep}{0pt}}

\crefname{equation}{}{}

\usepackage{placeins}

\newcommand{\bA}{\mathbb{A}}
\newcommand{\bB}{\mathbb{B}}
\newcommand{\bC}{\mathbb{C}}
\newcommand{\bE}{\mathbb{E}}
\newcommand{\bI}{\mathbb{I}}
\newcommand{\bK}{\mathbb{K}}
\newcommand{\bL}{\mathbb{L}}
\newcommand{\bN}{\mathbb{N}}
\newcommand{\bP}{\mathbb{P}}
\newcommand{\bR}{\mathbb{R}}
\newcommand{\bT}{\mathbb{T}}
\newcommand{\bV}{\mathbb{V}}
\newcommand{\bZ}{\mathbb{Z}}

\newcommand{\dummy}{\mathord{\textcolor[gray]{0.75}{\bullet}}}

\DeclareMathOperator*{\esssup}{ess\,sup}

\newcommand{\cA}{{\mathcal{A}}}
\newcommand{\cB}{{\mathcal{B}}}
\newcommand{\cD}{{\mathcal{D}}}
\newcommand{\cE}{{\mathcal{E}}}
\newcommand{\cF}{{\mathcal{F}}}
\newcommand{\cG}{{\mathcal{G}}}
\newcommand{\cH}{{\mathcal{H}}}
\newcommand{\cI}{{\mathcal{I}}}
\newcommand{\cK}{{\mathcal{K}}}
\newcommand{\cL}{{\mathcal{L}}}
\newcommand{\cM}{{\mathcal{M}}}
\newcommand{\cP}{{\mathcal{P}}}
\newcommand{\cR}{{\mathcal{R}}}
\newcommand{\cT}{{\mathcal{T}}}
\newcommand{\cV}{{\mathcal{V}}}
\newcommand{\cW}{{\mathcal{W}}}
\newcommand{\cX}{{\mathcal{X}}}

\newcommand{\sG}{{\mathscr{G}}}
\newcommand{\sL}{{\mathscr{L}}}
\newtheorem{theorem}{Theorem}[section]
\newtheorem{lemma}[theorem]{Lemma}
\newtheorem{proposition}[theorem]{Proposition}
\newtheorem{definition}[theorem]{Definition}

\newtheorem{remark}[theorem]{Remark}
\newtheorem{corollary}[theorem]{Corollary}
\theoremstyle{definition}
\newtheorem{example}[theorem]{Example}
\newtheorem{mainresult}{Result}
\theoremstyle{plain}

\AddToHook{env/theorem/begin}{\crefalias{section}{theorem}}
\AddToHook{env/lemma/begin}{\crefalias{section}{lemma}}
\AddToHook{env/proposition/begin}{\crefalias{section}{proposition}}
\AddToHook{env/definition/begin}{\crefalias{section}{definition}}
\AddToHook{env/assumption/begin}{\crefalias{section}{assumption}}
\AddToHook{env/remark/begin}{\crefalias{section}{remark}}
\AddToHook{env/corollary/begin}{\crefalias{section}{corollary}}
\AddToHook{env/example/begin}{\crefalias{section}{example}}

\AtBeginEnvironment{example}{\pushQED{\qed}}
\AtEndEnvironment{example}{\popQED}
\AtBeginEnvironment{mainresult}{\pushQED{\qed}}
\AtEndEnvironment{mainresult}{\popQED}
\newenvironment{namedparagraph}[1]
  {\paragraph{#1}\pushQED{\qed}}
  {\popQED\par}

\makeatletter
\let\kv@seccntformat\@seccntformat
\newcommand{\kv@appendixseccntformat}[1]{%
  \ifstrequal{#1}{section}{Appendix~}{}%
  \kv@seccntformat{#1}%
}
\renewcommand{\appendix}{%
  \par
  \setcounter{section}{0}%
  \setcounter{subsection}{0}%
  \gdef\thesection{\Alph{section}}%
  \gdef\thesubsection{\thesection.\arabic{subsection}}%
  \gdef\theHsection{appendix.\Alph{section}}%
  \gdef\theHsubsection{appendix.\Alph{section}.\arabic{subsection}}%
  \let\@seccntformat\kv@appendixseccntformat

  \crefname{section}{Appendix}{Appendices}%
  \Crefname{section}{Appendix}{Appendices}%
  \crefname{subsection}{Appendix}{Appendices}%
  \Crefname{subsection}{Appendix}{Appendices}%
}
\makeatother

\crefname{appendix}{Appendix}{Appendices}
\Crefname{appendix}{Appendix}{Appendices}

\begin{document}

\begin{frontmatter}
\title{Operator Learning of Microstructure-Dependent Internal-Variable Models for Inelastic Homogenization:\\ Approximation and Complexity}

\author{Lianghao Cao\corref{cor1}
\fnref{fn1}}\ead{lianghao@caltech.edu}

\author{George Stepaniants\fnref{fn2}}\ead{gs852@cam.ac.uk}

\author{Andrew M. Stuart
\fnref{fn1}}\ead{astuart@caltech.edu}

\cortext[cor1]{Corresponding author}

\fntext[fn1]{Department of Computing and Mathematical Sciences, California Institute of Technology, Pasadena CA, 91125, USA.}

\fntext[fn2]{Department of Applied Mathematics and Theoretical Physics, University of Cambridge, Wilberforce
Road, Cambridge CB3 0WA, United Kingdom.}

\begin{abstract}
    
We study the learning of homogenized constitutive laws for inelastic heterogeneous materials, aiming to replace costly concurrent unit-cell simulations with surrogates that generalize across loading histories and microstructures. We focus on the recurrent neural operator framework, motivated by the explicit form of homogenized Kelvin--Voigt linear viscoelasticity. In this framework, internal variables enable a Markovian approximation of loading-history dependence, and neural operators represent the microstructure-dependent constitutive response. For periodic Kelvin--Voigt viscoelasticity, we prove a universal approximation theorem showing that the model can approximate the homogenized constitutive law to any prescribed accuracy, uniformly in time and over a wide set of admissible strain trajectories and microstructural property fields, in two and three dimensions. This comprises a significant extension of the existing, somewhat limited, one-dimensional theory that uses explicit solutions for piecewise constant materials. Instead, in two and three dimensions, we adopt a more versatile approximation based on the spectral representation of the memory kernel that arises from homogenization. This also enables us to give an explicit bound on the internal-variable dimension, sufficient for approximation, exhibiting logarithmic growth in the inverse relative stress error and in material contrast. Numerical experiments on viscoelasticity and crystal plasticity complement the theory, illustrating the roles of internal-variable dimension and sample complexity in learning homogenized responses from unit-cell data. Theoretical results are formally verified in Lean 4.

\end{abstract}
\end{frontmatter}

\begin{figure}[H]
    \centering
    \includegraphics[width=0.9\textwidth]{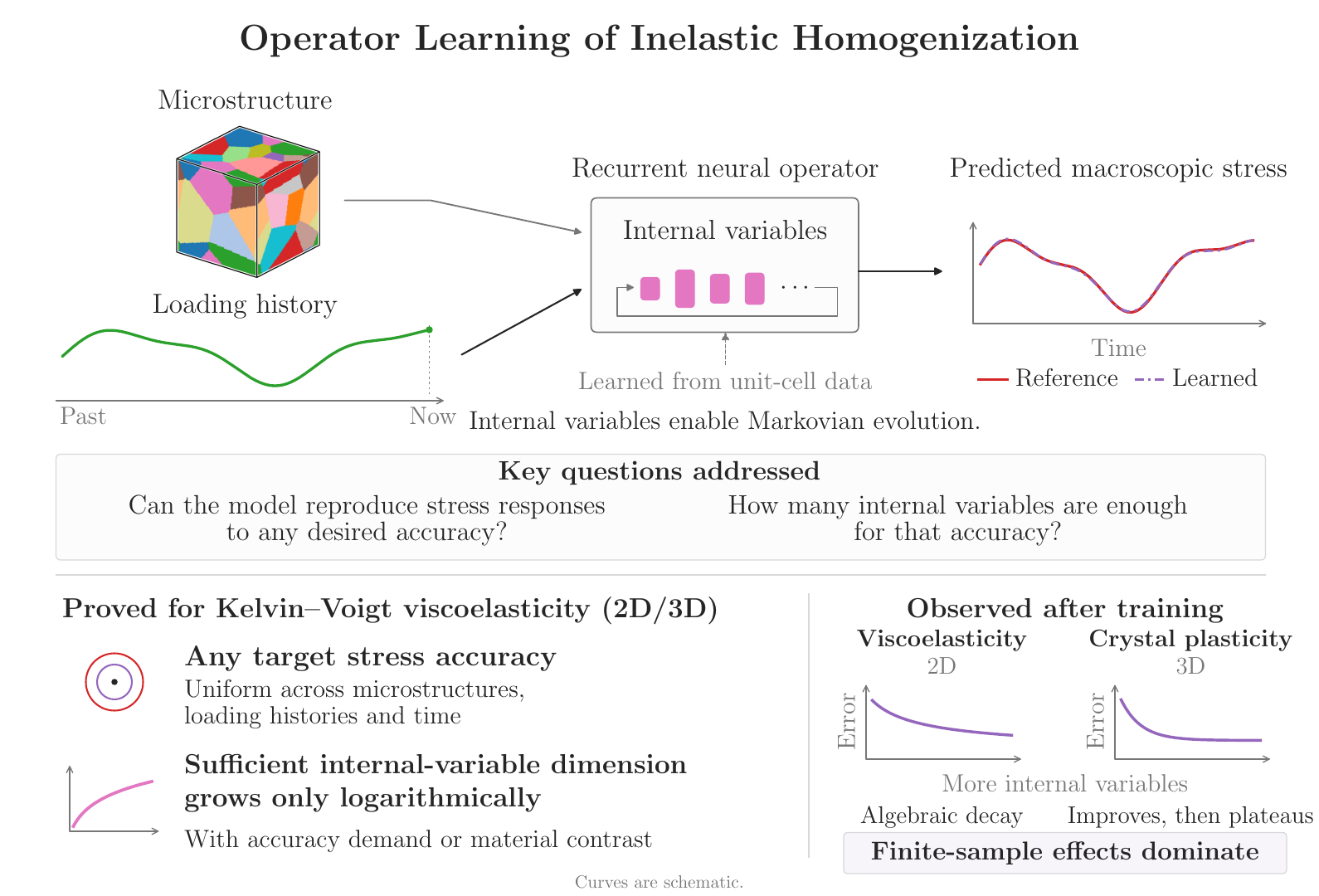}
    \caption{Graphical abstract.}
    \label{fig:graphical_abstract}
\end{figure}

\section{Introduction}

Macroscopic material behavior often reflects physical processes at much smaller length scales. Homogenization provides a systematic framework for deriving macroscopic constitutive models for multiscale materials by averaging the response of a unit cell that resolves the local microstructural behavior \cite{bensoussan1978asymptotic,sanchezpalencia1980nonhomogeneous,allaire1992homogenization,jikov1994homogenization, pavliotis2008multiscale}. For inelastic heterogeneous materials, the homogenized constitutive model $\cH$ depends both on the loading history and on the microstructural properties:
\begin{equation} \label{eq:one}
    \overline{P}(t) = \cH(\{\overline{F}(s)\}_{s\in[0,t]}; M),\quad t\geq 0,
\end{equation}
where $\overline{P}$ is the macroscopic first Piola--Kirchhoff stress tensor, $\overline{F}$ is the macroscopic deformation gradient, and $M$ denotes the microstructural property field on the unit cell. 

Typically, the form of $\cH$ is not known analytically, but it can be evaluated through unit-cell simulations. This makes macroscopic simulation via computational homogenization or FE$^2$ \cite{feyel2000fe2,miehe1999computational,geers2010multiscale} expensive, since the unit-cell problem must be repeatedly solved at macroscopic quadrature points and time steps. This has motivated recent efforts to learn $\cH$ from data \cite{ghavamian2019accelerating,logarzo2021smart,liu2022learning,liu2023learning,bhattacharya2023learning, karimi2024learning, raj2026neural}. Many existing approaches rely on standard machine-learning architectures applied to discretized loading histories and microstructural descriptions \cite{ghavamian2019accelerating,mozaffar2019deep,liu2019a,gorji2020on,rao2020threedimensional}, and therefore approximate discretized input-output relations rather than the underlying continuum micromechanics. Since loading histories and microstructural properties are functions of time and space, neural operators \cite{liu2022learning, kovachki2023neural, liu2023learning, bhattacharya2024learning, bhattacharya2026multiscale} offer a natural alternative: they learn mappings between function spaces, {are hence interpretable in the standard language of continuum mechanics} and can, in principle, generalize across discretizations. Following \cite{liu2023learning, bhattacharya2025learning}, we consider the neural operator model
\begin{equation*}
    \overline{P}(t)\approx \overline{P}_{\mathrm{NO}}(t) = \cH_{\mathrm{NO}}(\{\overline{F}(s)\}_{s\in[0,t]}; M),
\end{equation*}
defined through a microstructure-dependent internal-variable dynamical system:
\begin{subequations}
\label{eq:two}
\begin{align}
    \overline{P}_{\mathrm{NO}}(t) &= \cF_{\mathrm{NO}}(z(t), \xi(t); M),\\
    \dot{\xi}(t) &= \cG_{\mathrm{NO}}(z(t),\xi(t); M),\quad \xi(0)=0.
\end{align}
\end{subequations}
Here, $z(t)$ denotes the collection of relevant instantaneous macroscopic variables at time $t$, for example, $\overline{F}(t)$, $\dot{\overline{F}}(t)$, or both. In this model, $\xi\in\bR^{n_{\xi}}$ is a vector of internal variables that encodes the deformation history \cite{coleman1967thermodynamics,rice1971inelastic,halphen1975sur, he2022thermodynamically, jones2022neural, liu2023learning, bhattacharya2023learning, rosenkranz2024viscoelasticity, flaschel2025convex, jones2025attention, qiu2025bridging}, and {$\cF_{\rm NO}$ and $\cG_{\rm NO}$} are the stress and kinetic response, respectively. The internal-variable dimension $n_{\xi}$ and the response mappings {$\cF_{\rm NO},\cG_{\rm NO}$} are learned from data, and the response mappings are parameterized by neural operators with universal approximation properties \cite{lu2021learning, kovachki2023neural, lanthaler2025nonlocality}. In this formulation, each macroscopic material point stores only the current value of $\xi(t)\in\bR^{n_\xi}$ in place of the full loading history. This leads to an explicit, Markovian surrogate $\cH_{\rm NO}$ for the implicit, generally non-Markovian closure $\cH$. The dimension $n_{\xi}$ thus governs both the representational capacity of the model and the computational cost of the macroscopic simulation.

\begin{remark} \label{rem:one}
 {In the small strain setting the constitutive law \eqref{eq:one} is typically formulated in terms of the macroscopic infinitesimal strain $\overline{e} = \frac{1}{2}(\overline{F}+\overline{F}^{\top})-I$ and the Cauchy stress $\overline{\sigma}$ with $\overline{\sigma}\approx \overline{P}$. Then \eqref{eq:two} can be modified, replacing $\overline{F}$ with the macroscopic strain $\overline e$ and $\overline{P}$ with $\overline{\sigma}$.}
\end{remark}

This formulation raises the following questions for the neural operator model $\cH_{\rm NO}$ in relation to the homogenized model $\cH$:
\begin{enumerate}[label=(Q\arabic*), ref=Q\arabic*]
    \item\label{q:representation} (\textit{Approximation}) Can $\cH_{\rm NO}$ approximate $\cH$ to any prescribed accuracy?
    \item\label{q:model_complexity} (\textit{Internal-Variable Complexity}) How do the approximation properties of $\cH_{\rm NO}$ depend on the dimension of internal variables $n_{\xi}$?
\end{enumerate}
These questions are difficult to address in general because the unit-cell problem is typically not amenable to analysis. In this work, we answer these questions in two settings: Kelvin--Voigt viscoelasticity, for which we provide theory and numerical experiments, and finite-deformation crystal plasticity, for which we provide numerical experiments. 

For Kelvin--Voigt viscoelasticity, the fine-scale constitutive model is
\begin{equation} \label{eq:sigma}
    \sigma^{\epsilon}(x,t) = \bV(x/\epsilon)\partial_t e^{\epsilon}(x,t) + \bE(x/\epsilon)e^{\epsilon}(x,t),
\end{equation}
where $x\in\Omega$ denotes a macroscopic material point, $\epsilon>0$ is the microstructure length scale, $\sigma^{\epsilon}$ and $e^{\epsilon}$ are the $\epsilon$-scale stress and strain fields, and $\bV$ and $\bE$ are the microstructural viscosity and elasticity tensor fields that constitute $M$. {Kelvin--Voigt viscoelasticity results in a fading-memory model for the homogenized constitutive law in the limit $\epsilon\to 0$} \cite{coleman1960an,coleman1961foundations,tartar1990memory,ortiz2026linear}. Specifically, Francfort and Suquet \cite{francfort1986homogenization} show that at each material point, the macroscopic stress is governed by a constitutive law of the form
\begin{equation}
\label{eq:otf}
    \overline{\sigma}(t) = \bV^{\rm hom}\dot{\overline{e}}(t) + \bE^{\rm hom}\overline{e}(t) + \int_{0}^t \bK(t-s)\dot{\overline{e}}(s)\,\mathrm{d}s.
\end{equation}
Here, $\overline{\sigma}$ and $\overline{e}$ denote the macroscopic stress and strain, and $\dot{\overline{e}}\coloneqq\frac{\mathrm{d}}{\mathrm{d}t}\overline{e}$. The homogenized viscosity and elasticity tensors $\bV^{\rm hom}$ and $\bE^{\rm hom}$, as well as the memory kernel $\bK$, depend implicitly on the microstructural properties $M$. Thus, this homogenization result defines a constitutive law of the form
$\overline{\sigma}(t) = \cH(\{\overline{e}(s)\}_{s\in[0,t]}; M).$ Our interest is in approximating the memory
kernel to obtain an internal-variable Markovian model of the form
$\overline{\sigma}_{\rm NO}(t) = \cH_{\mathrm{NO}}(\{\overline{e}(s)\}_{s\in[0,t]}; M).$

\subsection{Main Results and Paper Overview}\label{subsec:overview_results}

In this subsection, we summarize our main results and outline the rest of the paper. We formally verified the theoretical results in Lean 4 \cite{demoura2021lean4} using mathlib \cite{mathlib2020}; details appear below. The formalization and verification instructions are available in the accompanying GitHub repository\footnotemark.

\footnotetext{\url{https://github.com/sgstepaniants/inelastic_hom_learn}}

\begin{mainresult}[Homogenized Model]\label{res:homogenized_model}
Periodic homogenization of Kelvin--Voigt viscoelasticity {gives rise to} a fading-memory constitutive law whose instantaneous response and memory kernel are determined by bounded linear operators arising from static unit-cell problems. The memory kernel admits a spectral representation whose relaxation-rate spectrum has bounded support and bounded total-variation norm.
\end{mainresult}

The full statement appears in \cref{thm:viscoelastic_homogenization,prop:spectral_rep}, and both are formally verified. 
Building on the framework of \cite{francfort1986homogenization}, this model form and the 
spectral representation of the kernel provide the basis for the approximation theory that follows.

\begin{mainresult}[Universal Approximation]\label{res:ua}
The parameters of the neural operator model $\cH_{\mathrm{NO}}$ may be chosen to approximate the true homogenized constitutive law $\cH$ to any prescribed accuracy, uniformly in time, over microstructures whose viscosity and elasticity tensor fields are uniformly bounded, coercive, and of bounded spatial variation, and over strain trajectories with uniformly bounded magnitudes and rates.

\end{mainresult}

The full statement appears in \cref{thm:ua} and is formally verified. This theorem addresses (\ref{q:representation}) and,
although similar in form to the one-dimensional approximation results of \cite{bhattacharya2025learning},
requires substantial new ideas in two and three spatial dimensions. The previous one-dimensional analysis relies on the fact that the Prony series approximation of the memory kernel {\cite{baumgaertel1989determination,tschoegl1989phenomenological,baumgaertel1992interrelation,park1999methods,beylkin2005approximation,ortiz2026linear}} is exact for piecewise-constant microstructures \cite{bhattacharya2023learning}. In two and three dimensions, we instead directly estimate the approximation error without resorting to a specific microstructure. The proof constructs a Prony-series internal-variable model whose stress response approximates the homogenized constitutive law uniformly over admissible strain trajectories and microstructures. Its continuous dependence on microstructure enables approximation by neural operators. Assuming that the neural operator model emulates this intermediate model, a quantitative answer to (\ref{q:model_complexity}) is obtained by the approximation rate of the Prony series:

\begin{mainresult}[Internal-Variable Complexity]\label{res:complexity}
{The sufficient internal-variable dimension grows logarithmically with the inverse relative stress error and with each material contrast when the remaining parameters are held constant.}
\end{mainresult}

The full statement appears in \Cref{thm:sufficient_dimension} and is formally verified.
This scaling suggests that the required internal-variable dimension remains computationally manageable for uniformly accurate homogenization across microstructural instances with large property variations or a small tolerance for approximation error. It also supports the practical use of the internal-variable memory representation in macroscopic simulations.

We also present empirical results on training and testing the neural-operator model on the viscoelasticity and crystal plasticity datasets. We report the supporting numerical study in \Cref{sec:numerical_examples}.

\begin{mainresult}[Empirical Complexity]\label{res:empirical} For viscoelasticity, the accuracy of the trained model improves algebraically when the internal-variable dimension increases. For crystal plasticity, the trained model's accuracy improves and then plateaus as the internal-variable dimension increases.
\end{mainresult}

The remainder of the paper is organized as follows. \Cref{sec:formal_homogenization} introduces the periodic homogenization setting for Kelvin--Voigt viscoelasticity and defines the neural operator constitutive model. \Cref{sec:homogenization_theory} presents the homogenization theory, identifies the operator form of the homogenized model, and provides the spectral representation of the memory kernel; together, these establish Result~\ref{res:homogenized_model}. \Cref{sec:universal_approximation} states the universal approximation theorem and the internal-variable complexity bound, corresponding to Results~\ref{res:ua} and~\ref{res:complexity}. \Cref{sec:numerical_examples} reports numerical experiments on viscoelasticity and crystal plasticity (Result~\ref{res:empirical}), and \cref{sec:conclusion} summarizes the findings.

\begin{remark} \label{rem:2}
Result~\ref{res:homogenized_model} underpins the approximation and complexity results. It is therefore valuable to give intuition about why this result holds. If the divergence of the stress in \eqref{eq:sigma} is balanced with an $\epsilon$-independent body force $f$ then this results in a dynamical system of the form
\begin{equation*}
    \partial_t\mathcal L_{\bV}^{\epsilon}u+\mathcal L_{\bE}^{\epsilon}u=f.
\end{equation*}
Here, $u=u^\epsilon$ is the displacement field. Based on standard assumptions concerning the material properties, the elliptic operators $\mathcal L_{\bV}^{\epsilon}$ and $\mathcal L_{\bE}^{\epsilon}$ are associated with coercive and bounded quadratic forms on the displacement space. The relaxation rates of the problem are governed by the bounded operator $(\mathcal L_{\bV}^{\epsilon})^{-1}\mathcal L_{\bE}^{\epsilon}$. If the coercivity and boundedness hold uniformly in $\epsilon$, then this leads to upper and lower bounds, uniformly in $\epsilon$, for the relaxation rates of the multiscale problem. Therefore, it is unsurprising that the relaxation rates appearing in $\bK$ should satisfy uniform positive lower and upper bounds. However, additional effects come into play in the memory and, in particular, memory vanishes if $\bV=\bE$. Hence the operator $(\mathcal L_{\bV}^{\epsilon})^{-1}\mathcal L_{\bE}^{\epsilon}$ does not by itself determine the memory kernel; rather the analysis also quantifies the discrepancy between $\bV$ and $\bE$.
\end{remark}

\section{Kelvin--Voigt Viscoelasticity: Homogenization and Neural Operator Model} \label{sec:formal_homogenization}

In this section, we present {the unit-cell formulation of the homogenized response of} a Kelvin--Voigt viscoelastic heterogeneous material and introduce a neural-operator constitutive model designed to learn {this response}. In \cref{subsec:formal_homogenization_settings}, we define the geometric and tensorial setting for the heterogeneous material and periodic unit cell. In \cref{subsec:unit_cell_formulation}, we pose the $\epsilon$-scale problem and {define the constitutive response through a periodic unit-cell problem, before coupling it to the macroscopic balance equations.} In \cref{subsec:viscoelastic_neural_operator_model}, we present the neural operator constitutive model.

\subsection{Settings}\label{subsec:formal_homogenization_settings}

\begin{namedparagraph}{Spatial and Temporal Domains}
Let $\Omega\subset\bR^{d}$ with $d\in\{1,2,3\}$ be a bounded open material domain and $\bT^d\coloneqq\bR^{d}/\bZ^d$ the $d$-dimensional unit torus with normalized Lebesgue measure. Throughout, we identify 1-periodic functions with functions on $\bT^d$. For the homogenization problem, choose any terminal time $T>0$ and let $t\in[0,T]$ denote time.
\end{namedparagraph}

\begin{namedparagraph}{Second and Fourth Order Tensors}
Let $\cV^2$ denote the set of symmetric 2nd-order tensors, equipped with the Frobenius inner product $(e_1, e_2)_F\coloneqq{\operatorname{tr}}(e_1^{\top}e_2)$ and the induced norm $\|e\|_{F}$. Let $\cV^4$ denote the set of 4th-order tensors with major and minor symmetries. For $\bA \in \cV^4$ and $e \in \cV^2$, the action $\bA e \in \cV^2$ is defined by $(\bA e)_{ij} = \sum_{k,l} \bA_{ijkl} e_{kl}$, and the Frobenius inner product on $\cV^4$ is $(\bA,\bB)_F = \sum_{i,j,k,l} \bA_{ijkl}\bB_{ijkl}$. The positive-definite subset $\cV^4_{+}\subset\cV^4$ consists of tensors $\bA$ with $(e,\bA e)_F\geq \alpha \|e\|_{F}^2$ for some $\alpha>0$ and all $e\in\cV^2$. Throughout, we assume tensor fields with values in $\cV^4_{+}$ are uniformly coercive, so $\alpha$ can be chosen independently of position.
\end{namedparagraph}

\subsection{\texorpdfstring{{Microscopic, Unit-Cell, and Macroscopic Problems}}{Microscopic, Unit-Cell, and Macroscopic Problems}}\label{subsec:unit_cell_formulation}

\begin{namedparagraph}{The $\epsilon$-Scale Problem}
Let $\epsilon>0$ be the microstructure length scale relative to the material domain $\Omega$. Let $\bV, \bE:\bT^d\to\cV^4_{+}$ be 1-periodic viscosity and elasticity tensor fields that characterize the microstructure. With $\bV^{\epsilon}(x) = \bV(x/\epsilon)$, $\bE^{\epsilon}(x) = \bE(x/\epsilon)$, and $\mathrm{sym}(\dummy)$ the symmetric part, the $\epsilon$-scale problem is:
\begin{subequations}\label{eqns:two_scale_problem}
\begin{align}
    e^{\epsilon}(x,t) &=\mathrm{sym}(\nabla_x u^{\epsilon}(x,t)),&& x \in \Omega, \; t \in [0,T], \label{eq:two_scale_strain}\\
    \sigma^{\epsilon}(x, t) &= \bV^{\epsilon}(x) \partial_t{e}^{\epsilon}(x, t) + \bE^{\epsilon}(x)e^{\epsilon}(x,t), && x \in \Omega, \; t \in [0,T], \label{eq:two_scale_stress}\\
    -\nabla_x \cdot \sigma^{\epsilon}(x,t) &=  f(x,t) ,&& x \in \Omega, \; t \in [0,T], \label{eq:two_scale_balance}\\
    u^{\epsilon}(x, 0) & = u^{\epsilon}_0(x), && x\in\Omega, \\
    u^{\epsilon}(x, t) & = 0, && x\in\partial\Omega.
\end{align}
\end{subequations}
Here, $u^{\epsilon}:\Omega\times[0,T]\to\bR^d$ is the displacement, $e^{\epsilon}:\Omega\times[0,T]\to\cV^2$ is the strain, $\sigma^{\epsilon}:\Omega\times[0,T]\to\cV^2$ is the stress, $u^{\epsilon}_0:\Omega\to\bR^d$ is the initial condition, and $f:\Omega\times[0,T]\to\bR^d$ is an external force.
\end{namedparagraph}

\begin{namedparagraph}{The Unit-Cell Problem}

{At a fixed macroscopic material point, we prescribe a strain history $\overline e:[0,T]\to\cV^2$ and consider the periodic unit cell with the same material fields $\bV,\bE$ that enter the $\epsilon$-scale problem \cref{eqns:two_scale_problem}. A periodic displacement correction allows the local strain to vary while preserving its prescribed spatial average $\overline e(t)$. Denote the displacement corrector by $\chi$ and the cell strain and stress fields by $e^{(0)}$ and $\sigma^{(0)}$, respectively. They satisfy the time-domain unit-cell problem}
\begin{subequations}\label{eqns:cell_problem}
\begin{align}
    \eta(y,t) & = \mathrm{sym}\left(\nabla_y \chi (y,t)\right), && y \in \bT^d,\; t \in [0,T], \label{eq:strain_corrector}\\
    e^{(0)}(y,t) &= \overline{e}(t) + \eta(y,t), && y \in \bT^d,\; t \in [0,T], \\
    \sigma^{(0)}(y,t) & = \bV(y) \partial_t e^{(0)}(y,t) + \bE(y) e^{(0)}(y,t), && y \in \bT^d,\; t \in [0,T], \label{eq:cell_stress}\\
    - \nabla_y \cdot \sigma^{(0)}(y,t) &= 0, && y \in \bT^d,\; t \in [0,T], \\
    \chi(y,0) & = \chi_{0}(y), && y\in\bT^d, \\
    \int_{\bT^d} \chi(y,t)\,\mathrm{d} y &= 0, && t\in [0,T],
\end{align}
\end{subequations}
Here, $\chi_0$ is chosen as the solution of the static elastic unit-cell problem corresponding to $\overline{e}(0)$. We call $\eta$ the strain corrector.
\end{namedparagraph}

\begin{namedparagraph}{Homogenized Model and the Macroscopic Problem}
{The averaged unit-cell stress defines the constitutive map $\cH$ by}
\begin{equation}\label{eq:cell_constitutive_map}
    {\overline\sigma(t)\coloneqq\int_{\bT^d}\sigma^{(0)}(y,t)\,\mathrm dy
    =\cH\bigl(\{\overline e(s)\}_{s\in[0,t]};\bV,\bE\bigr).}
\end{equation}
{This is the small-strain analog of \eqref{eq:one}; see Remark \ref{rem:one}.}
{Solving \cref{eqns:cell_problem} for sampled strain histories and material fields produces the stress trajectories used to train the neural operator model in \cref{subsec:2d_viscoelasticity}. The homogenization theorem in \cref{sec:homogenization_theory} identifies the constitutive law obtained as $\epsilon\to0$ in \cref{eqns:two_scale_problem}. With the initialization specified above, \cref{app:unit_cell_homogenization} shows that the averaged response \cref{eq:cell_constitutive_map} coincides with that law.}

{Let $u^{(0)}:\Omega\times[0,T]\to\bR^d$ denote the homogenized displacement. At each macroscopic point $x$, its strain history drives the unit-cell response in \cref{eq:cell_constitutive_map}. Coupling this constitutive response to macroscopic compatibility and force balance gives the following problem:}
\begin{subequations}\label{eqns:homogenized_problem}
\begin{align}
    \overline{e}(x,t) &= \mathrm{sym}\left(\nabla_x u^{(0)}(x,t)\right), && x \in \Omega,\; t \in [0,T], \\
    \overline{\sigma}(x,t) &= \cH\big(\{\overline{e}(x,s)\}_{s\in[0,t]};\bV,\bE\big), && x \in \Omega,\; t \in [0,T], \\
    -\nabla_x \cdot \overline{\sigma}(x,t) &= f(x,t), && x \in \Omega,\; t \in [0,T], \\
    u^{(0)}(x,0) &= u_0^{(0)}(x), && x\in\Omega, \\
    u^{(0)}(x,t) &= 0, && x\in\partial\Omega,\; t\in [0,T],
\end{align}
\end{subequations}
{The unit-cell fields depend parametrically on $x$ through the strain history $\overline e(x,\dummy)$; this parameter is suppressed in \cref{eqns:cell_problem,eq:cell_constitutive_map}, {because $\cH$ acts pointwise} in space. The homogenized initial condition $u_0^{(0)}$ is specified in \cref{thm:viscoelastic_homogenization}.}
\end{namedparagraph}

\subsection{Viscoelastic Neural Operator Constitutive Model}\label{subsec:viscoelastic_neural_operator_model}

Since $\cH$ {is defined pointwise} in space, we suppress $x$ and define the neural operator model
$\cH_{{\mathrm{NO}}}$ that approximates $\cH$ in \cref{eq:cell_constitutive_map} at a given material point. Specializing \eqref{eq:two} to the small strain setting, noting Remark \ref{rem:one}, we consider
\begin{equation}\label{def:learned_model}
    \overline{\sigma}_{\mathrm{NO}}(t) = \cH_{{\mathrm{NO}}}(\{\overline{e}(s)\}_{s\in[0,t]}; \bV, \bE),
\end{equation}
defined through the dynamical system
\begin{subequations}\label{eq:learned_model_full}
\begin{align}
    \overline{\sigma}_{\mathrm{NO}}(t) &= \cF_{\mathrm{NO}}(\overline{e}(t), \dot{\overline{e}}(t), \xi(t); \bV, \bE), && t\in [0,T],\label{eq:learned_model_stress}\\
    \dot{\xi}(t) &= \cG_{\mathrm{NO}}(\dot{\overline{e}}(t),\xi(t); \bV, \bE), && t\in [0,T],\label{eq:learned_model_kinetics}\\
    \xi(0)&=0,
\end{align} 
\end{subequations}
where $\xi:[0,T]\to\bR^{n_\xi}$ is an internal variable, and only the strain rate enters the kinetic response $\cG_{\mathrm{NO}}$. The choice of inputs to $\cF_{\mathrm{NO}}$ and $\cG_{\mathrm{NO}}$ made here reflects the structure of the Kelvin--Voigt homogenized model identified in \cref{sec:homogenization_theory}. For other materials, different choices may be appropriate: finite strain may be relevant, the kinetic response \cref{eq:learned_model_kinetics} may also depend on strain, or both \cref{eq:learned_model_stress,eq:learned_model_kinetics} may depend only on strain and not strain rate \cite{liu2023learning}.

We note that the response mappings defining \eqref{eq:learned_model_full} depend on \emph{functions} summarizing material properties. To handle this, we use neural operator concepts. As well as operating on functions, the chosen neural operator must also allow Euclidean vector-to-function lifting and function-to-vector projection layers to accommodate the dependence on strain, strain rate, and internal variables. In this work, we use the Fourier neural mapping (FNM) \cite{huang2025operator} to achieve this.

\begin{example}[Fourier Neural Mapping]
Let $\mathsf{vec}_2:\cV^2\to\bR^{d_{\cV^2}}$ and $\mathsf{vec}_4:\cV^4\to\bR^{d_{\cV^4}}$ denote the Mandel vectorizations \cite{helnwein2001some}, which are isometric isomorphisms, and denote $d_{\cV^2}={d(d+1)}/{2}$ and $d_{\cV^4}={d_{\cV^2}(d_{\cV^2}+1)}/{2}$. We parameterize $\cF_{\rm NO}$ and $\cG_{\rm NO}$ by two FNM instances $\Phi_{\cF}$ and $\Phi_{\cG}$:
\begin{subequations}\label{eq:F_and_G}
    \begin{align}
    \cF_{\rm NO} &= \mathsf{vec}_2^{-1}\circ\Phi_{\cF}\left(\bigl(\mathsf{vec}_2\overline{e},\mathsf{vec}_2\dot{\overline{e}},\xi\bigr), \left(\mathsf{vec}_4\circ \bV,\mathsf{vec}_4\circ\bE\right)\right),\label{eq:stress_response}\\
    \cG_{\rm NO} &= \Phi_{\cG}\left(\bigl(\mathsf{vec}_2\dot{\overline{e}},\xi\bigr), \left(\mathsf{vec}_4\circ \bV,\mathsf{vec}_4\circ\bE\right)\right).\label{eq:kinetic_response}
\end{align}
\end{subequations}
Each FNM $\Phi$ takes a vector input $v_{\rm in}$ and a function input $M=(\mathsf{vec}_4\circ \bV,\mathsf{vec}_4\circ\bE):\bT^d\to\bR^{2d_{\cV^4}}$ and returns $v_{\rm out}\in\bR^{d_{\mathrm{out}}}$. It has three components:

\begin{enumerate}[label=(\roman*)]
    \item \textit{Lifting}: A trainable pointwise map $S:\bR^{d_{\mathrm{in}}+2d_{\cV^4}}\to\bR^{d_{h^{(0)}}}$ defines $h^{(0)}(y)=S(v_{\rm in},M(y))$.
    \item \textit{Fourier layers}: For $\ell=1,\dots,L$, the layer $\sL_\ell$ maps $h^{(\ell-1)}$ to
    \begin{equation*}
    h^{(\ell)}(y)=\sigma_{\rm act}^{(\ell)}\Bigl(W_{\ell}h^{(\ell-1)}(y)+\sum_{k\in\Lambda_{\ell}} P_{\ell}^{(k)}\widehat{h^{(\ell-1)}}(k)\,e^{2\pi\mathrm{i} k^{\top}y}+b_{\ell}\Bigr),
    \end{equation*}
    where $\sigma_{\rm act}^{(\ell)}$ is the activation function at layer $\ell$, $\widehat{(\dummy)}$ denotes the Fourier transform, $W_{\ell}$ and $b_{\ell}$ are trainable weights and biases applied pointwise, $P_{\ell}^{(k)}$ are trainable Fourier-space weights satisfying conjugate symmetry, and $\Lambda_{\ell}=\{k\in\bZ^d : |k_j|\le m_{\ell,j}\}$ is the set of the retained Fourier modes.
    \item \textit{Projection}: A Fourier-space projection $\sG h^{(L)}=\sum_{k\in\Lambda_{\rm proj}} P_{\rm proj}^{(k)}\widehat{h^{(L)}}(k)\in\bR^{d_{\mathrm{proj}}}$ is followed by a trainable map $Q:\bR^{d_{\mathrm{proj}}}\to\bR^{d_{\mathrm{out}}}$. The weights $P_{\rm proj}^{(k)}$ satisfy conjugate symmetry.
\end{enumerate}
The resulting map is $\Phi(v_{\rm in},M)=(Q\circ\sG\circ\sL_L\circ\cdots\circ\sL_1)(h^{(0)})$. 
\end{example}

\begin{figure}[H]
    \centering
    \includegraphics[width=0.90\textwidth]{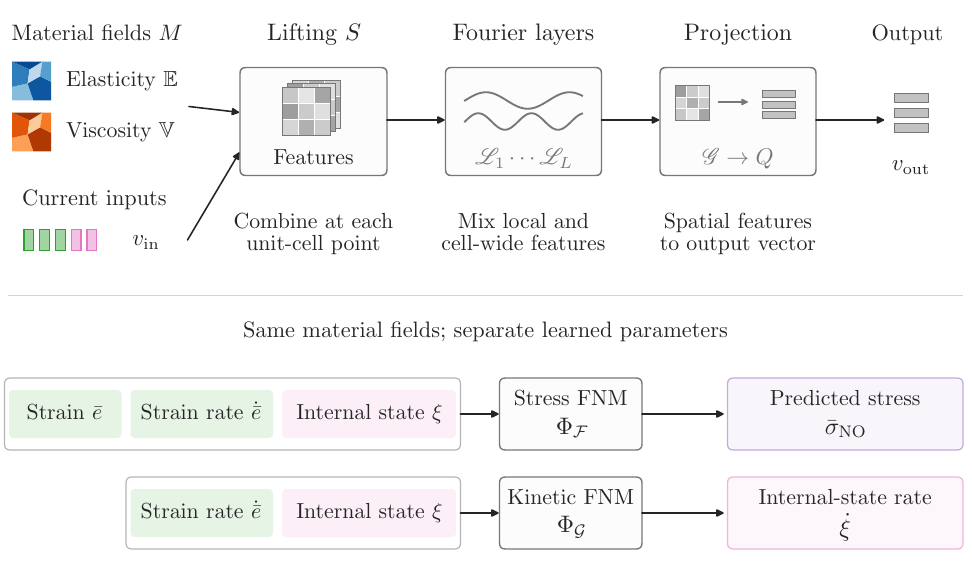}
    \caption{Fourier neural mappings for the viscoelastic constitutive model. At each unit-cell spatial point and time, we combine the relevant macroscopic states and material fields to form spatial features, process them with Fourier layers, and project them to an output. Two separately parameterized instances predict macroscopic stress and the internal-variable rate in \eqref{eq:F_and_G}. Tensor vectorizations are omitted in the schematic.}
    \label{fig:fnm_schematic}
\end{figure}

\section{Homogenization Theory and Memory Kernel}\label{sec:homogenization_theory}

We revisit the viscoelastic homogenization theory of Francfort and Suquet \cite{francfort1986homogenization} and identify the homogenized model $\cH$ with fading memory in a form suitable for analysis. {To identify the homogenization limit of the $\epsilon$-scale problem and justify the constitutive response introduced in \cref{subsec:unit_cell_formulation}, we use the Laplace transform to reduce the evolution equation to a family of static elliptic problems to which classical periodic homogenization applies.} In \cref{subsec:homogenization_theory_settings}, we introduce the functional setting. In \cref{subsec:unit_cell_operators}, we define the bounded linear operators that govern $\cH$. In \cref{subsec:viscoelastic_homogenization_fading_memory}, we state the homogenization theorem and sketch its proof. In \cref{subsec:spectral_rep}, we present the spectral representation of the memory kernel that arises from homogenization.

\subsection{Functional Setting}\label{subsec:homogenization_theory_settings}

We use standard notation for Lebesgue spaces $L^p$, Sobolev spaces $W^{k,p}$ (with $H^k\coloneqq W^{k,2}$). Periodic function spaces on $\bT^d$ carry the subscript $\mathrm{per}$. The space $H_0^1(\Omega;\bR^d)$ denotes functions with homogeneous Dirichlet boundary conditions.

Let $X$ and $Y$ be Banach spaces. We denote Bochner and Sobolev spaces for time-evolving functions as $L^p(0,T;X)$ and $W^{k,p}(0,T;X)$ (with $H^k\coloneqq W^{k,2}$). We write $\cL(X,Y)$ for the space of bounded linear operators from $X$ to $Y$.

For $\bA\in L^{\infty}_{\mathrm{per}}(\bT^d;\cV^4_{+})$, we define $\bA$-weighted inner products $(e_1, e_2)_{\bA} \coloneqq (\bA e_1, e_2)_{L^2_{\rm per}}$ on $L^{2}_{\mathrm{per}}(\bT^d;\cV^2)$, with induced norms equivalent to the base norm. For a macroscopic tensor field $\bB\in L^{\infty}(\Omega;\cV_{+}^4)$, we define $(u, v)_{\bB,\, \mathrm{sym}} \coloneqq (\bB\mathrm{sym}(\nabla_x u), \mathrm{sym}(\nabla_x v))_{L^2}$ on $H_0^1(\Omega;\bR^d)$, which is an equivalent inner product by Korn's inequality. Recall that the Frobenius inner product is used for second- and fourth-order tensors.

The displacement corrector in the unit-cell problem is sought in $\mathring{H}^1_{\mathrm{per}}(\bT^d; \bR^d)$, the subspace of $H^1_{\mathrm{per}}(\bT^d; \bR^d)$ with zero mean. The strain corrector
space is
\begin{equation}\label{eq:space_strain_corrector}
    \cX \coloneqq \left\{ \mathrm{sym}(\nabla_y \chi) \in L^2_{\mathrm{per}}(\bT^d; \cV^2) \mid \chi \in \mathring{H}^1_{\mathrm{per}}(\bT^d; \bR^d) \right\},
\end{equation}
a closed subspace of $L^2_{\mathrm{per}}(\bT^d; \cV^2)$ by Korn's inequality \cite{ciarlet2021mathematical}. The space $\cX$ contains no nonzero spatially constant tensors.

\subsection{Continuous Linear Operators for the Unit-Cell Problem}\label{subsec:unit_cell_operators}

We define the operators governing the form of $\cH$, assuming $\bV,\bE\in L^{\infty}_{\mathrm{per}}(\bT^d;\cV^4_{+})$. {Whenever a macroscopic tensor $\overline e\in\cV^2$ is an input to a function-space operator, we identify it with the constant field $y\mapsto\overline e$ on $\bT^d$.} \Cref{tab:operators} summarizes these operators.

\begin{namedparagraph}{Strain Corrector Projectors}
The orthogonal projectors $\cP_{\bV}, \cP_{\bE}\in \mathcal{L}(L^{2}_{\mathrm{per}}(\bT^d;\cV^2), \cX)$ are defined by
\begin{equation}\label{eqns:operator_P}
        \left(\cP_{\bV}e,\eta\right)_{\bV} = (e, \eta)_{\bV},\quad \left(\cP_{\bE}e, \eta\right)_{\bE} = (e, \eta)_{\bE}\quad \forall \eta\in \cX.
\end{equation}
For a static macroscopic strain $\overline{e}\in\cV^2$, $\cP_{\bV}\overline{e}$ and $\cP_{\bE}\overline{e}$ are the negative strain correctors solving the purely viscous or elastic static unit-cell problems. For example, the action of the viscous projector $\cP_{\bV}$ on $\overline e$ is defined by \cref{eq:viscous_cell_projector} where $\chi\in\mathring H^1_{\mathrm{per}}(\bT^d;\bR^d)$ is the weak solution of \cref{eq:viscous_cell_balance,eq:viscous_cell_mean}:
\begin{subequations}\label{eqns:viscous_cell}
\begin{align}
[\cP_{\bV}\overline{e}](y) &= -{\rm sym}(\nabla_y \chi(y)),&& y \in \bT^d,\label{eq:viscous_cell_projector}\\
    \nabla_y \cdot \bV(y){\rm sym}(\nabla_y \chi(y)) &= -\nabla_y \cdot \bV(y) \overline{e}, && y \in \bT^d, \label{eq:viscous_cell_balance}\\
    \int_{\bT^d} \chi(y)\,\mathrm{d} y &= 0.\label{eq:viscous_cell_mean}
\end{align}
\end{subequations}
where $\chi:\bT^d\to\bR^{d}$ is the displacement corrector given the macroscopic strain $\overline{e}\in\cV^2$. {Replacing $\bV$ by $\bE$ gives the corresponding elastic cell problem. The tensors $\bV^{\rm hom}$ and $\bE^{\rm hom}$ are the effective tensors obtained from these separate static problems by standard periodic elliptic homogenization.}
\end{namedparagraph}

\begin{namedparagraph}{Strain Localization and Mismatch Operators}
The complementary projectors $\cP^{\perp}_{\bV} = \mathrm{Id}_{L^2_{\rm per}} - \cP_{\bV}$ and $\cP^{\perp}_{\bE} = \mathrm{Id}_{L^2_{\rm per}} - \cP_{\bE}$ return the corrected microscopic strain fields, where $\mathrm{Id}$ denotes the identity on $L^2_{\rm per}(\bT^d;\cV^2)$. The mismatch operator $\cD_{\bV,\bE} \coloneqq \cP_{\bV} - \cP_{\bE}$ quantifies the incompatibility between viscous flow and elastic equilibrium and drives memory effects in the homogenized model.
\end{namedparagraph}

\begin{namedparagraph}{Relaxation Operator}
The operator $\cR_{\bV,\bE}\in\mathcal{L}(\cX)$ defined by
\begin{equation}\label{eq:operator_L}
    (\cR_{\bV,\bE} \eta_1, \eta_2)_{\bV} = (\eta_1, \eta_2)_{\bE}\quad\forall \eta_1, \eta_2\in\cX
\end{equation}
is the Riesz representation map associated with $(\dummy,\dummy)_{\bE}$ relative to $(\dummy,\dummy)_{\bV}$. {Equivalently, $\cR_{\bV,\bE}=\left.\cP_{\bV}W_{\bV,\bE}\right|_{\cX}$, where $W_{\bV,\bE}$ denotes pointwise multiplication by $\bV^{-1}\bE$.} It governs relaxation dynamics by balancing elastic restoring forces against viscous resistance.

\begin{samepage}
Given $w:\bT^d\to\bR^d$ in $H^1_{\mathrm{per}}(\bT^d;\bR^d)$ with $\int_{\bT^d} w(y)\,\mathrm{d} y = 0$, it is defined by {\cref{eq:relaxation_cell_operator}} where $v:\bT^d\to\bR^d$ is the weak solution of \cref{eq:relaxation_cell_balance,eq:relaxation_cell_mean}:\par
\end{samepage}
\begin{subequations}\label{eqns:relaxation_cell}
\begin{align}
[\cR_{\bV,\bE}{\rm sym}(\nabla_y w)](y) &= {\rm sym}(\nabla_y v(y)),&& y \in \bT^d,\label{eq:relaxation_cell_operator}\\
    \nabla_y \cdot \bV(y){\rm sym}(\nabla_y v(y)) &= \nabla_y \cdot \bE(y) {\rm sym}(\nabla_y w(y)), && y \in \bT^d, \label{eq:relaxation_cell_balance}\\
    \int_{\bT^d} v(y)\,\mathrm{d} y &= 0.\label{eq:relaxation_cell_mean}
\end{align}
\end{subequations}
The family $\{\exp(-t\cR_{\bV,\bE})\}_{t\ge0}$ is a strongly continuous semigroup on $(\cX,(\dummy,\dummy)_{\bV})$, generated by $-\cR_{\bV,\bE}$. Note that the mapping from $w$ to $v$ given by \cref{eq:relaxation_cell_balance,eq:relaxation_cell_mean} is defined by an operator $(\mathcal{L}_{\bV})^{-1}\mathcal{L}_{\bE}$, with $\mathcal{L}_{\bV}$ and $\mathcal{L}_{\bE}$ elliptic operators; this is the cell problem analog of the family of operators $(\mathcal{L}_{\bV}^{\epsilon})^{-1}\mathcal{L}_{\bE}^{\epsilon}$ discussed in Remark \ref{rem:2} and the same remarks concerning boundedness apply here.
\end{namedparagraph}

\begin{table}[htb]
    \centering
    \scalebox{0.95}{
    \renewcommand*{\arraystretch}{1.5}
    \begin{tabular}{
    | >{\centering\arraybackslash} m{0.18\linewidth}  
    | >{\centering\arraybackslash} m{0.15\linewidth}
    | >{\centering\arraybackslash} m{0.15\linewidth}
    |  >{\centering\arraybackslash} m{0.22\linewidth}
    |  >{\centering\arraybackslash} m{0.18\linewidth}|}\hline
        Name & Input & Output & \makecell{Physical\\Interpretation} & \makecell{Mathematical\\Interpretation} \\\hline\hline
        \makecell{Strain Corrector\\ Projectors\\$\cP_{\bV}$ and $\cP_{\bE}$} & $ L^{2}_{\mathrm{per}}(\bT^d;\cV^2)$ & $\cX$ & Identify purely viscous or elastic microscopic corrector (negative) & Orthogonal Projectors \\\hline

        \makecell{Strain \\ Localization\\$\cP^{\perp}_{\bV}$ and $\cP^{\perp}_{\bE}$} & $ L^{2}_{\mathrm{per}}(\bT^d;\cV^2)$ & $ L^{2}_{\mathrm{per}}(\bT^d;\cV^2)$ & Determine purely viscous or elastic microscopic strain & Complementary Projectors \\\hline
        
        \makecell{Strain Corrector\\Mismatch\\ $\cD_{\bV,\bE}$} & $L^{2}_{\mathrm{per}}(\bT^d;\cV^2)$ & $\cX$ & Difference between purely viscous and elastic microscopic correctors & Difference of Projectors\\\hline
     \makecell{Relaxation\\$\cR_{\bV,\bE}$} & $\cX$ & $\cX$ &  Dictates the relaxation timescales & Riesz Representation Map\\\hline
    \end{tabular}
    \renewcommand*{\arraystretch}{1.0}
    }
    \caption{Continuous linear operators governing the homogenized constitutive model.}
    \label{tab:operators}
\end{table}

\subsection{Viscoelastic Homogenization with Fading Memory}\label{subsec:viscoelastic_homogenization_fading_memory}

We identify the homogenized constitutive model for Kelvin--Voigt viscoelasticity. We first define the weak $\epsilon$-scale problem and admissible initial conditions, then state the homogenization theorem adapted from \cite{francfort1986homogenization} and sketch its proof.

\begin{namedparagraph}{Weak $\epsilon$-Scale Problem} Given $\bV,\bE\in L^{\infty}_{\mathrm{per}}(\mathbb{T}^d;\cV^4_{+})$, $f \in L^2(0,T; L^2(\Omega; \bR^d))$, and initial condition $u_0^\epsilon$, find $u^{\epsilon} \in H^1(0,T; H^1_0(\Omega;\bR^d))$ such that $u^{\epsilon}(0)=u_0^\epsilon$ in $H_0^1(\Omega;\bR^d)$ and, for a.e.\ $t\in [0,T]$,
\begin{equation}\label{eq:two_scale_weak_form}
    (\dot{u}^{\epsilon}(t), v)_{\bV^{\epsilon},\, \mathrm{sym}} + (u^{\epsilon}(t), v)_{\bE^{\epsilon},\,\mathrm{sym}} = (f(t), v)_{L^2}\quad \forall v\in H^1_0(\Omega;\bR^d).
\end{equation}
\end{namedparagraph}

\begin{namedparagraph}{Admissible Initial Conditions}
A family $\{u_0^\epsilon\}_{\epsilon>0}\subset H_0^1(\Omega;\bR^d)$ is admissible if there exists $g_0\in L^2(\Omega;\bR^d)$ such that each $u_0^\epsilon$ uniquely solves $(u_0^\epsilon ,v)_{\bE^{\epsilon},\, \mathrm{sym}} =(g_0, v)_{L^2}$ for all $v\in H_0^1(\Omega;\bR^d)$.
\end{namedparagraph}

\begin{theorem}[Viscoelastic Homogenization] \label{thm:viscoelastic_homogenization}
Suppose $\Omega$ is Lipschitz, $\{u_0^\epsilon\}_{\epsilon>0}$ is an admissible family generated by $g_0\in L^2(\Omega;\bR^d)$, $\bV,\bE\in L^{\infty}_{\mathrm{per}}(\mathbb{T}^d;\cV^4_{+})$, and $f \in L^2(0,T; L^2(\Omega; \bR^d))$. Define the homogenized viscosity and elasticity $\bV^{\mathrm{hom}}(\bV), \bE^{\mathrm{hom}}(\bE)\in\cV^4_{+}$ through the symmetric positive-definite bilinear forms:
\begin{align}\label{eq:homogenized_tensors}
    \left(\overline{e}_1,\bV^{\mathrm{hom}}(\bV)\overline{e}_2\right)_F = \left(\overline{e}_1, \cP_{\bV}^{\perp} \overline{e}_2\right)_{\bV}, \quad \left(\overline{e}_1, \bE^{\mathrm{hom}}(\bE) \overline{e}_2\right)_F = \left(\overline{e}_1, \cP_{\bE}^{\perp} \overline{e}_2\right)_{\bE}\quad \forall \overline{e}_1, \overline{e}_2 \in \cV^2,
\end{align}
and for each $t\ge0$, define the memory kernel $\bK(t;\bV,\bE)\in \cV^4$ through the symmetric positive semidefinite bilinear form:
\begin{align}\label{eq:memory_kernel}
    \left(\overline{e}_1, \bK(t; \bV, \bE) \overline{e}_2\right)_F &= \left(\cR_{\bV,\bE}\exp(-t\cR_{\bV,\bE})\cD_{\bV,\bE} \overline{e}_1, \cD_{\bV,\bE} \overline{e}_2\right)_{\bV}\quad\forall \overline{e}_1, \overline{e}_2 \in \cV^2.
\end{align}

The unique weak solution $u^{\epsilon}$ of \cref{eq:two_scale_weak_form} converges weakly in $H^1(0,T; H^1_0(\Omega; \bR^d))$ to $u^{(0)}$, the unique weak solution of the macroscopic problem \cref{eqns:homogenized_problem}, closed by the homogenized model with memory
\begin{align}\label{eq:homogenized_model_identified}
    \cH(\{\overline{e}(s)\}_{s\in[0, t]};\bV, \bE) \coloneqq \bV^{\mathrm{hom}}(\bV) \dot{\overline{e}}(t) + \bE^{\mathrm{hom}}(\bE)\overline{e}(t) + \int_{0}^{t} \bK(t-s; \bV, \bE)\dot{\overline{e}}(s) \, \mathrm{d}s,
\end{align}
with homogenized initial condition $u_0^{(0)}$ solving $(u^{(0)}_0,v)_{\bE^{\rm hom}(\bE),\, \mathrm{sym}}=(g_0,v)_{L^2}$ for all $v\in H_0^1(\Omega;\bR^d)$.
\end{theorem}

\begin{proof}[\bf Proof Sketch]
We follow \cite[Theorem 3.1]{francfort1986homogenization}. First, consider $f=0$. For admissible initial data, the semigroup estimate \cite[(3.20)]{francfort1986homogenization} gives an $\epsilon$-independent bound for $u^\epsilon$ in $W^{1,\infty}(0,\infty;H_0^1(\Omega;\bR^d))$. After passing to a subsequence, $u^\epsilon$ converges in weak-* to $u^{(0)}$ on this space. Their Laplace transforms therefore satisfy $\hat u^\epsilon(s)\rightharpoonup\hat u^{(0)}(s)$ in $H_0^1(\Omega;\bR^d)$ for every $s>0$. By admissibility, $\hat v^\epsilon(s)\coloneqq s\hat u^\epsilon(s)-u_0^\epsilon$ solves a static elliptic problem with coefficient $\bA_s^\epsilon(x)\coloneqq s^{-1}\bE(x/\epsilon)+\bV(x/\epsilon)$ and load $-g_0/s$. Periodic homogenization identifies the limits of these problems and of the initial data.

Next, we identify the constitutive law. Expressing the dissipative corrector of \cite{francfort1986homogenization} in terms of $\cR_{\bV,\bE}$ and $\cD_{\bV,\bE}$ gives the kernel \cref{eq:memory_kernel}. The cell-problem identity \cite[(3.12)]{francfort1986homogenization} then reads $\bA_s^{\mathrm{hom}}=s^{-1}\bE^{\mathrm{hom}}(\bE)+\bV^{\mathrm{hom}}(\bV)+\hat{\bK}(s)$, where $\hat{\bK}$ is the Laplace transform of $\bK$. Uniqueness of the homogenized elliptic problems and injectivity of the Laplace transform recover \cref{eq:homogenized_model_identified} and identify the limit uniquely, giving convergence of the whole family. We establish the properties of the homogenized tensors and kernel in \cref{app:homogenized_tensors_kernel_properties}.

Finally, for general $f\in L^2(0,T;L^2(\Omega;\bR^d))$, Duhamel's formula expresses the forcing contribution as a time integral of homogeneous evolutions whose initial data solve viscous elliptic problems. Applying the transformed elliptic argument to these evolutions and using the uniform semigroup bounds to pass to the limit under the integral yields the stated result on $[0,T]$; see \cite[(3.34)--(3.39)]{francfort1986homogenization}.
\end{proof}

\begin{remark}
    {The identity between the unit-cell representation \cref{eq:cell_constitutive_map} and the constitutive law \cref{eq:homogenized_model_identified} is proved in \cref{app:unit_cell_homogenization} by writing the unit-cell problem as an ordinary differential equation on $\cX$.}
\end{remark}

\begin{remark}
    The kernel representation shows that the mismatch operator $\mathcal D_{\mathbb V,\mathbb E}$ drives microscopic strain relaxation, while $\mathcal R_{\mathbb V,\mathbb E}$ determines its decay rates. The memory kernel vanishes when the viscous and elastic strain correctors coincide for every macroscopic strain. In particular, this holds when $\mathbb V=c\mathbb E$ for a spatially constant $c>0$,
\end{remark}

\subsection{Spectral Representation of the Memory Kernel}\label{subsec:spectral_rep}

The operator form of the memory kernel in \cref{eq:memory_kernel} leads to a representation of $\bK$ as a Laplace--Stieltjes integral over relaxation rates. The next proposition formalizes this representation together with two useful properties: the rates lie in an interval determined by the coercivity and boundedness constants of the microstructural viscosity and elasticity, and the total variation of the spectral measure is bounded in terms of the same constants.

\begin{proposition}[Spectral Representation of the Memory Kernel]
\label{prop:spectral_rep}
Let $\bV,\bE\in L^{\infty}_{\mathrm{per}}(\bT^d;\cV^4_{+})$, and let $\alpha_{\bV},\beta_{\bV},\alpha_{\bE},\beta_{\bE}>0$ satisfy
\begin{equation}\label{eq:constants_spectral_rep}
\alpha_{\bV}\|e\|_{L^2_{\mathrm{per}}}^2 \le \|e\|_{\bV}^2 \le \beta_{\bV}\|e\|_{L^2_{\mathrm{per}}}^2,
\qquad
\alpha_{\bE}\|e\|_{L^2_{\mathrm{per}}}^2 \le \|e\|_{\bE}^2 \le \beta_{\bE}\|e\|_{L^2_{\mathrm{per}}}^2
\end{equation}
for all $e\in L^2_{\mathrm{per}}(\bT^d;\cV^2)$. Then there exists a finite positive semidefinite $\cV^4$-valued Borel measure $\nu_{\bV,\bE}$, supported on $\cI_{\bV,\bE}\coloneqq[\alpha_{\bE}/\beta_{\bV},\beta_{\bE}/\alpha_{\bV}]$, such that
\begin{align}
\bK(t;\bV,\bE)&=\int_{\cI_{\bV,\bE}} \lambda \exp(-\lambda t)\,\mathrm{d}\nu_{\bV,\bE}(\lambda)
\quad \forall t\ge 0,\label{eq:K_laplace}\\
\|\nu_{\bV,\bE}\|_{\mathrm{TV}(\cI_{\bV,\bE})}
&\le
d_{\cV^2}\frac{\beta_{\bV}\beta_{\bE}}{\alpha_{\bE}}, \label{eq:uniform_TV_bound_nu_explicit}
\end{align}
where the total variation of the measure is defined as
\begin{align}\label{eq:tv_norm_def}
\|\nu_{\bV,\bE}\|_{\mathrm{TV}(\cI_{\bV,\bE})}
\coloneqq
\sup\bigg\{\sum_{i=1}^{m}\|\nu_{\bV,\bE}(B_i)\|_{F}\ \Big\vert\ &m\in\bN,\,
B_1,\dots,B_m\ \text{are Borel subsets of }\cI_{\bV,\bE},\nonumber\\
&B_i\cap B_j=\emptyset\ \forall i\ne j,\,\bigcup_{i=1}^{m}B_i=\cI_{\bV,\bE}\bigg\}.
\end{align}
\end{proposition}
The proof is given in \cref{app:spectral_rep}. The interval in the proposition gives the range of relaxation rates that are relevant to the memory kernel, while the total-variation norm bound controls the total amount of memory effects arising from homogenization. While the constants $\alpha_{\bV},\beta_{\bV},\alpha_{\bE},\beta_{\bE}$ are properties of the specific microstructural instances, they are shared among a wide class of microstructures that satisfy \cref{eq:constants_spectral_rep}. \Cref{prop:spectral_rep} shows that the spectral support and the total variation bound are common to this class; this is a key component of the approximation theory presented in the next section.

\section{Universal Approximation and Internal-Variable Complexity}\label{sec:universal_approximation}

We now use the spectral representation of the memory kernel from \cref{subsec:spectral_rep} to address the approximation and internal-variable complexity questions for Kelvin--Voigt viscoelasticity. The key step is to introduce a Prony-series internal-variable model that approximates the homogenized constitutive law $\cH$ in \cref{eq:homogenized_model_identified} uniformly over admissible microstructures and strain trajectories and depends continuously on the microstructure. This intermediate model links the fading-memory representation of $\cH$ to the neural operator internal-variable model $\cH_{\rm NO}$ in \cref{def:learned_model}.

In \cref{subsec:ua_settings}, we specify the admissible microstructures, strain trajectories, and neural-operator response classes. In \cref{subsec:ua_main_result}, we state the universal approximation theorem and sketch its proof, with the full argument deferred to \cref{app:ua}. In \cref{subsec:ua_internal_variable_dimension}, we derive the sufficient internal-variable dimension bound, and in \cref{subsec:ua_limits}, we discuss continuity constraints on microstructure-dependent relaxation rates.

In the weak macroscopic problem, the strain and stress are Bochner functions $\overline e,\overline\sigma:[0,T]\to L^2(\Omega;\cV^2)$. Since $\cH$ is local in space, we pass to a material-point trajectory formulation: for a.e.\ $x\in\Omega$, we write $\overline{e},\overline{\sigma}:[0,T]\to\cV^2$. {Moreover, the constitutive law \cref{eq:homogenized_model_identified} is independent of the terminal time $T$. For the approximation results in this section, we take $\cT=[0,\infty)$ as the time domain for this law and the neural operator model \cref{eq:learned_model_full}.}

\subsection{Settings}\label{subsec:ua_settings}

We specify the admissible classes of microstructures and strain trajectories, together with the neural-operator response classes that enter the approximation theorem. 

\begin{namedparagraph}{Admissible Microstructural Properties}
For $\bA\in L^{1}_{\mathrm{per}}(\bT^d;\mathcal{V}^{4})$, we define the total variation seminorm by duality \cite{ambrosio2000functions},
\begin{equation}\label{eq:TV_tensor}
|\bA|_{\mathrm{TV}}
:=\sup\left\{
\int_{\bT^d}\left(\bA(y),\nabla_y\cdot \phi(y)\right)_F\,\mathrm{d}y
\;\middle|\;
\phi\in C^{1}_{\mathrm{per}}\big(\bT^d;\mathcal{V}^{4}\otimes\bR^{d}\big),\ 
\|\phi\|_{L^\infty_{\mathrm{per}}}\leq 1
\right\},
\end{equation}
and we set $\|\bA\|_{\mathrm{BV}}:=\|\bA\|_{L^{1}_{\mathrm{per}}}+|\bA|_{\mathrm{TV}}$ and $\mathrm{BV}:=\{\bA\in L^{1}_{\mathrm{per}}(\bT^d;\cV^{4})\mid\|\bA\|_{\mathrm{BV}}<\infty\}$. Given constants $\alpha_{\rm visc},\beta_{\rm visc},\gamma_{\rm visc},\alpha_{\rm elas},\beta_{\rm elas},\gamma_{\rm elas}>0$, the admissible set
\begin{equation}\label{eq:admissible_microstructure}
    \mathcal{M}_{\substack{\alpha_{\rm visc},\beta_{\rm visc},\gamma_{\rm visc}\\\alpha_{\rm elas},\beta_{\rm elas},\gamma_{\rm elas}}}\subset 
\Big(\mathrm{BV}\cap L^\infty_{\mathrm{per}}(\bT^d;\mathcal{V}^{4})\Big)^2
\end{equation}
consists of pairs $(\bV,\bE)$ satisfying:
\begin{subequations} 
\begin{gather}
\alpha_{\rm visc}\|e\|_{L^2_{\mathrm{per}}}^{2} \leq \|e\|_{\bV}^{2} \leq \beta_{\rm visc}\|e\|_{L^2_{\mathrm{per}}}^{2},\,
\alpha_{\rm elas}\|e\|_{L^2_{\mathrm{per}}}^{2} \leq \|e\|_{\bE}^{2} \leq \beta_{\rm elas}\|e\|_{L^2_{\mathrm{per}}}^{2} ,\label{eq:uniform_ellipticity_microstructure}\\
\|\bV\|_{\mathrm{BV}} \leq \gamma_{\rm visc}, \quad \|\bE\|_{\mathrm{BV}} \leq \gamma_{\rm elas}.\label{eq:bv_bound_microstructure}
\end{gather}
for all $e\in L^2_{\mathrm{per}}(\bT^d;\cV^{2})$.
\end{subequations}
\end{namedparagraph}
\begin{remark}[Bounded Variation]
 
 Bounded variation is a natural assumption on the material microstructure which gives compactness in $(L^2_{\mathrm{per}})^2$, but smoother classes (e.g., uniformly bounded in $W^{1,p}$) also have this property. We consider this class of materials as it is natural for composite materials to exhibit discontinuous jumps in microstructural properties, but fractal-type oscillations are excluded. The bounded variation assumption is not needed in the proof of the spectral representation theorem above. However, the compactness it provides in $(L^2_{\mathrm{per}})^2$ is essential for establishing the universal approximation theorem below; see \cref{lem:bv_compact}.
 
\end{remark}

\begin{namedparagraph}{Admissible Strain Trajectories}
    Since the constitutive model depends on both $\overline e$ and $\dot{\overline e}$, we impose a uniform bound on both quantities. Let $\delta_{\rm stress}>0$ denote the prescribed stress scale. We define the set of admissible macroscopic strain trajectories by:
    \begin{equation}\label{eq:admissible_strain}
        \mathcal{E}_{\delta_{\rm stress}, \beta_{\rm visc},\beta_{\rm elas}} := \left\{ \overline{e} \in W^{1, \infty}(\mathcal{T}; \mathcal{V}^2) \mid \ \beta_{\rm elas}\|\overline{e}\|_{L^{\infty}} + \beta_{\rm visc}\|\dot{\overline{e}}\|_{L^{\infty}} \leq\delta_{\rm stress} \right\}.
    \end{equation}
We emphasize that this admissible class excludes idealized step-strain tests.
\end{namedparagraph}

\begin{namedparagraph}{Neural-Operator Response Classes}
We state the theorem for an abstract class of neural operators, rather than for a particular architecture. For integers $d_{\mathrm{in}},d_{\mathrm{out}}\in\bN$, let
\begin{equation*}
    \mathfrak{N}\subset\left\{\Phi:\bR^{d_{\mathrm{in}}}\times L^2_{\mathrm{per}}(\bT^d;\bR^{2d_{\cV^4}}) \to \bR^{d_{\mathrm{out}}}\right\}
\end{equation*}
be a class of continuous mappings with the following approximation properties.

\begin{enumerate}[label=(A\arabic*)]
    \item\label{no:ua} \textit{Universal Approximation.} For every compact $\cK_{\mathrm{vec}}\subset \bR^{d_{\mathrm{in}}}$, every compact $\cK_{\mathrm{fun}}\subset L^2_{\mathrm{per}}(\bT^d;\bR^{2d_{\cV^4}})$, every continuous $g:\cK_{\mathrm{vec}}\times \cK_{\mathrm{fun}}\to\bR^{d_{\mathrm{out}}}$, and every $\varepsilon>0$, there exists $\Phi\in\mathfrak{N}$ such that $\sup_{\cK_{\mathrm{vec}}\times \cK_{\mathrm{fun}}}\|\Phi-g\|\le\varepsilon$.

    \item\label{no:lip} \textit{Locally Lipschitz Universal Approximation.} The class $\mathfrak{N}$ satisfies~\ref{no:ua}. In addition, the approximating $\Phi$ can be chosen so that, for every $M\in \cK_{\mathrm{fun}}$, the map $v\mapsto \Phi(v,M)$ is locally Lipschitz.
\end{enumerate}
For $n_{\xi}\in\bN$, we define the stress and kinetic response classes as follows.
\begin{itemize}
    \item \textit{Stress Response Class $\mathfrak{F}(n_{\xi})$}. This class consists of all $\cF_{\mathrm{NO}}$ of the form \cref{eq:stress_response} with $\Phi_{\cF}\in\mathfrak{N}_{\cF}$, where $\mathfrak{N}_{\cF}$ satisfies~\ref{no:ua} with $(d_{\mathrm{in}},d_{\mathrm{out}})=(2d_{\cV^2}+n_{\xi},\,d_{\cV^2})$.
    \item \textit{Kinetic Response Class $\mathfrak{G}(n_{\xi})$}. This class consists of all $\cG_{\mathrm{NO}}$ of the form \cref{eq:kinetic_response} with $\Phi_{\cG}\in\mathfrak{N}_{\cG}$, where $\mathfrak{N}_{\cG}$ satisfies~\ref{no:lip} with $(d_{\mathrm{in}},d_{\mathrm{out}})=(d_{\cV^2}+n_{\xi},\,n_{\xi})$.
\end{itemize}
Only the kinetic response requires~\ref{no:lip}, as this local Lipschitz property is used to control the internal-variable evolution in \cref{eq:learned_model_full}.
\end{namedparagraph}
\begin{remark}[Architecture]\label{rmk:general_admissibility}
    The set $\mathfrak{N}$ contains the maps realizable by a chosen neural-operator architecture as its weights and hyperparameters vary. Condition~\ref{no:ua} is an expressivity requirement: the architecture approximates continuous maps uniformly on compact sets \cite{kovachki2023neural}. Approximation results for a range of neural-operator architectures and function-space settings are available; see, e.g., \cite{chen1995universal,bhattacharya2021pca,kovachki2021universal,nelsen2021random,luo2025dimension,lanthaler2025nonlocality}. Condition~\ref{no:lip} strengthens~\ref{no:ua} by requiring local Lipschitz dependence on the vector input. In practice,~\ref{no:lip} holds when the relevant architectural components are locally Lipschitz, for example, when using ReLU or GELU activations.
\end{remark}

\begin{remark}[Fourier Neural Mapping]\label{rmk:fnm_admissible}
    The FNM architecture in \cref{subsec:viscoelastic_neural_operator_model} satisfies both~\ref{no:ua} and~\ref{no:lip} when equipped with smooth, globally Lipschitz, non-polynomial activations. Universal approximation follows from~\cite{huang2025operator} for linear lifting and projection maps $S$ and $Q$, while the Lipschitz condition follows from the bounded linear layers and Lipschitz activations. If $S$ and $Q$ are instead chosen as feedforward neural networks, locally Lipschitz activations preserve~\ref{no:lip}.
\end{remark}

\subsection{Universal Approximation Theorem}\label{subsec:ua_main_result}

We state the main result, which establishes the existence of a neural operator constitutive model $\cH_{\rm NO}$, as defined in \cref{def:learned_model}, that approximates the homogenized constitutive law of Kelvin--Voigt viscoelasticity to any prescribed accuracy, uniformly in time and over a wide range of microstructural properties and strain trajectories.

\begin{theorem}[Universal Approximation of Viscoelastic Homogenization]\label{thm:ua}
Let $d\in\{1,2,3\}$. For any constants $\alpha_{\rm visc},\beta_{\rm visc},\gamma_{\rm visc},\alpha_{\rm elas},\beta_{\rm elas}, \gamma_{\rm elas},\delta_{\rm stress}>0$ such that $0<\alpha_{\rm elas}\leq \beta_{\rm elas}$ and $0<\alpha_{\rm visc}\leq \beta_{\rm visc}$,  let
\begin{equation*}
\cM\coloneqq\cM_{\substack{\alpha_{\rm visc},\beta_{\rm visc},\gamma_{\rm visc}\\\alpha_{\rm elas},\beta_{\rm elas},\gamma_{\rm elas}}},\qquad
\cE\coloneqq\cE_{\delta_{\rm stress},\beta_{\rm visc},\beta_{\rm elas}}.
\end{equation*}
For any $\varepsilon>0$, there exist $n_\xi\in\bN$ and a pair of stress and kinetic response mappings
\begin{equation*}
    (\cF_{\mathrm{NO}}, \cG_{\mathrm{NO}})\in \mathfrak{F}(n_\xi)\times \mathfrak{G}(n_\xi),
\end{equation*}
such that
\begin{equation}\label{eq:ua}
\sup_{ (\bV,\bE)\in\cM}\;
\sup_{ \overline e\in\cE}\;
\esssup_{  t\in\cT}\;
\left\|
\cH_{\mathrm{NO}}(\{\overline e(s)\}_{s\in[0,t]};\bV,\bE)
-
\cH(\{\overline e(s)\}_{s\in[0,t]};\bV,\bE)
\right\|_F
\le \varepsilon.
\end{equation}
\end{theorem}

\begin{proof}[\bf Proof Sketch]
The proof proceeds in four steps; full details are in \cref{app:ua}.
The proof follows the approach introduced in \cite{bhattacharya2023learning,bhattacharya2025learning}, but significantly extends it to facilitate approximation of the spectral measure $\nu_{\bV,\bE}$ underlying the formula for the memory kernel appearing in \Cref{prop:spectral_rep}; the preceding papers used approximation by piecewise-constant materials, whose memory kernels admit exact finite Prony-series representations in one dimension.

\begin{enumerate}[label=(S\arabic*), ref=(S\arabic*)]
    \item\label{item:ua_proof_s1} Establish the compactness of $\cM$ in $(L^2_{\rm per})^2$ and the continuity of the homogenized tensors with respect to $(\bV,\bE)$.
    \item\label{item:ua_proof_s2} Construct an intermediate Prony-series internal-variable model $\cH_{\mathrm{IV}}$ by discretizing the spectral representation of the memory kernel on a common interval of relaxation rates.
    \item\label{item:ua_proof_s3} Approximate the associated stress and kinetic responses of $\cH_{\mathrm{IV}}$ by those of the neural operator model $\cH_{\mathrm{NO}}$ on compact input sets and control the induced trajectory error.
    \item\label{item:ua_proof_s4} Choose $n_{\xi}$ and the error tolerances for $\cF_{\mathrm{NO}}$ and $\cG_{\mathrm{NO}}$ to obtain \cref{eq:ua}.
\end{enumerate}

\paragraph{\ref{item:ua_proof_s1}: Compactness and Continuity}
The uniform BV and $L^\infty_{\rm per}$ bounds make $\cM$ compact in $(L^2_{\rm per})^2$ (\cref{lem:bv_compact}). Convergence of the material fields in $L^2_{\rm per}$ implies strong convergence of the associated unit-cell operators (\cref{lem:strong_conv_tensor,prop:strong_conv}). Consequently, the homogenized tensors $\bV^{\rm hom}$ and $\bE^{\rm hom}$ depend continuously on the microstructure, and the spectral measure in the memory representation is weakly continuous (\cref{cor:continuity_tensors,lem:weak_convergence}).

{
\paragraph{\ref{item:ua_proof_s2}: Internal-Variable Approximation}
We replace the spectral integral in \cref{prop:spectral_rep} by a finite sum of exponentials, or Prony series. All relaxation rates lie in the common interval $\cI=[\lambda_{\min},\lambda_{\max}]$, where $\lambda_{\min}=\alpha_{\rm elas}/\beta_{\rm visc}$ and $\lambda_{\max}=\beta_{\rm elas}/\alpha_{\rm visc}$.

We introduce the following weighted spectral measure and relaxation response:
\begin{align*}
\mathrm d\mu_{\bV,\bE}(\lambda)\coloneqq\frac{\lambda}{\beta_{\rm elas}/\beta_{\rm visc}+\lambda}\,\mathrm d\nu_{\bV,\bE}(\lambda),\qquad\cW_{\overline e,t}(\lambda)\coloneqq\left(\frac{\beta_{\rm elas}}{\beta_{\rm visc}}+\lambda\right)\int_0^t\exp(-\lambda(t-s))\dot{\overline e}(s)\,\mathrm ds.
\end{align*}
This weighting will allow us to bound the stress error introduced by the Prony-series approximation with a constant independent of the elasticity contrast $\beta_{\rm elas}/\alpha_{\rm elas}$ appearing in \cref{eq:uniform_TV_bound_nu_explicit} (\cref{lem:weighted_spectral_measure,lem:weighted_response}). The memory contribution to stress can then be written as
\begin{equation}\label{eq:weighted_kernel_representation}
\int_0^t\bK(t-s)\dot{\overline e}(s)\,\mathrm ds
=\int_{\cI}\mathrm d\mu_{\bV,\bE}(\lambda)\cW_{\overline e,t}(\lambda).
\end{equation}
We interpolate $\cW_{\overline e,t}$ in $\ln\lambda$ on overlapping subintervals of bounded logarithmic width, using Chebyshev nodes on each subinterval (\cref{fig:local_chebyshev,def:prony_series_approximation}). Analyticity of the response and the bounded widths give a uniform geometric convergence rate (\cref{lem:uniform_prony_discretization_beta}). The overlap permits continuous nonnegative weights summing to one, which combine the local interpolants without accumulating their errors.

Integrating the combined interpolant against $\mu_{\bV,\bE}$ gives a model with memory kernel $\bK^{(n)}(t)=\sum_{k=1}^n\bK_k\exp(-\lambda_kt)$. The rates $\lambda_k$ are shared by all admissible microstructures, while the coefficients $\bK_k(\bV,\bE)$ depend continuously on the microstructure (\cref{prop:uniform_prony_homogenized_model}), as needed in \ref{item:ua_proof_s3}. Each exponential term introduces one tensorial internal variable $q_k$, giving the model $\cH_{\rm IV}$:
\begin{subequations}\label{eqns:prony_internal_varaible_model}
\begin{align}
    \overline\sigma_{\rm IV}(t)&=\bV^{\rm hom}\dot{\overline e}(t)+\bE^{\rm hom}\overline e(t)+\sum_{k=1}^n q_k(t),\\
    \dot q_k(t)&=\bK_k\dot{\overline e}(t)-\lambda_kq_k(t),\qquad q_k(0)=0,\qquad k=1,\dots,n.
\end{align}
\end{subequations}

\begin{figure}[H]
    \centering
    \includegraphics[width=0.9\linewidth]{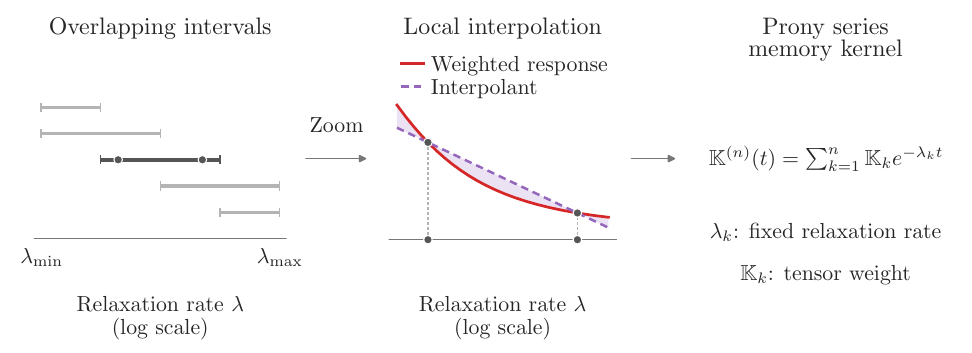}
    \caption{Prony-series construction by local Chebyshev interpolation \ref{item:ua_proof_s2}. The middle panel shows one component of the weighted relaxation response as $\lambda$ varies, for a given strain history and time. Two nodes illustrate interpolation; shading marks the error. The resulting relaxation rates are shared across microstructures, while the tensor coefficients depend continuously on the microstructure.}
    \label{fig:local_chebyshev}
\end{figure}

Define the logarithmic relaxation span by $\mathsf R\coloneqq\ln(\lambda_{\max}/\lambda_{\min})$. For $\mathsf R>0$, we use $J=1+\lceil\mathsf R\rceil$ subintervals and $m$ nodes per subinterval, giving $n=Jm$ and $n_\xi=Jmd_{\cV^2}$. \Cref{prop:uniform_prony_homogenized_model} gives the uniform stress error bound
\begin{equation*}
    \sup_{(\bV,\bE)\in\cM}\,\sup_{\overline e\in\cE}
    \|\overline\sigma-\overline\sigma_{\rm IV}\|_{L^\infty(\cT;\cV^2)}<\kappa_{\rm IV}2^{-(m-1)},
\end{equation*}
where $\kappa_{\rm IV}$ depends only on $d$ and $\delta_{\rm stress}$. Increasing $m$ therefore gives any prescribed accuracy, uniformly over admissible microstructures, strain trajectories, and time. In the next subsection, we use this estimate to bound the number of tensorial internal variables sufficient for a prescribed stress accuracy. \Cref{rmk:prony_construction_choices} explains how the weighted representation \cref{eq:weighted_kernel_representation} and the local interpolation introduced in this step improve that bound.
}

\paragraph{\ref{item:ua_proof_s3}: Neural Operator Approximation}
Uniform bounds on strain, strain rate, and the internal variables of $\cH_{\rm IV}$ let us choose a compact set of vector inputs with a margin for internal-variable errors (\cref{lem:ua_step3_targets}). The continuity established in \ref{item:ua_proof_s1}--\ref{item:ua_proof_s2} and compactness of $\cM$ permit uniform approximation of the stress and kinetic response mappings of $\cH_{\rm IV}$ by neural operators, with errors $\varepsilon_{\cF}$ and $\varepsilon_{\cG}$ (\cref{lem:ua_step3_no}). Choosing $\cG_{\rm NO}$ locally Lipschitz in its vector inputs gives a unique local solution $\xi_{\rm NO}$.

{The terms $-\lambda_kq_k$ in \cref{eqns:prony_internal_varaible_model} yield a bound, uniform in time, on the difference between the internal variables of $\cH_{\rm NO}$ and $\cH_{\rm IV}$. For sufficiently small $\varepsilon_{\cG}$, the norm of this difference remains below the chosen margin, so the inputs remain in this compact set and the solution extends to all of $\cT$. Thus \cref{prop:ua_step3_stability} gives
\begin{equation*}
    \sup_{(\bV,\bE)\in\cM}\,\sup_{\overline e\in\cE}
    \|\overline\sigma_{\rm NO}-\overline\sigma_{\rm IV}\|_{L^\infty(\cT;\cV^2)}
    \le\varepsilon_{\cF}+\frac{\sqrt n}{\lambda_{\min}}\varepsilon_{\cG}.
\end{equation*}
}
The first term bounds the error in approximating the stress response mapping; the second accounts for the error in the internal variables.

\paragraph{\ref{item:ua_proof_s4}: Combination of the Approximation Errors}
Choose any $\varepsilon_{\rm IV}\in(0,\varepsilon)$ and then $n$ so that the stress error in \ref{item:ua_proof_s2} is at most $\varepsilon_{\rm IV}$. For this choice of $n$, choose $\varepsilon_{\cF}$ and $\varepsilon_{\cG}$ so that the error in \ref{item:ua_proof_s3} is at most $\varepsilon-\varepsilon_{\rm IV}$ and the stability condition holds. The triangle inequality gives \cref{eq:ua}.
\end{proof}

\subsection{On the Sufficient Dimension of the Internal Variable}\label{subsec:ua_internal_variable_dimension}

{The proof of \cref{thm:ua} constructs the intermediate model $\cH_{\rm IV}$ from a Prony series, then approximates it by the neural operator model $\cH_{\rm NO}$ {without increasing the dimension of the internal variable model implied by $\cH_{\rm IV}$}. Each term in the series introduces one tensorial internal variable, so bounding $n_\xi$ reduces to bounding the number of terms needed for the target accuracy. The local Chebyshev interpolation estimate in \ref{item:ua_proof_s2} provides this bound. For a prescribed stress error tolerance $\varepsilon>0$, the following theorem expresses it in terms of the logarithmic relaxation span $\mathsf R$ and the relative stress error $\mathsf E\coloneqq\varepsilon/\delta_{\rm stress}$.}

{
\begin{theorem}[Sufficient Internal-Variable Dimension]\label{thm:sufficient_dimension}
    In the setting of \cref{thm:ua}, use the parameters $\mathsf R$ and $\mathsf E$ defined above. If $\mathsf R=0$, set $N\coloneqq d_{\cV^2}$. If $\mathsf R>0$, define
    \begin{equation*}
        N\coloneqq d_{\cV^2}(1+\lceil\mathsf R\rceil)
        \left(1+\left\lceil
        \log_2^+\left(\frac{8d_{\cV^2}(1+\sec(3/4))}{\mathsf E}\right)
        \right\rceil\right).
    \end{equation*}
    Here $\log_2^+(x)\coloneqq\max\{0,\log_2 x\}$. Then there exist response mappings $(\cF_{\rm NO},\cG_{\rm NO})\in\mathfrak F(N)\times\mathfrak G(N)$ attaining \cref{eq:ua}. Thus $N$ is a sufficient internal-variable dimension.
\end{theorem}

\begin{proof}[\bf Proof]
If $\mathsf R=0$, take $n=1$; the approximation by $\cH_{\rm IV}$ is exact. Otherwise, choose $m=1+\lceil\log_2^+(\kappa_{\rm IV}/\varepsilon)\rceil$ and $n=(1+\lceil\mathsf R\rceil)m$, with $\kappa_{\rm IV}$ as in \cref{eq:prony_homogenized_error_uniform}. The strict bound in \cref{eq:prony_homogenized_error_uniform} ensures that the uniform error in the stress predicted by $\cH_{\rm IV}$ is less than $\varepsilon$. In either case, choose $\varepsilon_{\rm IV}$ strictly between this error and $\varepsilon$. The tolerances in \cref{eq:ua_error_tolerances} then bound the remaining stress error by $\varepsilon-\varepsilon_{\rm IV}$ through \cref{prop:ua_step3_stability}. The triangle inequality yields \cref{eq:ua} with $nd_{\cV^2}=N$.
\end{proof}
}

\begin{table}[!htb]
    \centering
\setlength{\tabcolsep}{5pt}
\renewcommand{\arraystretch}{2.5}

\begin{tabular}{ | c || c | c |}\hline
    Name & Definition & Interpretation \\\hline\hline
    \makecell{Logarithmic\\Relaxation Span} &
    $\mathsf{R}\coloneqq \ln\left(\dfrac{\beta_{\rm visc}\beta_{\rm elas}}{\alpha_{\rm visc}\alpha_{\rm elas}}\right)
    $ & \makecell{Logarithm of the ratio of the largest \\ to the smallest relaxation rate}\\\hline
    \makecell{Viscosity Contrast} &
    $\mathsf{C}_{\rm visc}\coloneqq \dfrac{\beta_{\rm visc}}{\alpha_{\rm visc}}
    $ & \makecell{Ratio between the maximum and minimum \\ viscosity}\\\hline
    \makecell{Elasticity Contrast} &
    $\mathsf{C}_{\rm elas}\coloneqq \dfrac{\beta_{\rm elas}}{\alpha_{\rm elas}}
    $ & \makecell{Ratio between the maximum and minimum \\ elasticity}\\\hline
    Relative Stress Error & $\mathsf{E}\coloneqq\dfrac{\varepsilon}{\delta_{\rm stress}}$ & \makecell{Ratio between the stress error tolerance \\ and the prescribed stress scale}\\\hline
\end{tabular}
    \caption{Non-dimensional groups for the sufficient internal-variable dimension.}
\label{tab:non-dimensional}
\end{table}

{Define the viscosity and elasticity contrasts by $\mathsf C_{\rm visc}\coloneqq\beta_{\rm visc}/\alpha_{\rm visc}$ and $\mathsf C_{\rm elas}\coloneqq\beta_{\rm elas}/\alpha_{\rm elas}$, respectively. The non-dimensional parameters are summarized in \cref{tab:non-dimensional}. The bound in \cref{thm:sufficient_dimension} has the form
\begin{equation*}
    N=\mathcal O_d\left((1+\mathsf R)\left(1+\ln^+(1/\mathsf E)\right)\right),
    \qquad
    \mathsf R=\ln(\mathsf C_{\rm visc}\mathsf C_{\rm elas}),
\end{equation*}
where $\ln^+x=\max\{0,\ln x\}$ and the implicit constant depends only on $d$. In particular,
\begin{equation}\label{eq:sufficient_dimension_rates}
\begin{aligned}
N&=\mathcal O\bigl(\ln(1/\mathsf E)\bigr), &&\mathsf E\to0,\quad \mathsf C_{\rm visc},\mathsf C_{\rm elas}\text{ held constant},\\
N&=\mathcal O\bigl(\ln\mathsf C_{\rm visc}\bigr), &&\mathsf C_{\rm visc}\to\infty,\quad \mathsf E,\mathsf C_{\rm elas}\text{ held constant},\\
N&=\mathcal O\bigl(\ln\mathsf C_{\rm elas}\bigr), &&\mathsf C_{\rm elas}\to\infty,\quad \mathsf E,\mathsf C_{\rm visc}\text{ held constant}.
\end{aligned}
\end{equation}
}

{
\begin{remark}[Choices in the Prony-Series Construction]\label{rmk:prony_construction_choices}
Simpler Prony-series constructions also establish \cref{thm:ua}, but reweighting and local interpolation improve the dependence of the sufficient internal-variable dimension on the material contrasts.

Without reweighting, interpolation with $m$ Chebyshev nodes per subinterval and the total-variation bound for $\nu_{\bV,\bE}$ in \cref{eq:uniform_TV_bound_nu_explicit} give a relative stress error bound of order $\mathsf C_{\rm elas}2^{-(m-1)}$. At fixed relative accuracy $\mathsf E$ and viscosity contrast $\mathsf C_{\rm visc}$, choosing $m$ from this estimate gives a sufficient internal-variable dimension of order $(\ln\mathsf C_{\rm elas})^2$.

For the weighted representation, consider a single interpolant with $n$ Chebyshev nodes on $[\ln\lambda_{\min},\ln\lambda_{\max}]$. The Chebyshev interpolation estimate for analytic functions \cite[Theorem 8.2]{trefethen2019approximation} yields a relative stress error bound of order $\mathsf R(1+1/\mathsf R)^{-(n-1)}$ for $\mathsf R\ge1$. At fixed $\mathsf E$, this estimate gives a sufficient internal-variable dimension of order $\mathsf R\ln\mathsf R$.
\end{remark}
}

\subsection{Continuity Constraints on Adaptive Relaxation Rates}\label{subsec:ua_limits}

Relaxation rates adapted to each microstructure's spectral measure may yield Prony-series models with fewer internal variables. Such a construction must preserve continuous dependence on $(\bV,\bE)\in\cM$ to support the neural-operator approximation step in \cref{thm:ua}. Our construction uses the same relaxation rates throughout $\cM$ and lets only the Prony-series coefficients $\{\bK_i\}_{i=1}^n$ vary. The following proposition shows that adapting the interpolation interval to the spectral edge can fail to preserve continuity with respect to the microstructure. We emphasize that this does not exclude placing the adaptive relaxation rate continuously based on alternative spectral information.

\begin{proposition}[Spectral Edge Discontinuity]
\label{prop:no_spectral_stability}
There is a choice of the constants defining the admissible set $\cM$
in \cref{eq:admissible_microstructure} for which the following holds.
For each $(\bV,\bE)\in\cM$, let the largest value of the spectrum of
$\cR_{\bV,\bE}$ be $\Lambda_+(\bV,\bE)$. Then the map $(\bV,\bE)\mapsto\Lambda_+(\bV,\bE)$ is not continuous on $\cM$.
\end{proposition}

\begin{proof}[\bf Proof]
Choose $\alpha_{\rm visc}=\beta_{\rm visc}=\alpha_{\rm elas}=1$, $\beta_{\rm elas}=2$, and $\gamma_{\rm visc},\gamma_{\rm elas}$ large enough to contain the sequence below. Set $B_m\subset\bT^d$ to be the strip $\{0<y_1<1/m\}$ and let $\mathbf{1}_{B_m}$ denote the indicator function on $B_m$. Set $\bV_m\equiv \bE\equiv \bV\equiv \bI$ and $\bE_m=(1+\mathbf{1}_{B_m})\bI$, where $\bI\in\cV^4\text{ is the identity operator on }\cV^2$. Then $(\bV_m,\bE_m)\in\cM$ and $(\bV_m,\bE_m)\to(\bV,\bE)$ in $L^1_{\rm per}\text{ and }L^2_{\rm per}$. We have $\cR_{\bV,\bE}=\mathrm{Id}_{\cX}$, $\Lambda_+(\bV,\bE)=1$, where $\mathrm{Id}_{\cX}$ is the identity on $\cX$. Furthermore, $\cR_{\bV_m,\bE_m}\eta=
\eta+\cP_{\bI}(\mathbf{1}_{B_m}\eta)$ for all $\eta\in\cX$. Choose normalized strain correctors $\eta_m\in\cX$ supported in $B_m$. For example, take a nonconstant smooth function $\psi_m$ supported in $B_m$ and a nonzero $y_*\in\bR^d$; after subtracting the mean of $\psi_m y_*$ and normalizing its symmetric gradient in $L^2_{\rm per}$, one obtains such an $\eta_m$. Then $\mathbf{1}_{B_m}\eta_m=\eta_m$ and $\cP_{\bI}\eta_m=\eta_m$, so $\cR_{\bV_m,\bE_m}\eta_m=2\eta_m$. Thus $\Lambda_+(\bV_m,\bE_m)\ge 2$. Conversely, since $1\le 1+\mathbf{1}_{B_m}\le 2$, we have $(\cR_{\bV_m,\bE_m}\eta,\eta)_{L^2_{\rm per}}\le 2\|\eta\|_{L^2_{\rm per}}^2$ for all $\eta\in\cX$. Therefore, $\Lambda_+(\bV_m,\bE_m)\le 2$, and $\Lambda_+(\bV_m,\bE_m)=2$ for all $m$, whereas $\Lambda_+(\bV,\bE)=1$.
\end{proof}

\section{Numerical Examples}\label{sec:numerical_examples}

This section presents two sets of numerical experiments that complement the theory. In \cref{subsec:2d_viscoelasticity}, we study Kelvin--Voigt viscoelasticity and examine whether the numerical results corroborate our theory. In \cref{subsec:3d_crystal_plasticity}, we turn to crystal plasticity and show that the neural operator model can be applied in a nonlinear inelastic setting beyond the scope of the theory. {For consistency, all numerical datasets comprise trajectories restricted to the time interval $[0,1]$.}

\subsection{Viscoelasticity}\label{subsec:2d_viscoelasticity}

We consider learning 2D plane-stress Kelvin--Voigt viscoelasticity. This subsection first describes the dataset, model architecture, and training objective in \cref{sssec:kv2d_setup}. In \cref{sssec:kv2d_model_complexity_numerical}, we link the numerical results to our theory: the stress prediction error decreases as $n_{\xi}$ increases, increases with material contrast, and also depends on the number of training samples. \Cref{sssec:kv2d_prony_diagnostics} shows that the internal variables of the trained neural operator model admit a tensorial interpretation after a linear change of coordinates, but the fitted dynamics remain coupled and are not necessarily dissipative. Finally, \cref{sssec:kv2d_particulate_transfer} evaluates models trained on Voronoi microstructures on a separate particulate dataset and shows that they retain predictability under this input distributional shift.

\subsubsection{Setup}\label{sssec:kv2d_setup}
\begin{figure}[htp]
    \centering
\addtolength{\tabcolsep}{-5 pt}
    {\bf Sample \#1}\\[5pt]
    
    \scalebox{0.85}{
    \begin{tabular}{>{\centering\arraybackslash} m{0.05\linewidth} >{\centering\arraybackslash} m{0.15\linewidth} >{\centering\arraybackslash} m{0.15\linewidth}}
       & Elasticity & Viscosity  \\
       \raisebox{5pt}{\rotatebox[origin=c]{90}{Modulus}}& \includegraphics[width=\linewidth]{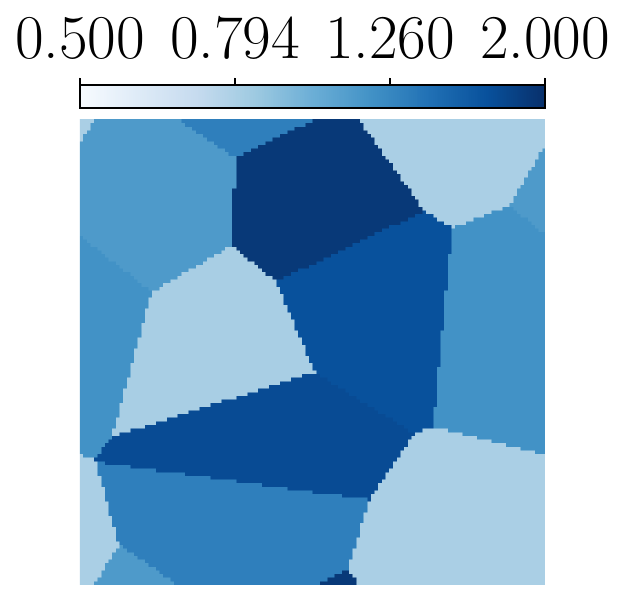}  & \includegraphics[width=\linewidth]{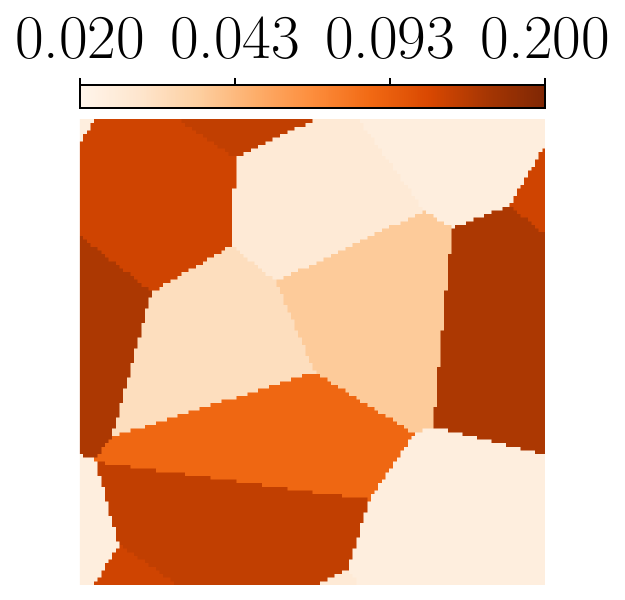} \\
       \raisebox{5pt}{\rotatebox[origin=c]{90}{\makecell{Poisson\\ratio}}}& \includegraphics[width=\linewidth]{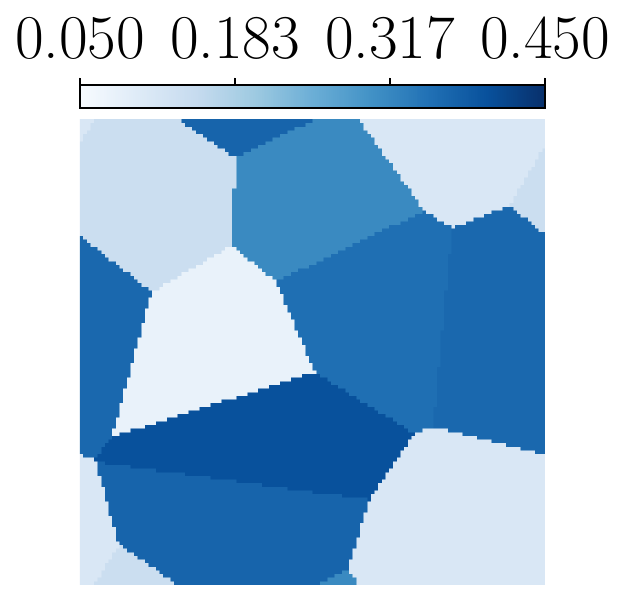} & \includegraphics[width=\linewidth]{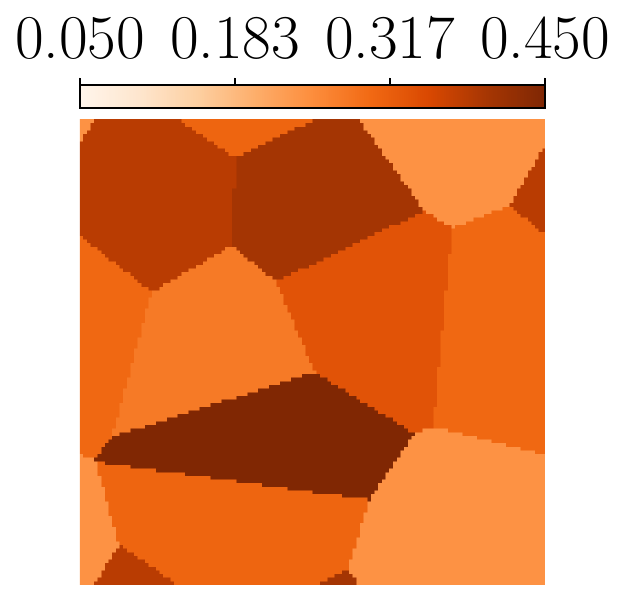}
    \end{tabular}
    }
    {
\renewcommand{\arraystretch}{0.0}
    \begin{tabular}{c c}
        \includegraphics[width=0.25\linewidth, trim={0 0.9in 0 0},clip]{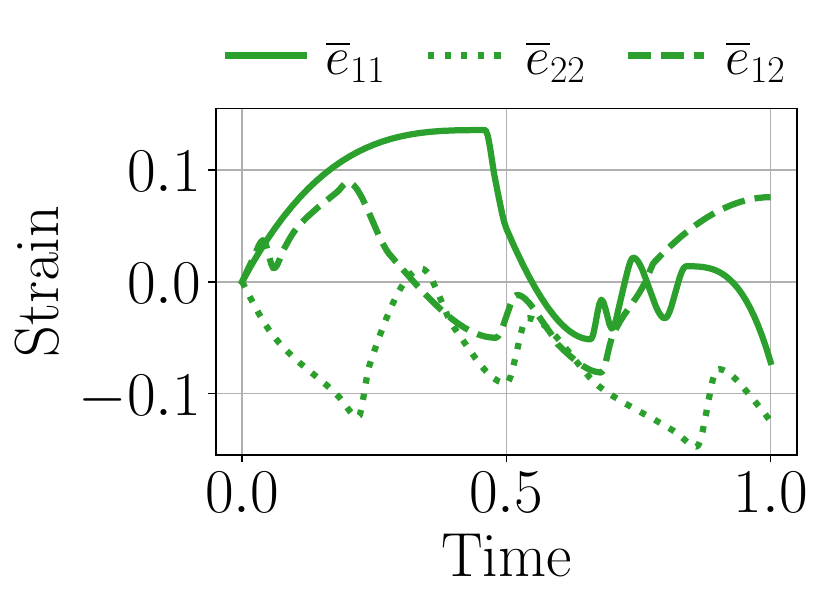} & \includegraphics[width=0.26\linewidth, trim={0 0.9in 0 0},clip]{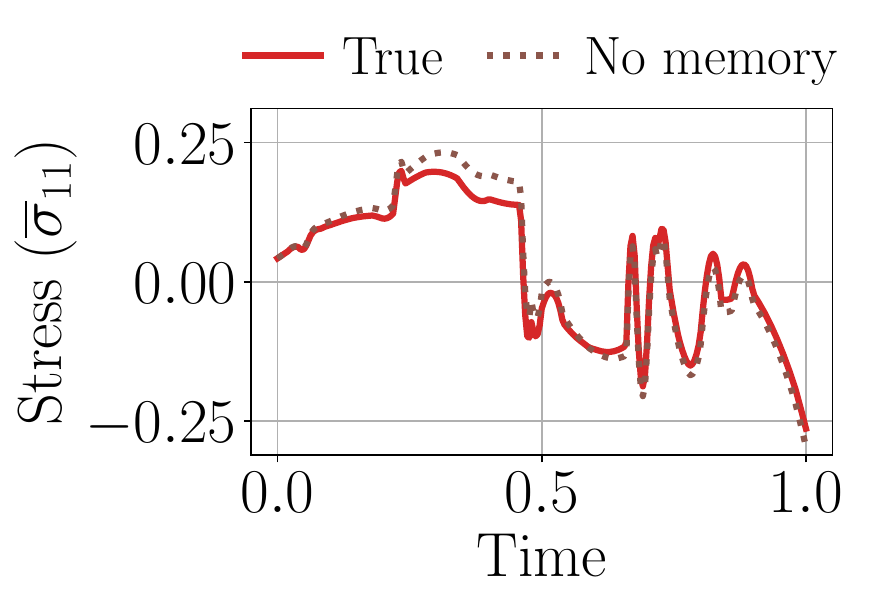}  \\
        \hspace{-0.005\linewidth}\includegraphics[width=0.267\linewidth, trim={0 0 0 0.6in},clip]{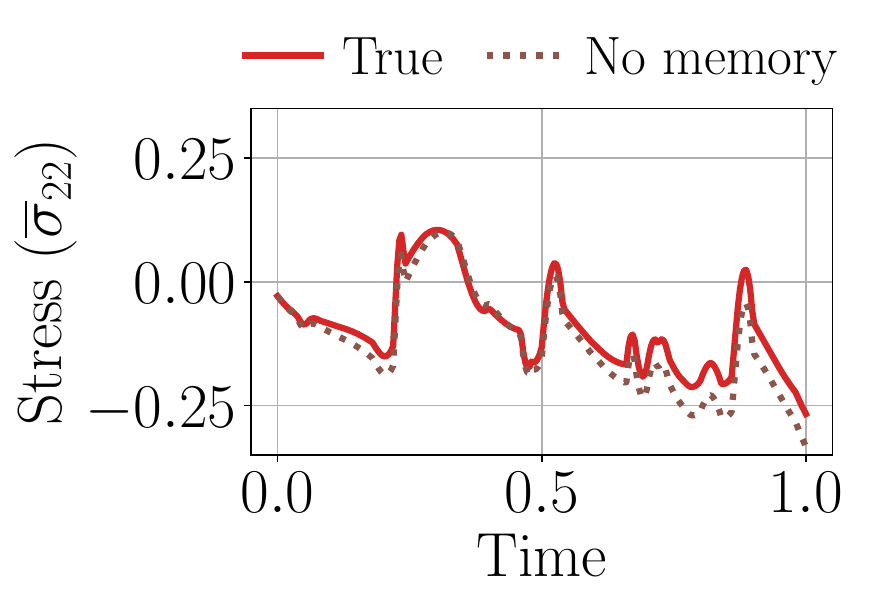} &\hspace{0.003\linewidth} \includegraphics[width=0.25\linewidth, trim={0 0 0 0.6in},clip]{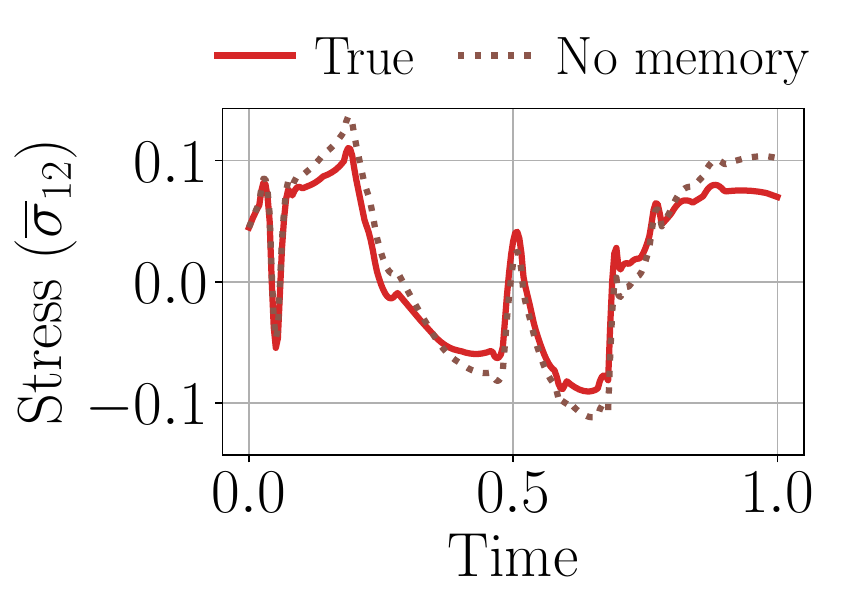} 
    \end{tabular}
\renewcommand{\arraystretch}{1.0}
}\\[5pt]

    {\bf Sample \#2}\\[5pt]
    
    \scalebox{0.85}{
    \begin{tabular}{>{\centering\arraybackslash} m{0.05\linewidth} >{\centering\arraybackslash} m{0.15\linewidth} >{\centering\arraybackslash} m{0.15\linewidth}}
       & Elasticity & Viscosity  \\
       \raisebox{5pt}{\rotatebox[origin=c]{90}{Modulus}}& \includegraphics[width=\linewidth]{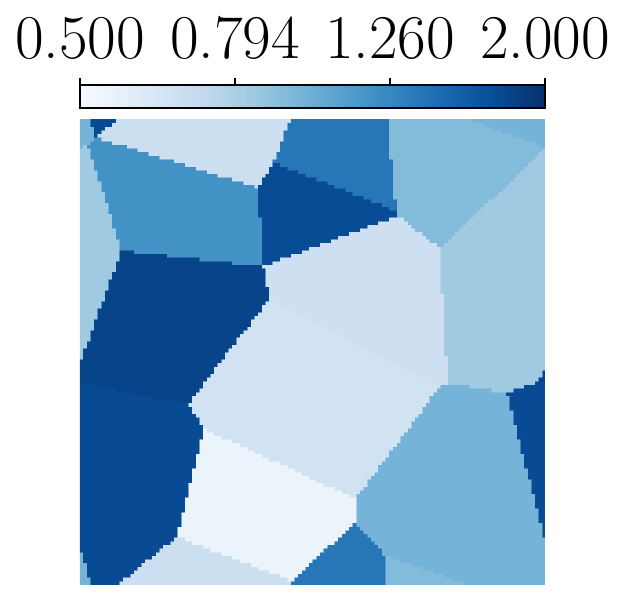}  & \includegraphics[width=\linewidth]{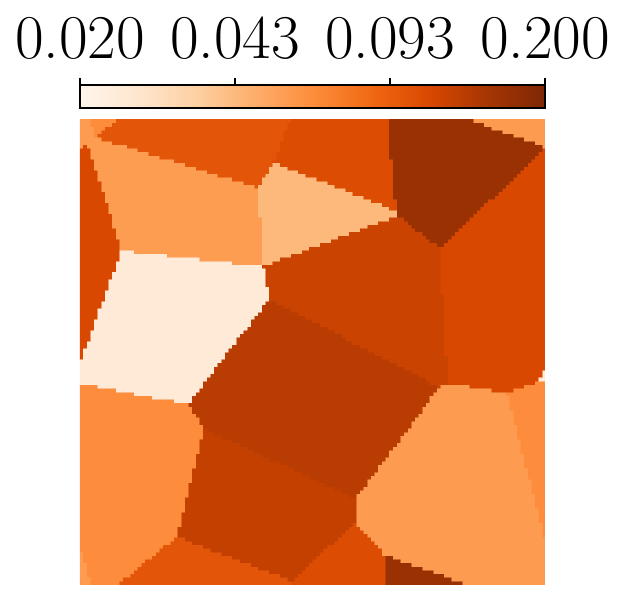} \\
       \raisebox{5pt}{\rotatebox[origin=c]{90}{\makecell{Poisson\\ratio}}}& \includegraphics[width=\linewidth]{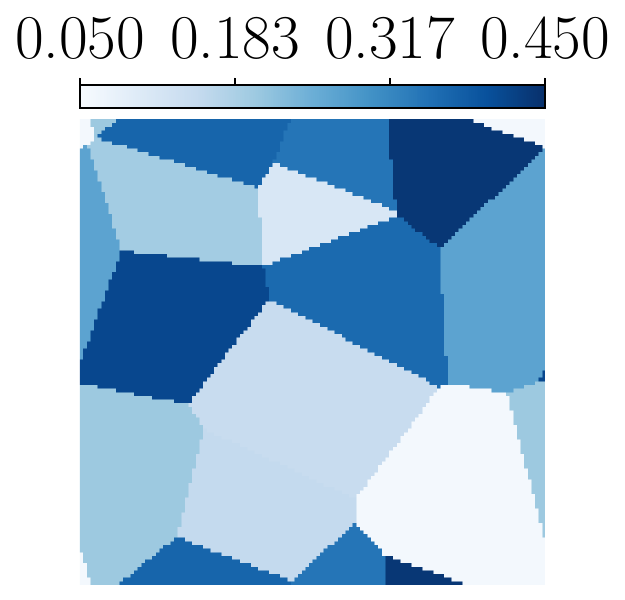} & \includegraphics[width=\linewidth]{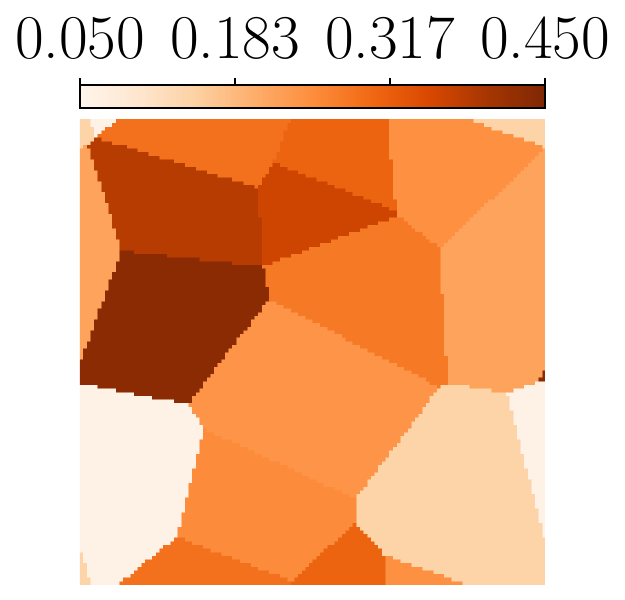}
    \end{tabular}
    }
    {
\renewcommand{\arraystretch}{0.0}
    \begin{tabular}{c c c}
        \includegraphics[width=0.25\linewidth, trim={0 0.9in 0 0},clip]{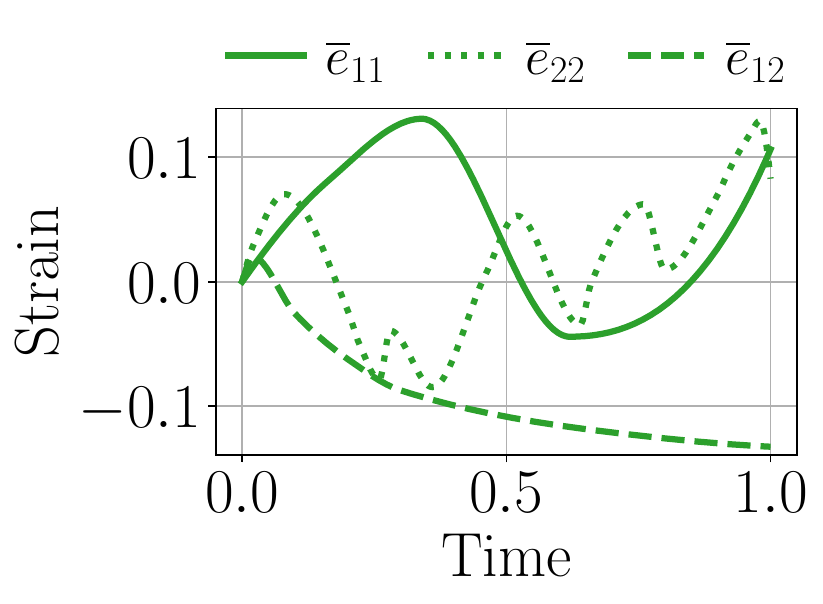} & \includegraphics[width=0.26\linewidth, trim={0 0.9in 0 0},clip]{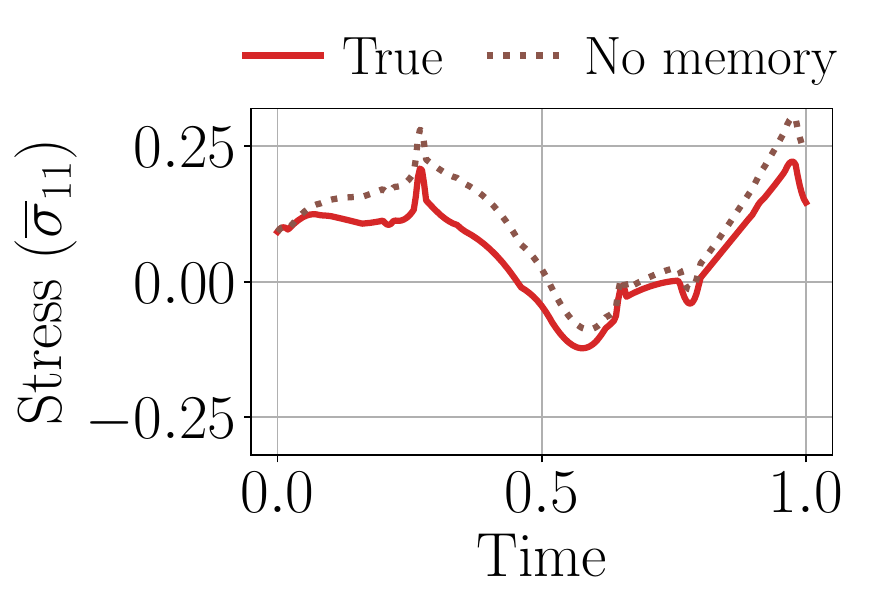}  \\
        \hspace{-0.005\linewidth}\includegraphics[width=0.265\linewidth, trim={0 0 0 0.6in},clip]{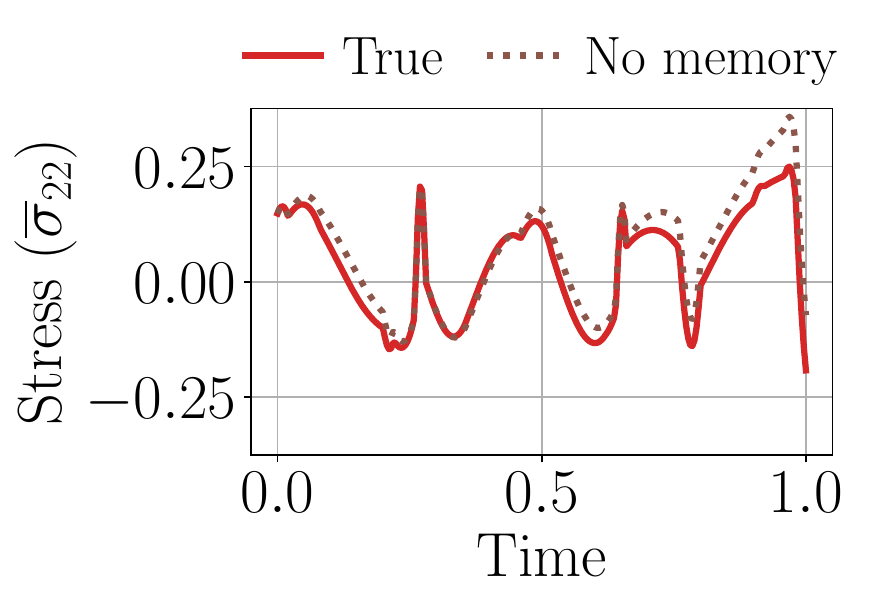} & \hspace{0.01\linewidth}\includegraphics[width=0.255\linewidth, trim={0 0 0 0.6in},clip]{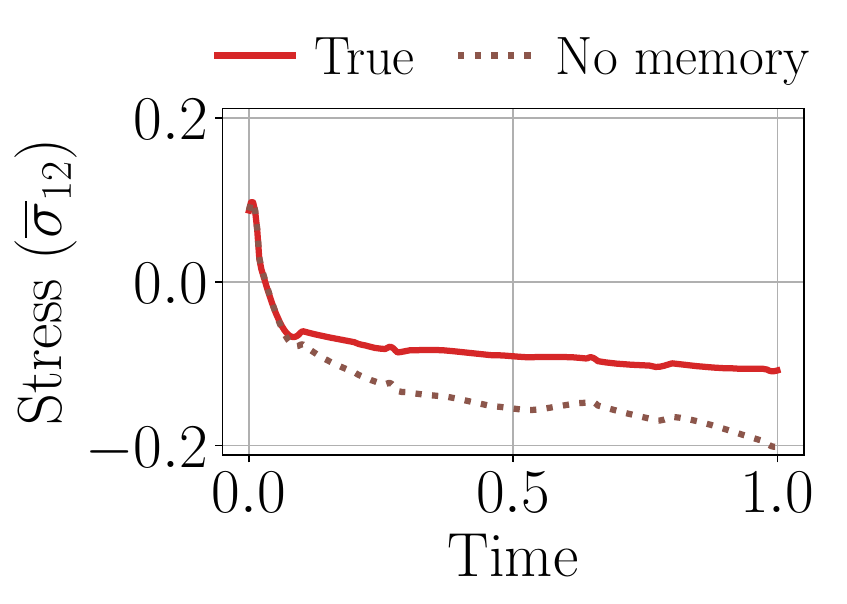} 
    \end{tabular}
\renewcommand{\arraystretch}{1.0}
}
\addtolength{\tabcolsep}{5 pt}
    \caption{Representative samples from the viscoelasticity dataset. For each sample, the left panels show the elastic and viscous microstructural fields, with the modulus values plotted on a logarithmic scale. The right panels show the prescribed macroscopic strain trajectory, the resulting homogenized stress components, and the stress predicted by the memoryless model ($\bK=0$).}
    \label{fig:kv2d_data_visual}
\end{figure}

The viscous and elastic tensors are parameterized by Young's modulus $E$ and Poisson ratio $\nu$, i.e., $\bV(y)=\bA(E_{\bV}(y),\nu_{\bV}(y))$ and $\bE(y)=\bA(E_{\bE}(y),\nu_{\bE}(y))$, where $\bA(E,\nu)e = Ee/(1+\nu) + (E\nu{\operatorname{tr}}(e)/(1-\nu^2)){I}$ for $e\in\cV^2$, where {$I$} is the second-order identity tensor. Thus, the microstructural input consists of the four scalar fields $M=(E_{\bV},\nu_{\bV},E_{\bE},\nu_{\bE})$, so that $d_M=4$. The baseline microstructures are periodic Voronoi tessellations with randomly sampled 5--15 seeds, whose regions carry independently sampled modulus and Poisson-ratio values. Details of the Voronoi dataset are given in \cref{app:kv2d_voronoi_dataset,app:kv2d_modulus_ranges}. For each microstructure, we sample a macroscopic strain trajectory and compute the corresponding averaged stress from numerical solutions of the unit-cell problem. We provide details on strain-trajectory sampling and the numerical solver in \cref{app:kv2d_strain_sampling} and \cref{app:kv2d_cell_solver}, respectively.  In \cref{fig:kv2d_data_visual}, we visualize samples from the Voronoi dataset.

{Using a Fourier neural mapping architecture for both the stress and kinetic response mappings in \cref{eq:F_and_G}, we} solve the internal-variable dynamics with the forward Euler method. Let $\overline\sigma_{\mathrm{NO}}^{(j)}$ denote the predicted stress trajectory for the $j$th training sample. We train the model by minimizing the loss function given by
\begin{equation}\label{eq:loss_function}
    \mathrm{Loss}
    =
    \frac{1}{n_{\mathrm{train}}}
    \sum_{j=1}^{n_{\mathrm{train}}}
    \left(
        {\operatorname{err}_{\rm rel}}\bigl(\overline\sigma_{\mathrm{NO}}^{(j)},\overline\sigma^{(j)}\bigr)^2
        +
        \alpha\left\|\cG_{\mathrm{NO}}\bigl(0,0;M^{(j)}\bigr)\right\|^2
    \right),
\end{equation}
where the relative stress prediction error is defined by
\begin{equation}\label{eq:relative_stress_error}
    {\operatorname{err}_{\rm rel}}(\overline\sigma_{\mathrm{NO}},\overline\sigma)
    =
    \left(
        \frac{\int_0^T \|\overline\sigma_{\mathrm{NO}}(t)-\overline\sigma(t)\|_F^2\,\mathrm{d}t}{\int_0^T \|\overline\sigma(t)\|_F^2\,\mathrm{d}t}
    \right)^{1/2}.
\end{equation}
Here $T=1$. The penalty term enforces that the kinetic response vanishes at zero strain rate and zero internal state; see \cite[Appendix E]{bhattacharya2025learning} for a discussion of this penalty term. We tune the weighting constant $\alpha$ along with the other architectural hyperparameters. Details of the architecture and training are given in \cref{app:kv2d_training}.

\subsubsection{Internal-Variable Complexity and Sample Complexity}\label{sssec:kv2d_model_complexity_numerical}
We examine how test error depends on internal-variable dimension, material contrast, and the number of training samples. The first two quantities are related to the internal-variable complexity estimate in \cref{thm:sufficient_dimension}. The sample-complexity study is empirical, since the approximation theorem does not address finite-sample effects. Our theory relies on an existence result for the internal-variable neural-operator model, while the experiments use a prescribed architecture, finite training samples, and non-convex stochastic optimization. The empirical rates reported below therefore measure the trained neural operator model's behavior under these additional constraints.

\begin{figure}[htbp]
    \centering
\addtolength{\tabcolsep}{-5 pt}
    \begin{tabular}{c c c}
    \hspace{0.05\linewidth}\makecell{\small Internal-variable\\ complexity} & \hspace{0.05\linewidth}\makecell{\small  Material-contrast\\ complexity} & \hspace{0.05\linewidth}{\small Sample complexity} \\
       \includegraphics[width=0.325\linewidth]{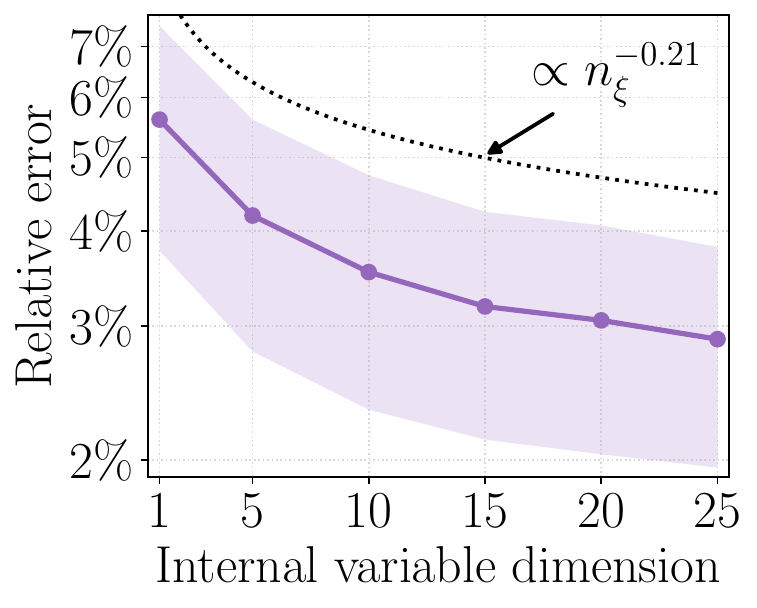}  &  \includegraphics[width=0.325\linewidth]{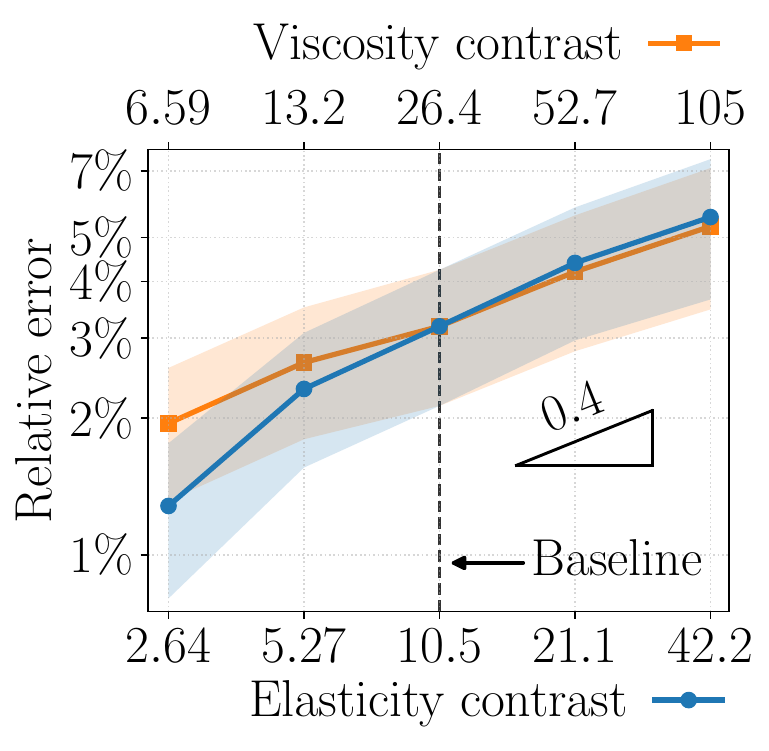} & \includegraphics[width=0.325\linewidth]{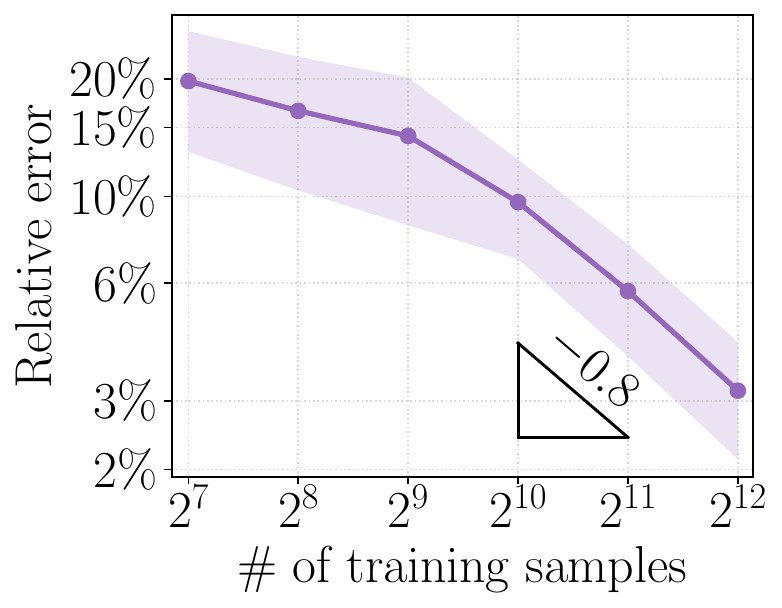}
    \end{tabular}
\addtolength{\tabcolsep}{5 pt}
    \caption{Internal-variable and sample-complexity trends for learning viscoelasticity on the Voronoi dataset. (\emph{Left}) Relative stress error versus internal-variable dimension $n_{\xi}$ for testing samples in the baseline dataset. (\emph{Middle}) Relative stress error versus material contrast at $n_{\xi}=15$, varying one material contrast at a time from the baseline. (\emph{Right}) Relative stress error versus the number of training samples for $n_{\xi}=15$. Markers show the mean; shaded bands show one standard deviation.}
    \label{fig:ve_model_complexity}
    \label{fig:test_error_vs_internal_variables}
    \label{fig:ve_complexity}
\end{figure}

The left panel of \cref{fig:ve_model_complexity} shows relative test error versus internal-variable dimension for a model trained on 4096 samples from the baseline dataset, where the elastic and viscous moduli in each Voronoi region are sampled from $[0.5,2.0]$ and $[0.02,0.2]$, respectively. The memoryless model with $\bK\equiv0$ has about $25\%$ relative error, showing that memory effects {should be accounted for in this experimental setup}. For the trained neural operator models, the mean error decreases from $5.6\%$ at $n_{\xi}=1$ to $2.9\%$ at $n_{\xi}=25$, with an empirical algebraic decay exponent of $-0.21$. We show a median-error prediction example for the baseline Voronoi dataset at $n_{\xi}=25$ in \cref{fig:kv2d_prediction_median}.

We next test the effect of material contrast at internal-variable dimension $n_{\xi}=15$. Starting from the baseline dataset, we generate new datasets by varying one elasticity or viscosity contrast at a time while keeping the other unchanged, and train a neural operator model using 4096 samples from each dataset. The corresponding elastic and viscous modulus ranges are listed in \cref{app:kv2d_modulus_ranges}. The middle panel of \cref{fig:ve_model_complexity} shows that the mean test error increases monotonically in both sweeps, with an empirical algebraic exponent of about $0.4$.

The number of training samples required to achieve a given prediction accuracy is also practically important, since each sample requires solving a unit-cell problem. We additionally report a sample-complexity study at $n_{\xi}=15$, varying the training set size from $128$ to $4096$ while keeping the architecture and training procedure unchanged. As shown in the right panel of \cref{fig:ve_complexity}, the mean relative stress error decreases from $19.7\%$ at $128$ samples to $3.2\%$ at $4096$ samples, with an algebraic decay exponent of $-0.8$ over the range $1024$--$4096$ samples. Over this range, the decay with training samples is steeper than the dependence on the internal-variable dimension in the left panel of \cref{fig:ve_complexity}, indicating that finite-sample effects remain dominant in the studied regime.

\subsubsection{Connection with a Prony-series Internal-Variable Model}\label{sssec:kv2d_prony_diagnostics}

\begin{figure}
    \centering
\addtolength{\tabcolsep}{-5 pt}
\begin{minipage}{0.4\linewidth}
\centering
        \hspace{0.05\linewidth}\makecell{The distribution of eigenvalues \\ across  64 microstructures}\\
        \includegraphics[width=\linewidth]{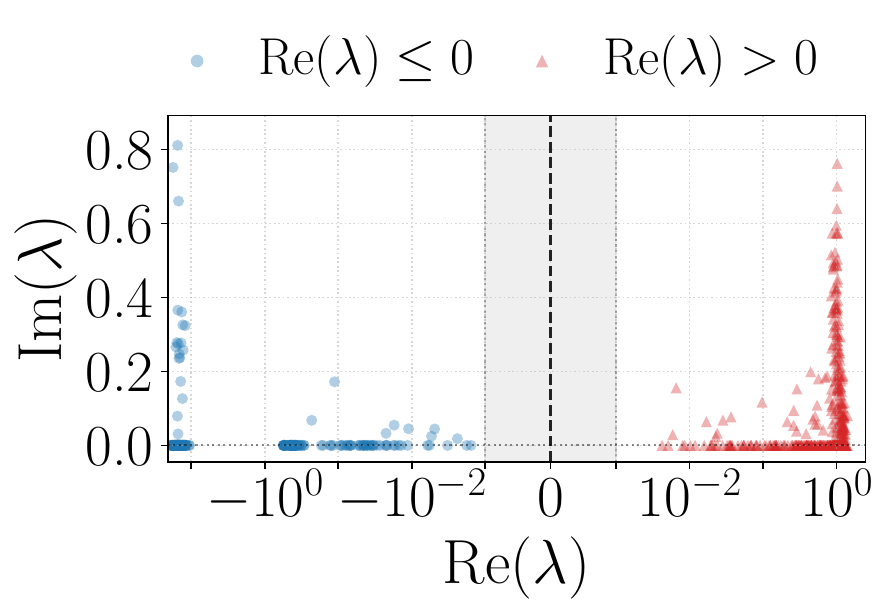}
\end{minipage}
\begin{minipage}{0.58\linewidth}
    \begin{tabular}{c c}
        \multicolumn{2}{c}{\hspace{0.075\linewidth}\makecell{Eigenvalues for 6 microstructures\\ shown separately}}\\
        \includegraphics[width=0.48\linewidth, trim={0 0.973in 0 0},clip]{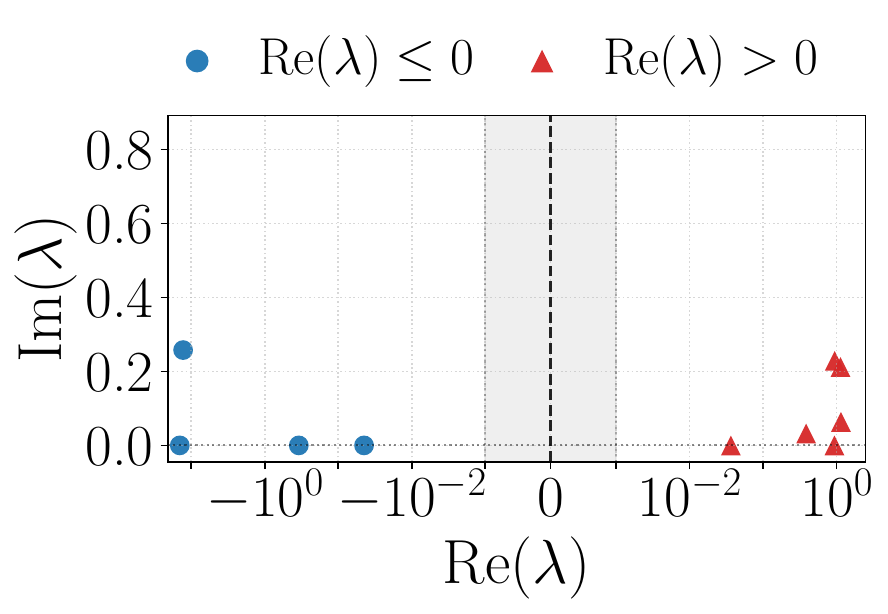} & \includegraphics[width=0.4\linewidth, trim={1.05in 0.973in 0 0.76in},clip]{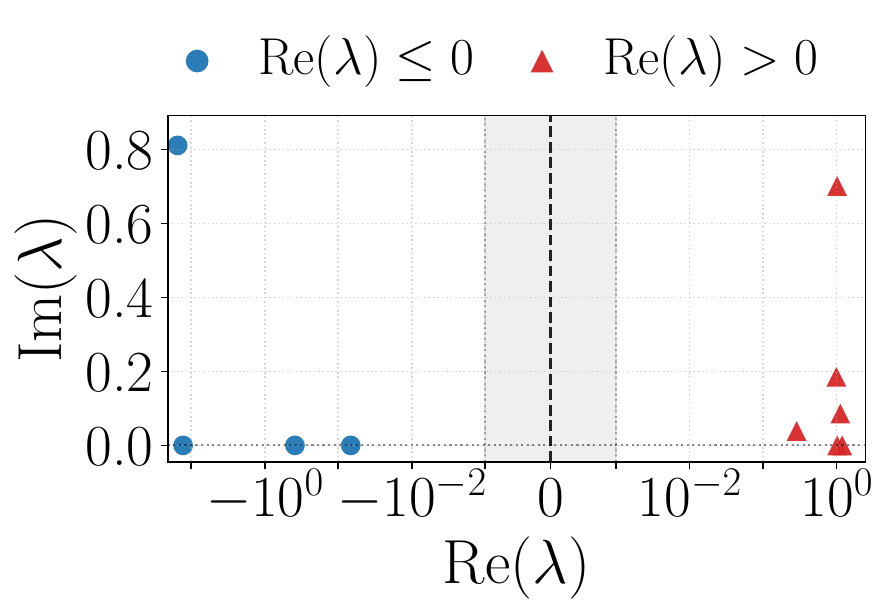}  \\
        \includegraphics[width=0.48\linewidth, trim={0 0.973in 0 0.76in},clip]{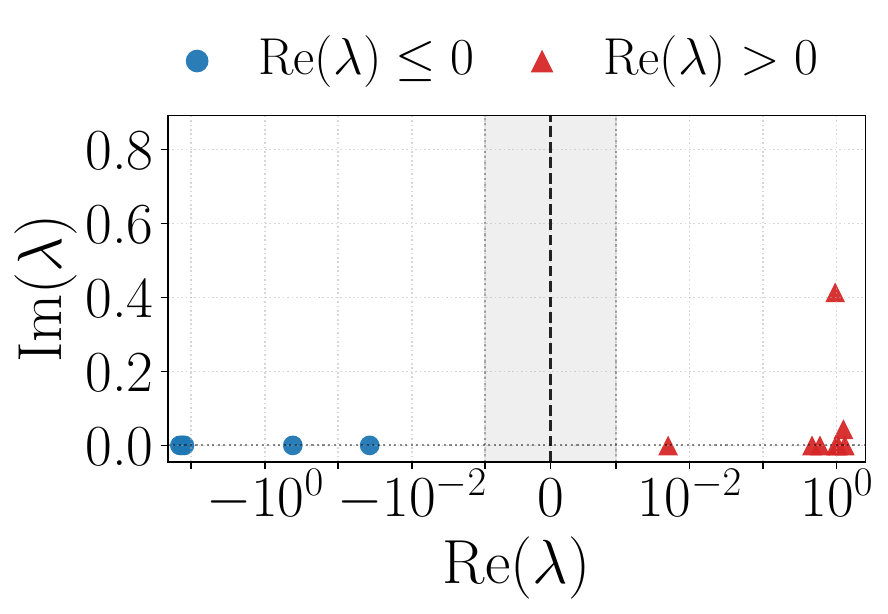} & \includegraphics[width=0.4\linewidth, trim={1.05in 0.973in 0 0.76in},clip]{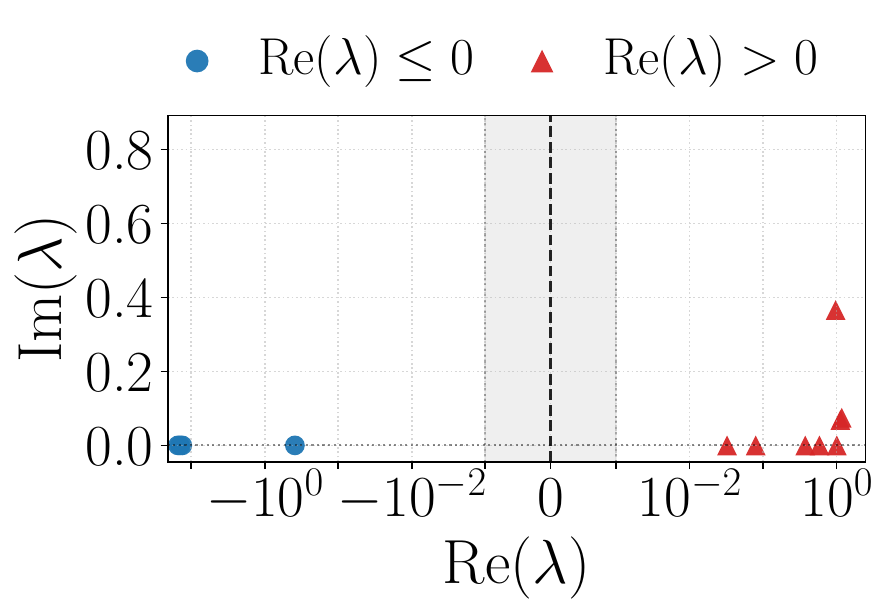} \\
        \includegraphics[width=0.48\linewidth, trim={0 0 0 0.76in},clip]{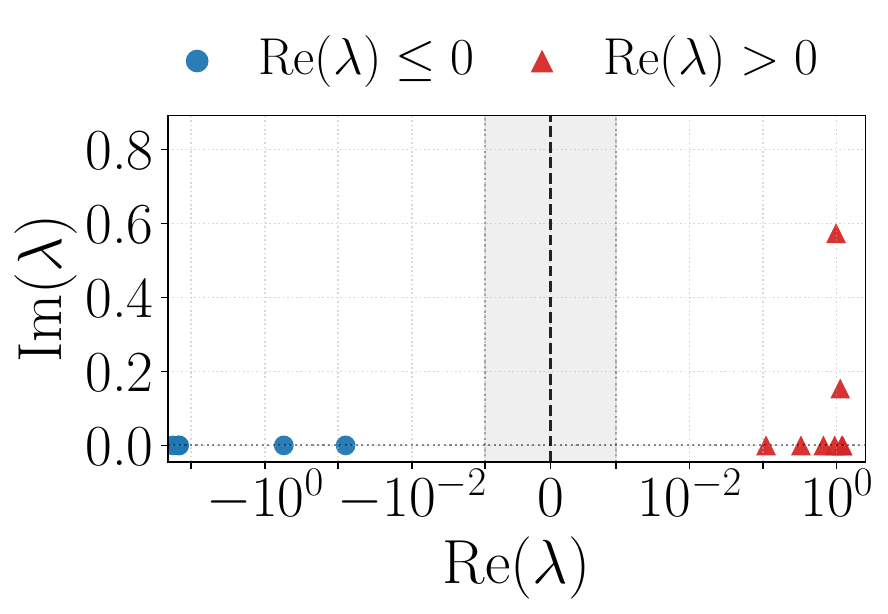} & \includegraphics[width=0.4\linewidth, trim={1.05in 0 0 0.76in},clip]{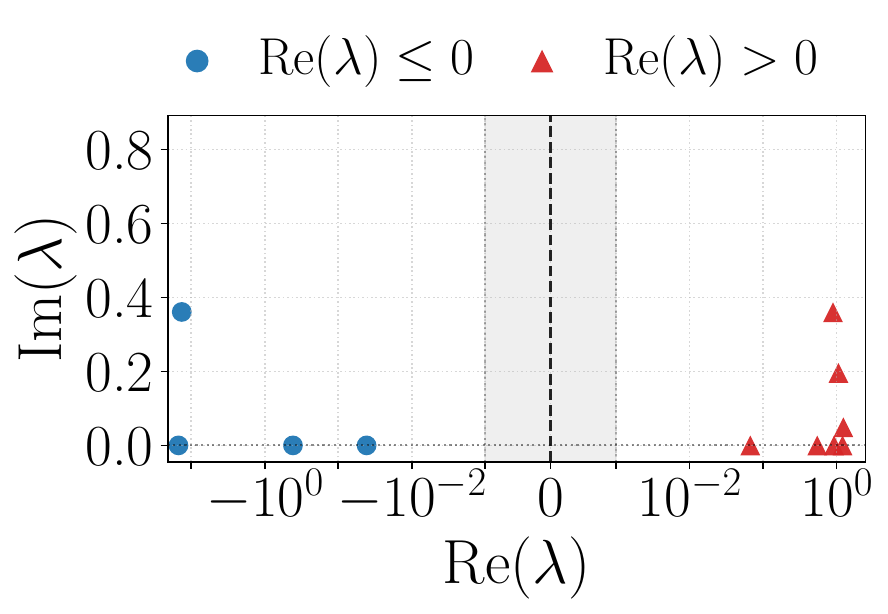} \\
    \end{tabular}
\end{minipage}
    
\addtolength{\tabcolsep}{5 pt}
    \caption{Visualizations of the eigenvalues ($\lambda$) of the {fitted internal-variable block operator} $\bL_{\rm NO}$ across microstructure samples. The {operator} $\bL_{\rm NO}$ is defined by its $(i,j)$-block $\bL_{\rm NO}^{(ij)}$ as in \cref{eqns:learned_fitted_dynamics}. (\emph{left}) Eigenvalue distribution across 64 microstructure samples. (\emph{right}) The eigenvalues obtained from 6 different microstructure samples are shown separately. The shaded area uses a linear scale for the real part of the eigenvalue, whereas the unshaded area uses a log scale.}
    \label{fig:kv2d_iv_diagnostics}
    \label{fig:eigenvalues_iv_dynamics}
\end{figure}

The approximation theory constructs the neural operator model in \cref{eq:learned_model_full} through an intermediate Prony-series model of the form \cref{eqns:prony_internal_varaible_model}. We examine whether the trained neural operator model exhibits this structure. For the trained model with $n_{\xi}=15$, we fit an affine approximation to its stress and kinetic responses $(\cF_{\rm NO}, \cG_{\rm NO})$ as functions of strain, strain rate, and internal-variable inputs. We perform the fit across 64 microstructures, with 1024 samples of macroscopic strain trajectories per microstructure, and achieve an average coefficient of determination $R^2$ of $99\%$. {In two dimensions, a symmetric tensor has $d_{\cV^2}=3$ independent components. An invertible linear change of coordinates therefore expresses the $15$ scalar internal variables as $n_\xi/d_{\cV^2}=5$ tensorial variables. The fitted linear part takes the form}
\begin{subequations}\label{eqns:learned_fitted_dynamics}
\begin{align}
    \overline{\sigma}_{\rm NO}(t)&=\bV^{\rm hom}_{\rm NO}\dot{\overline{e}}(t) + \bE^{\rm hom}_{\rm NO}\overline{e}(t)+\sum_{i=1}^{5} q_{\rm NO}^{(i)}(t),\\
    \dot{q}_{\rm NO}^{(i)}(t)
    &=\bK_{\rm NO}^{(i)}\dot{\overline{e}}(t) + \sum_{j=1}^5\bL_{\rm NO}^{(ij)} q_{\rm NO}^{(j)}(t),
    \quad q_{\rm NO}^{(i)}(0)=0,
\end{align}
\end{subequations}
where $q_{\rm NO}^{(i)}:[0,1]\to\cV^2$ are tensorial internal variables of the trained neural operator model and {$\bV^{\rm hom}_{\rm NO}, \bE^{\rm hom}_{\rm NO}, \bK_{\rm NO}^{(i)},\bL_{\rm NO}^{(ij)}\in\cL(\cV^2,\cV^2)$ are fitted, microstructure-dependent linear maps}, with $i,j\in\{1,\dots, 5\}$. Thus, after the coordinate transformation, these internal variables admit a stress-like tensorial interpretation. Their fitted dynamics, however, remain coupled and need not be dissipative, unlike the decoupled, dissipative dynamics in \cref{eqns:prony_internal_varaible_model}. In \cref{fig:kv2d_iv_diagnostics}, we visualize the eigenvalues of the {block operator} $\bL_{\rm NO}$ across the sampled microstructures and for 6 microstructures separately. The $(i,j)$ block of this {operator} is $\bL_{\rm NO}^{(ij)}$ in \cref{eqns:learned_fitted_dynamics}. The positive real parts show that some modes are not dissipative, while the nonzero imaginary parts indicate coupling-induced oscillations. In general, the negative real eigenvalues are separated on a logarithmic scale for each microstructure, efficiently capturing a wide range of relaxation timescales. Details of the linear model fitting are provided in \cref{app:kv2d_linear_fit}.

\subsubsection{Testing on Particulate Microstructures}\label{sssec:kv2d_particulate_transfer}

Finally, we test whether models trained on Voronoi tessellations transfer to a different class of microstructures. To this end, we generate a dataset of particulate microstructures: randomly placed 1--20 circular inclusions that are relatively stiff and elastic, with a relatively soft, viscous background. We keep the loading-path sampling and the overall range of material properties similar to those of the baseline Voronoi dataset. Details of the particulate dataset are described in \cref{app:kv2d_particulate_dataset}.

\begin{figure}
    \centering
\addtolength{\tabcolsep}{-5 pt}

        \scalebox{0.85}{
    \begin{tabular}{>{\centering\arraybackslash} m{0.05\linewidth} >{\centering\arraybackslash} m{0.15\linewidth} >{\centering\arraybackslash} m{0.15\linewidth}}
    \multicolumn{3}{c}{\makecell{\large A sample of \\ \large the particulate microstructure}} \\
       & Elasticity & Viscosity  \\
       \raisebox{5pt}{\rotatebox[origin=c]{90}{Modulus}}& \includegraphics[width=\linewidth]{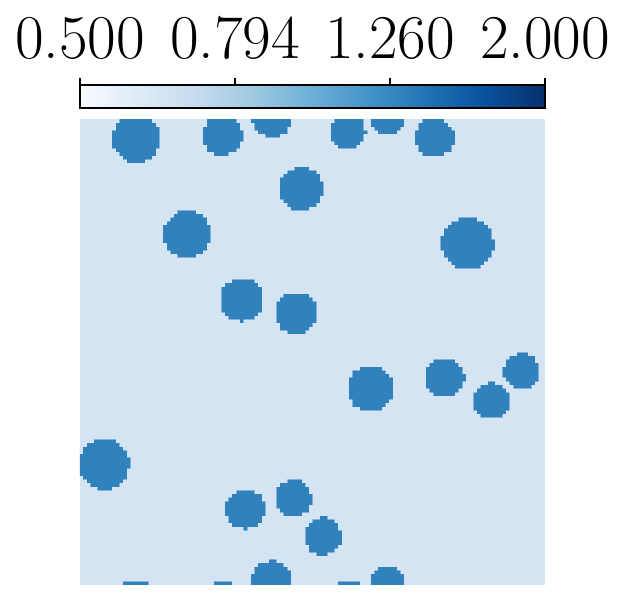}  & \includegraphics[width=\linewidth]{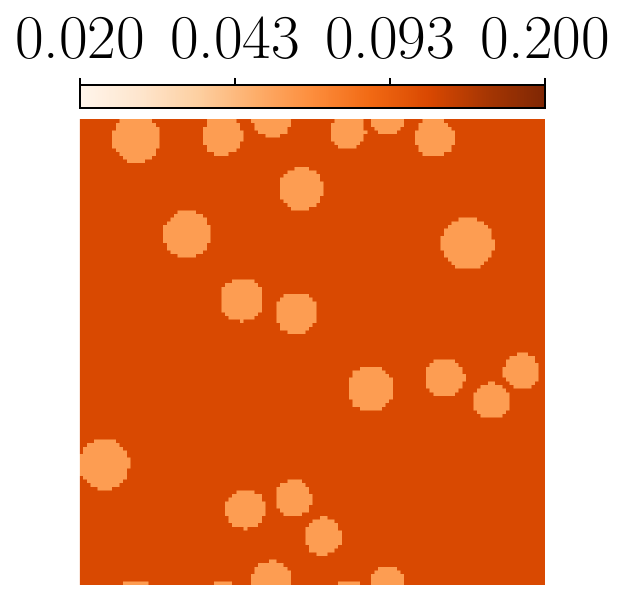} \\
       \raisebox{5pt}{\rotatebox[origin=c]{90}{\makecell{Poisson\\ratio}}}& \includegraphics[width=\linewidth]{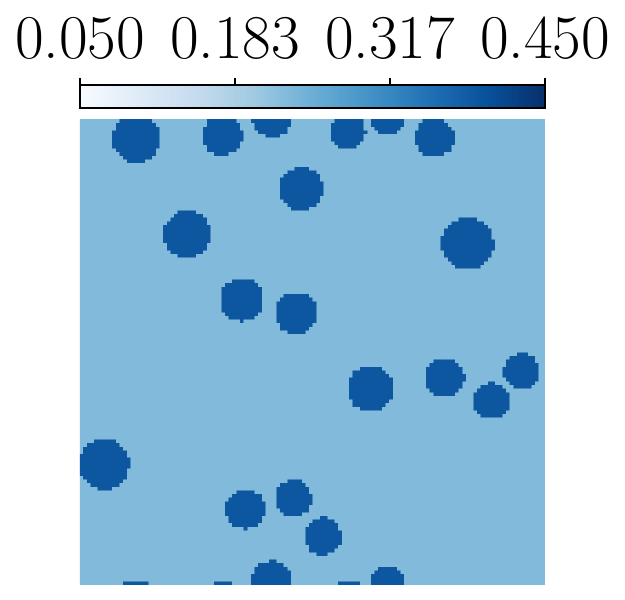} & \includegraphics[width=\linewidth]{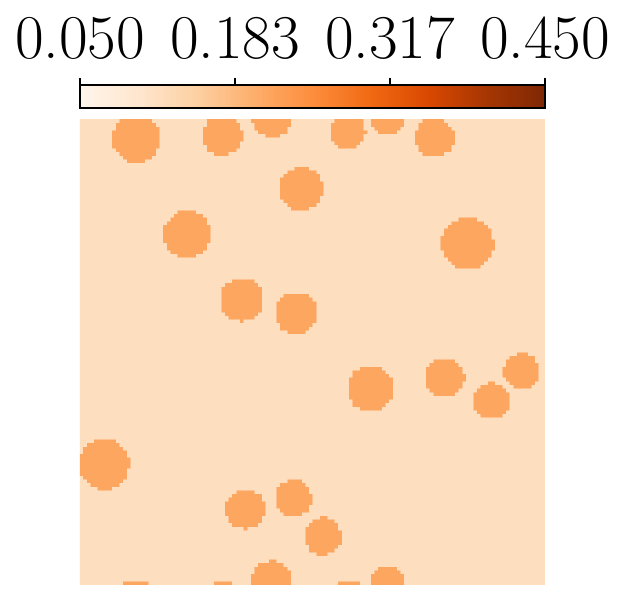}
    \end{tabular}
    }
    \begin{tabular}{c}
    \hspace{0.05\linewidth}\makecell{\small Trained on Voronoi microstructures and \\
    tested on particulate microstructures}\\
        \includegraphics[width=0.35\linewidth]{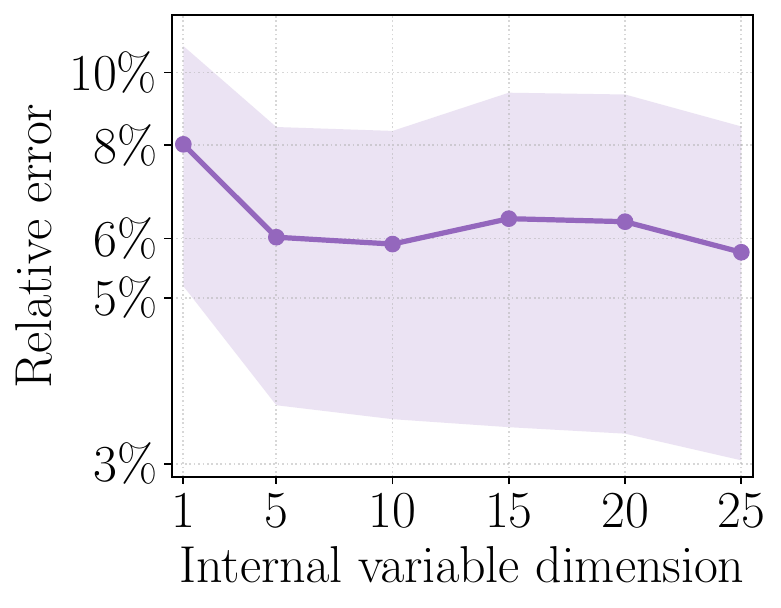}
    \end{tabular}
\addtolength{\tabcolsep}{5 pt}
    \caption{(\emph{left}) A representative sample of the particulate microstructural property field. (\emph{right}) Internal-variable complexity trend for models trained on Voronoi microstructures and evaluated on particulate microstructures. Markers show the mean; shaded bands show one standard deviation.}
    \label{fig:kv2d_particulate_iv_complexity}
\end{figure}

In \cref{fig:kv2d_particulate_iv_complexity}, we plot the test error computed on the particulate dataset as a function of the internal-variable dimension. For $n_{\xi}=25$, the mean relative stress error increases from $2.9\%$ on the in-distribution Voronoi test samples to $5.8\%$ on the particulate samples. Across all tested internal-variable dimensions, mean errors range from $5.8\%$ to $8.0\%$ on the particulate dataset, substantially lower than the memoryless model's mean error of $16.9\%$ ($\bK=0$). \Cref{fig:kv2d_particulate_median} shows a median-error prediction example for the particulate dataset at $n_{\xi}=25$. These results indicate that the trained neural operator model retains a substantial advantage over the memoryless model despite the shift in the microstructure geometry.

\subsection{Crystal Plasticity}\label{subsec:3d_crystal_plasticity}
{We next apply the neural operator model to homogenized crystal plasticity, a nonlinear inelastic setting. An explicit homogenized constitutive law exists for a one-dimensional, small-strain model without hardening \cite{liu2023learning,bhattacharya2025learning}; here we consider three-dimensional, finite-strain polycrystals, which require numerical homogenization. We describe the governing cell problem, dataset, and learning setup in \cref{sssec:cp3d_setup}, then examine empirical internal-variable and sample complexity in \cref{sssec:cp3d_model_complexity_numerical}.}

\subsubsection{Setup}\label{sssec:cp3d_setup}
We consider single-phase, face-centered cubic (FCC) Voronoi polycrystals on the unit cell $\bT^3$. {The local response combines the multiplicative decomposition, cubic Hookean elasticity, and rate-dependent slip on the twelve FCC $\{111\}\langle110\rangle$ systems \cite{asaro1983crystal,roters2010overview}}, implemented in DAMASK 3.0.2 \cite{roters2019damask}.

{
For a prescribed macroscopic deformation-gradient history $\overline F:[0,T]\to\bR^{3\times3}$ with $\det\overline F>0$ and $\overline F(0)=I$, the quasistatic cell problem seeks a periodic displacement corrector $\chi$ and the local plastic state satisfying
\begin{subequations}\label{eq:cp_cell_problem}
\begin{align}
    F(y,t) &= \overline F(t)+\nabla_y\chi(y,t), && y\in\bT^3,\; t\in[0,T],\\
    -\nabla_y \cdot P(y,t) &= 0, && y\in\bT^3,\; t\in[0,T],\\
    \int_{\bT^3}\chi(y,t)\,\mathrm{d}y &= 0, && t\in[0,T].
\end{align}
\end{subequations}
Here $F$ and $P$ are the deformation gradient and the first Piola--Kirchhoff stress on the unit cell, respectively. We specify the constitutive relation that maps $F$ to $P$ on the unit cell and the initially stress-free state in \cref{app:cp3d_constitutive}. The averaged stress $\overline P(t) = \int_{\bT^3}P(y,t)\,\mathrm{d}y$ defines the history-dependent constitutive response to be learned.
}

The dataset varies three microstructural features: grain geometry, crystallographic texture, and grainwise slip-law parameters. For each sample, the periodic Voronoi tessellation has $N_g\in\{16,\dots,32\}$ grains with randomly sampled seed locations. Grain orientations are drawn from a mixture of random, single-mode, and two-mode texture families, resulting in samples with both clustered and diffuse orientation distributions. The microstructural input $M:\bT^3\to\bR^{d_M}$ has $d_M=7$ channels: four orientation-quaternion components and the initial slip resistance, saturation slip resistance, and hardening scale. We randomly sample the three slip-law parameters for each grain around the copper reference values, while keeping the elastic constants and the remaining power-law parameters unchanged. Details of the microstructure sampling are given in \cref{app:cp3d_microstructure_sampling}.

\begin{figure}[!htp]
    \centering
    {\bf Samples \#1}\\[5pt]
    \scalebox{0.7}{
\addtolength{\tabcolsep}{-5 pt}
    \begin{tabular}{>{\centering\arraybackslash} m{0.18\linewidth} >{\centering\arraybackslash} m{0.06\linewidth}
    >{\centering\arraybackslash} m{0.16\linewidth} 
    >{\centering\arraybackslash} m{0.16\linewidth}
    >{\centering\arraybackslash} m{0.16\linewidth}}
        \raisebox{-40 pt}{\multirow{4}{*}{\includegraphics[width=\linewidth]{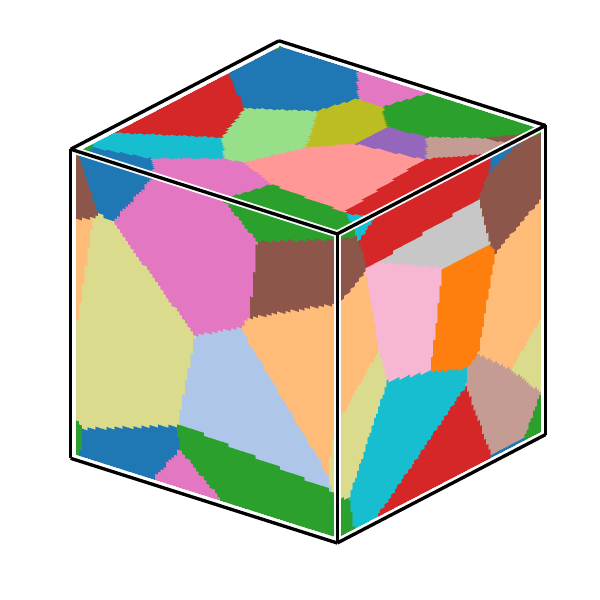}}} & \raisebox{45 pt}{\multirow{2}{*}{\rotatebox[origin=c]{90}{\makecell{Slip-law parameters \\ (MPa)}}}} &
        \includegraphics[width=\linewidth]{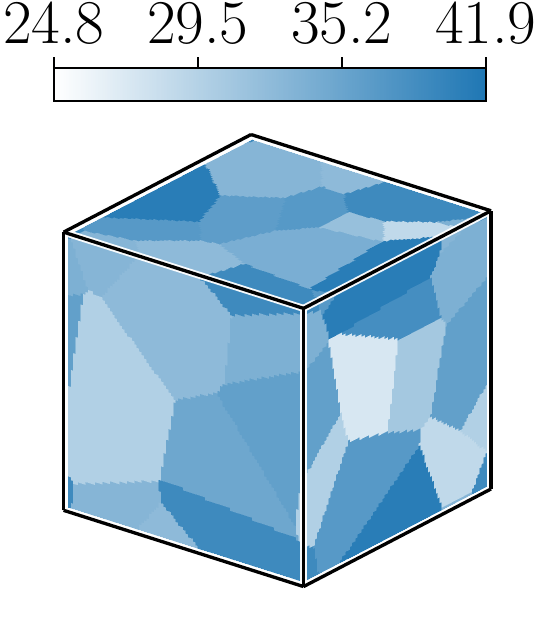}  & \includegraphics[width=\linewidth]{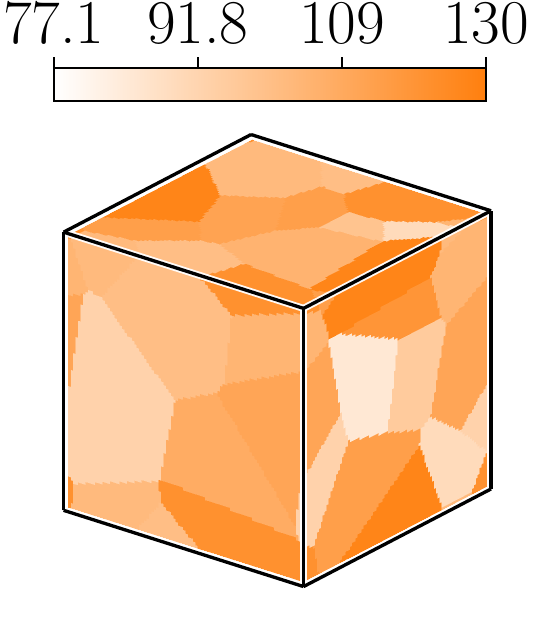} & \includegraphics[width=\linewidth]{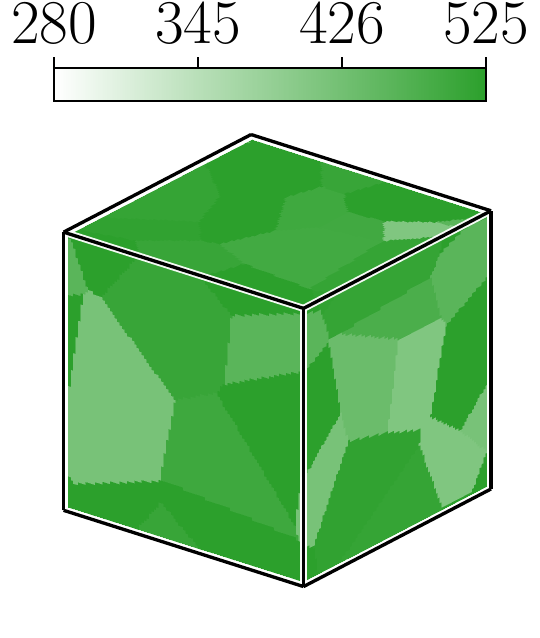} \\
         & & Initial slip resistance &Saturation slip resistance & Hardening scale \\
            & \raisebox{40 pt}{\rotatebox[origin=b]{90}{Pole figures}} & \multicolumn{3}{c}{\parbox[c]{0.48\linewidth}{\centering\includegraphics[width=0.875\linewidth]{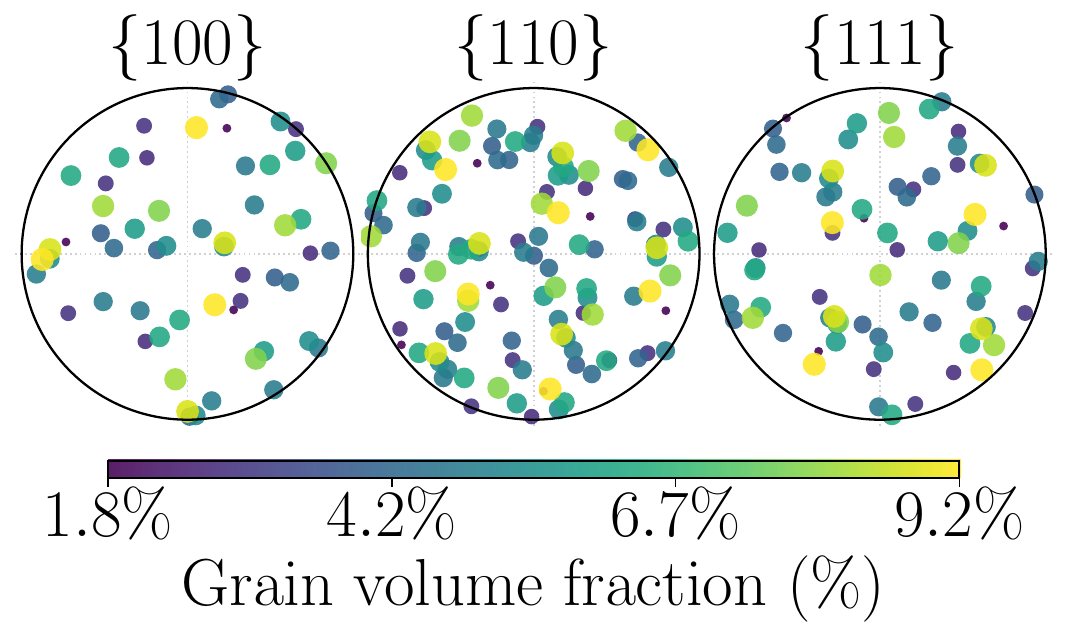}}}  \\
    \end{tabular}
    \begin{tabular}{c}
         \includegraphics[width=0.4\linewidth, trim={0 1.0in 0 0},clip]{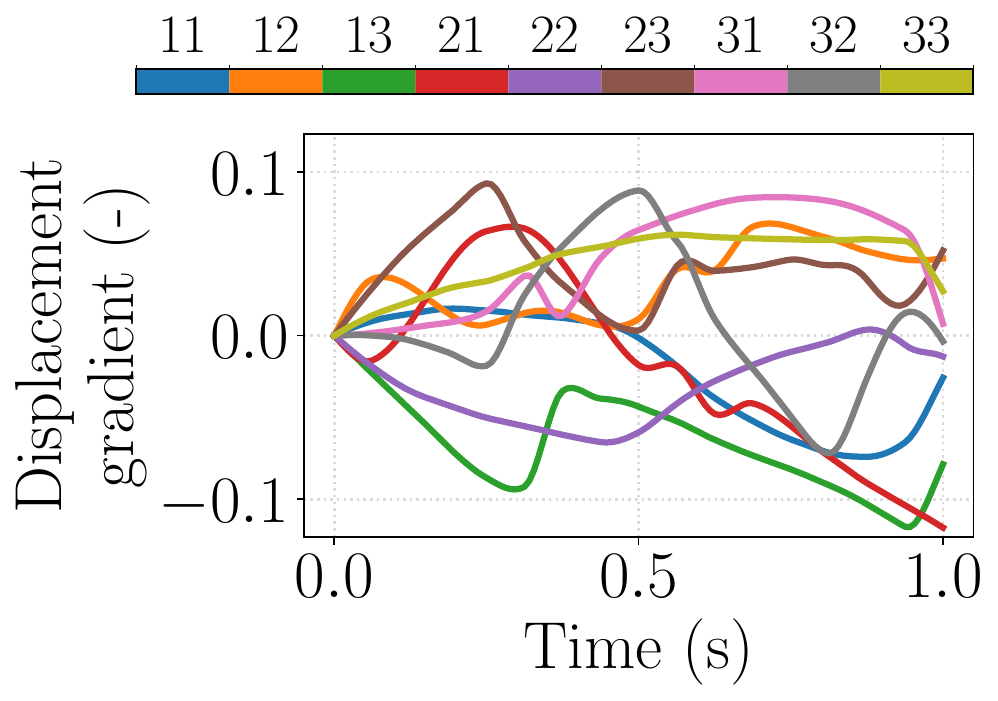}  \\
        \hspace{0.02\linewidth}\includegraphics[width=0.385\linewidth, trim={0 0 0 0.8in},clip]{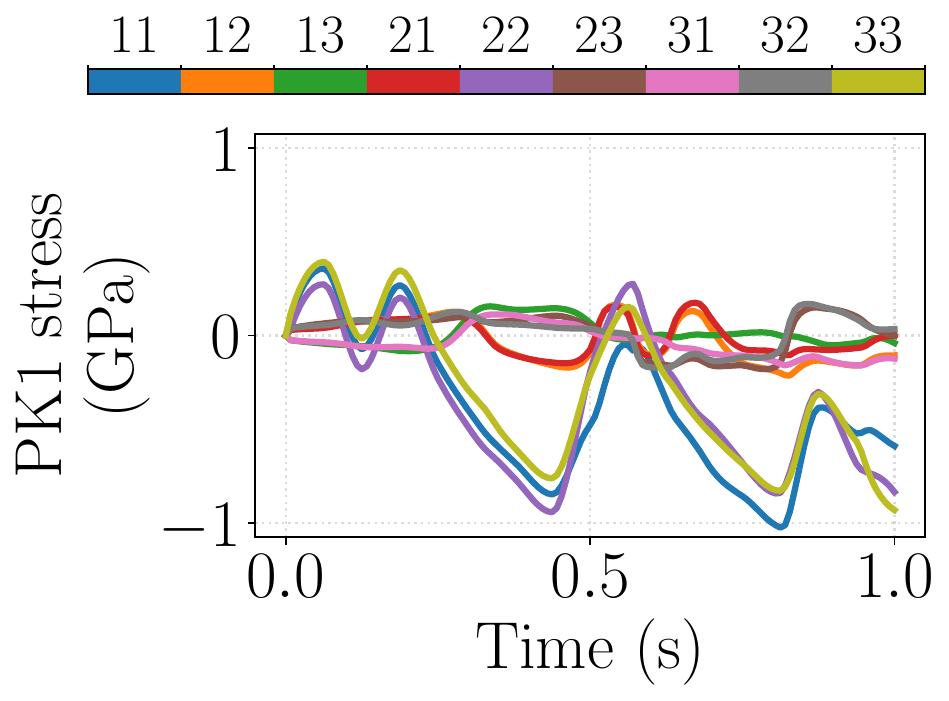}
    \end{tabular}
\addtolength{\tabcolsep}{5 pt}
    }\\[5pt]
    
    {\bf Samples \#2}\\[5pt]
    \scalebox{0.7}{
\addtolength{\tabcolsep}{-5 pt}
    \begin{tabular}{>{\centering\arraybackslash} m{0.18\linewidth} >{\centering\arraybackslash} m{0.06\linewidth}
    >{\centering\arraybackslash} m{0.16\linewidth} 
    >{\centering\arraybackslash} m{0.16\linewidth}
    >{\centering\arraybackslash} m{0.16\linewidth}}
        \raisebox{-40 pt}{\multirow{4}{*}{\includegraphics[width=\linewidth]{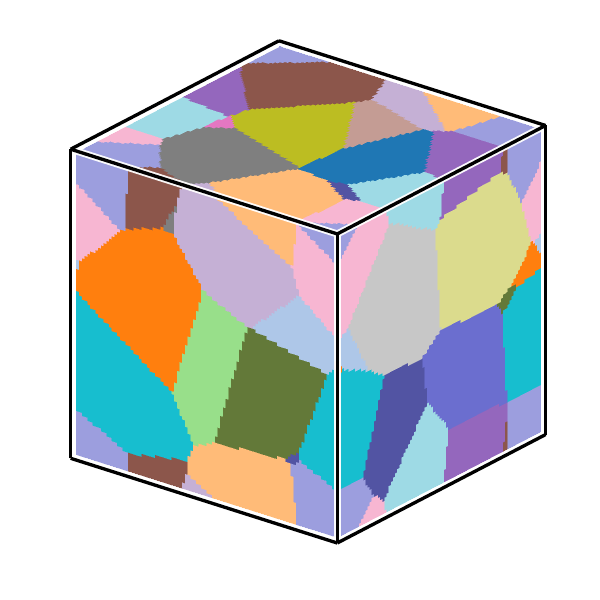}}} & \raisebox{45 pt}{\multirow{2}{*}{\rotatebox[origin=c]{90}{\makecell{Slip-law parameters \\ (MPa)}}}} &
        \includegraphics[width=\linewidth]{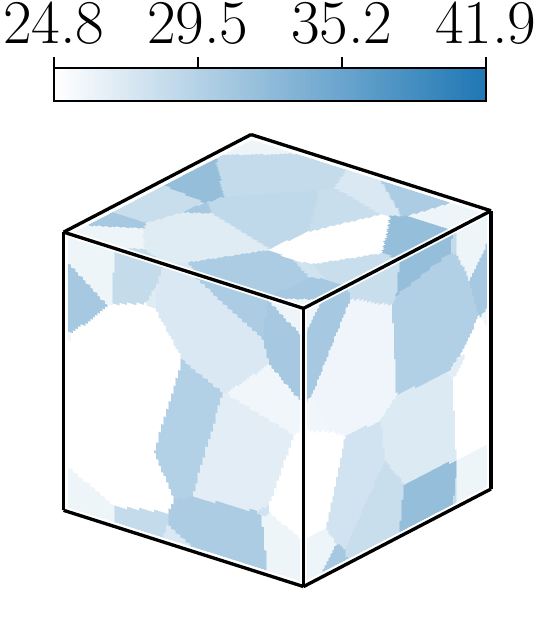}  & \includegraphics[width=\linewidth]{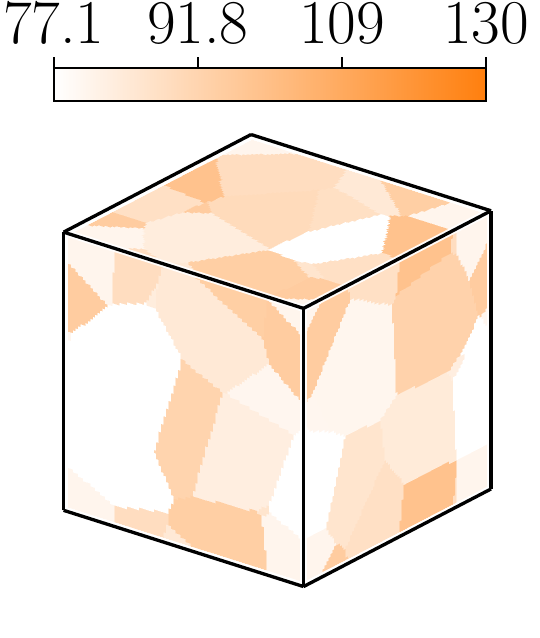} & \includegraphics[width=\linewidth]{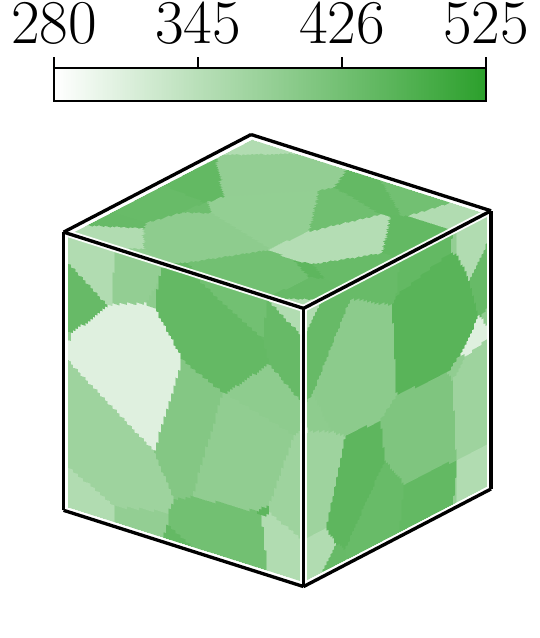} \\
         & & Initial slip resistance &Saturation slip resistance & Hardening scale \\
            & \raisebox{40 pt}{\rotatebox[origin=b]{90}{Pole figures}} & \multicolumn{3}{c}{\parbox[c]{0.48\linewidth}{\centering\includegraphics[width=0.875\linewidth]{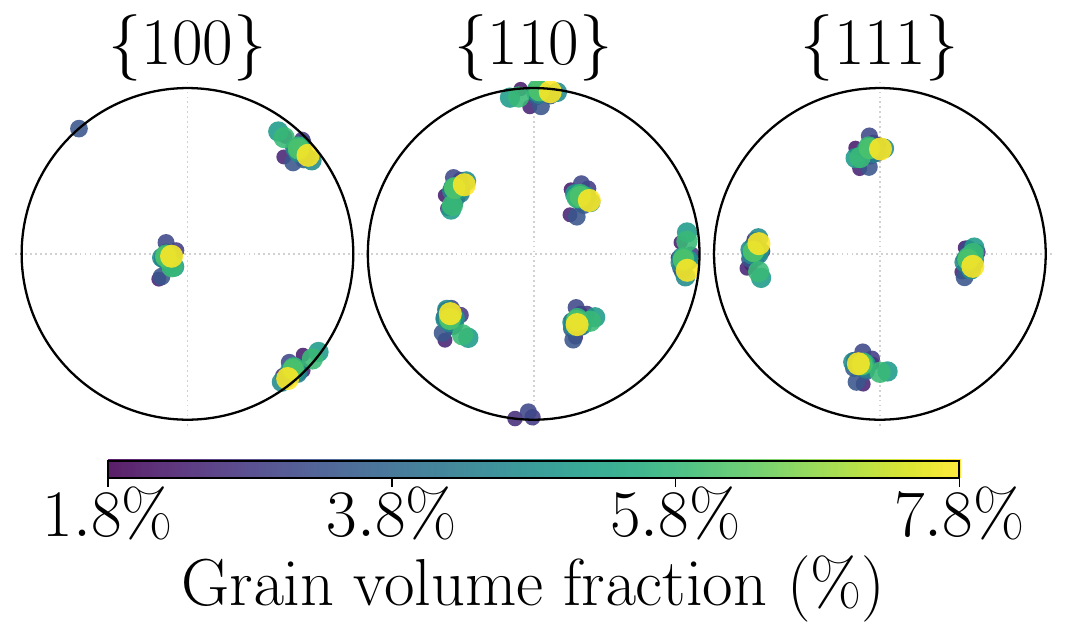}}}  \\
    \end{tabular}
    \begin{tabular}{c}
         \includegraphics[width=0.4\linewidth, trim={0 1.0in 0 0},clip]{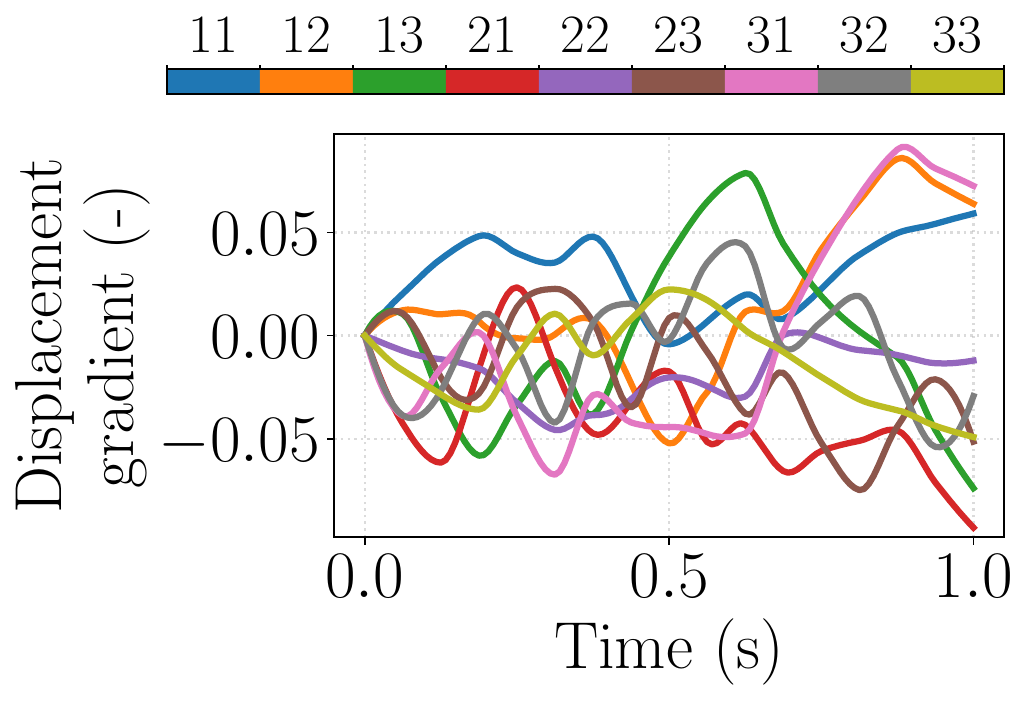}  \\
        \hspace{0.012\linewidth}\includegraphics[width=0.39\linewidth, trim={0 0 0 0.8in},clip]{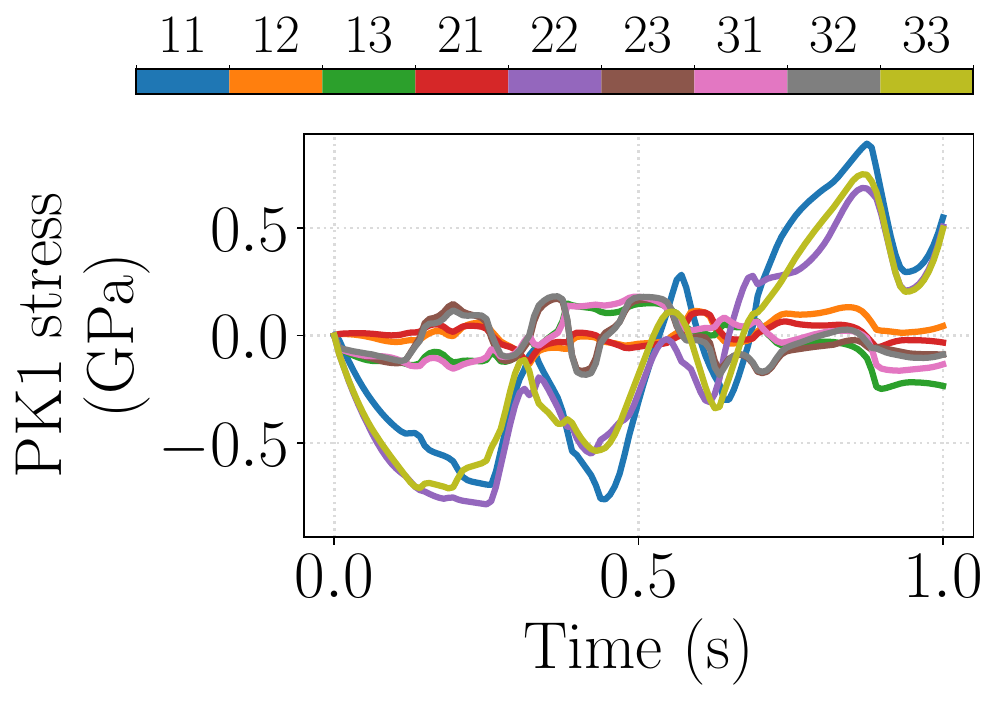}
    \end{tabular}
\addtolength{\tabcolsep}{5 pt}
    }
    \caption{Representative samples from the crystal plasticity dataset. For each sample, the left panels show the grain geometry, the grainwise slip-law parameters, and the pole figures corresponding to the crystallographic texture. The right panels show the prescribed macroscopic displacement-gradient trajectory $(\overline{F}-I)$ and the resulting averaged first Piola--Kirchhoff stress $(\overline{P})$. Sample \#1 is drawn from the random-texture family, whereas \#2 is drawn from the single-mode-texture family.}
    \label{fig:cp_samples}
\end{figure}

For each sample, we construct a loading history $\overline F$ by sampling the volumetric logarithmic strain, deviatoric logarithmic strain, and rotation vector, interpolating them in time, and combining the resulting stretch and rotation as described in \cref{app:cp3d_loading_sampling}. This construction guarantees a positive determinant and produces multiaxial, nonproportional, and sign-changing loading paths. {In the numerical implementation, we set $T=1$ and solve \cref{eq:cp_cell_problem}} using DAMASK's spectral solver on a $32\times32\times32$ voxel grid with 128 loading steps to obtain the averaged stress history $\overline P$. Representative microstructures, loading paths, and stress responses are shown in \cref{fig:cp_samples}.

We approximate this history-dependent response with the {following neural operator model, extended from \cite{liu2023learning} to include material dependence:}
\begin{align*}
    \overline{P}_{\mathrm{NO}}(t) &= \cF_{\mathrm{NO}}(\overline{F}(t), \xi(t); M),\\
    \dot{\xi}(t) &= \cG_{\mathrm{NO}}(\overline{F}(t),\xi(t); M),\quad \xi(0)=0{.}
\end{align*}
Here, we represent both the stress and kinetic responses $\cF_{\mathrm{NO}}$ and $\cG_{\mathrm{NO}}$ using Fourier neural mappings. We use the forward Euler method to solve the internal-variable dynamics. We train the model by minimizing a loss function similar to \cref{eq:loss_function}, which sums the mean relative stress-prediction error and a penalty term enforcing $\cG_{\rm NO}(I, 0;M)=0$. Details of the architecture and training are given in \cref{app:cp_training}.

\subsubsection{Internal-Variable Complexity and Sample Complexity}\label{sssec:cp3d_model_complexity_numerical}

\begin{figure}[htp]
    \centering
\addtolength{\tabcolsep}{-5 pt}
    \begin{tabular}{c c}
    \hspace{0.05\linewidth}\makecell{\small Internal-variable\\ complexity} & \hspace{0.05\linewidth}{\small Sample complexity} \\
       \includegraphics[width=0.325\linewidth]{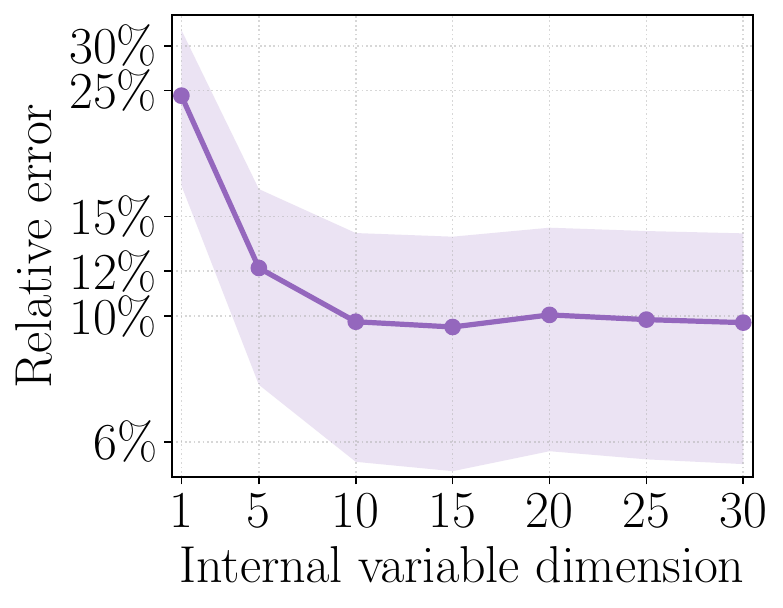}  & \includegraphics[width=0.325\linewidth]{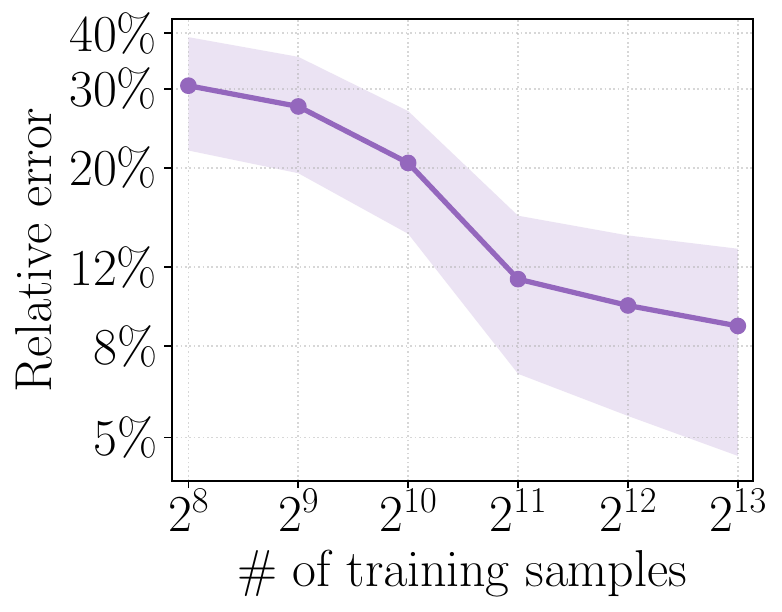}
    \end{tabular}
\addtolength{\tabcolsep}{5 pt}
    \caption{Internal-variable and sample-complexity trends for learning homogenized crystal plasticity. (\emph{Left}) Relative stress error versus internal-variable dimension $n_{\xi}$. (\emph{Right}) Relative stress error versus the training sample size for $n_{\xi}=25$. Markers show the mean; shaded bands show one standard deviation.}
    \label{fig:cp_complexity}
\end{figure}

The left panel of \cref{fig:cp_complexity} shows relative test error versus internal-variable dimension for a model trained on 4096 samples. For the trained neural operator models, the mean error decreases from $25\%$ at $n_{\xi}=1$ to $9.5\%$ at $n_{\xi}=15$. We observe a plateau in testing error after reaching $n_{\xi}=10$. A median-error prediction example for $n_{\xi}=25$ is shown in \cref{fig:cp_median_prediction}.

We additionally report a sample-complexity study at $n_{\xi}=25$, varying the training set size from $512$ to $8192$ while keeping the architecture and training procedure unchanged. As shown in the right panel of \cref{fig:cp_complexity}, the mean relative stress error decreases from $30.9\%$ at $512$ samples to $8.9\%$ at $8192$ samples. While finite-sample effects remain dominant relative to the internal variable dimension, the error decays more slowly with training samples than in viscoelasticity, indicating that homogenized crystal plasticity requires more data to learn well. This may partly reflect limitations of the data-generation strategy for learning shear-dominated modes. The normal stress components have large magnitudes, and the mean relative test errors are about $6\%$ for normal components and $30\%$ for shear components, compared with an aggregate error of about $9\%$. The larger shear errors are consistent with the visualization of the median-error test sample shown in \cref{fig:cp_median_prediction}. Further investigation of how these modes are learned could guide the design of more informative loading paths and active learning strategies to reduce data requirements; see, e.g., \cite{bhattacharya2026optimal,zhang2024iterated}.

\section{Conclusion}\label{sec:conclusion}

This paper studied microstructure-dependent recurrent neural operator models for learning inelastic homogenization. In these models, neural operators represent microstructure-dependent stress and kinetic responses, with a finite-dimensional internal variable that encodes the loading history, giving a Markovian representation of a history-dependent constitutive law. We proved that this model can approximate the homogenized constitutive law of Kelvin--Voigt viscoelasticity and quantified the internal-variable dimension sufficient to achieve a prescribed accuracy.

The analysis relied on an operator formulation of the homogenized Kelvin--Voigt law and its spectral representation of the memory kernel. { Based on this, we construct a Prony-series internal-variable model with relaxation rates shared across microstructures and continuous microstructure-dependent coefficients. Combining this model with neural-operator approximation of the stress and kinetic responses yielded uniform approximation over time, admissible microstructures, and strain trajectories. The resulting sufficient internal-variable dimension grows logarithmically with the inverse relative stress error and with each material contrast, holding the remaining parameters constant.}

The numerical experiments complement these conclusions. For Kelvin--Voigt viscoelasticity, the learned models substantially improve on the naive Markovian model without memory. The test error decreases as the internal-variable dimension or training dataset size increases, and increases as the material contrast grows. Our results indicate that the learned internal variable can be transformed into tensorial, stress-like variables, although their dynamics are generally coupled rather than exhibiting the decoupled structure arising in the theory. The models trained on Voronoi microstructures also retain predictive accuracy for particulate microstructures. For crystal plasticity, we show that the same data-driven modeling framework can learn nonlinear inelastic homogenization, demonstrating the framework's practicality beyond linear viscoelasticity theory.

There are numerous directions for future study: theoretical, methodological, and applied. A potential improvement of the complexity result involves sharpening the internal-variable complexity estimates using exponential-sum approximation of completely monotone functions \cite{kammler1976chebyshev,braess2005approximation,beylkin2010revisited,koyama2023exponential,darrow2026spectral}. Related history-variable constructions \cite{romero2026identification} motivate alternative memory representations, for which continuous dependence on microstructure remains a requirement. Other directions include parameter and sample complexity analyses \cite{lanthaler2026parametric,kovachki2024data} and extensions to nonlinear inelasticity \cite{masi2023evolution,ortiz2026encoding}. Finally, an important future direction is to develop methods that enforce thermodynamical consistency and well-posedness of the macroscopic balance equation. A good starting point is learning energetic formulations \cite{huang2022variational,flaschel2025convex, raj2026neural} with microstructure dependence.

\section*{Acknowledgment}

This work was supported by the U.S. Office of Naval Research through the Multidisciplinary University Research Initiative (MURI) program on Data-Driven Closure Relations under award number N00014-23-1-2654. George Stepaniants was supported by an NSF Mathematical Sciences Postdoctoral Research Fellowship (MSPRF) under award number 2402074. Andrew Stuart was also supported by the U.S. Department of Defense through a Vannevar Bush Faculty Fellowship under award number N00014-22-1-2790; this award also partially supported Lianghao Cao. The authors thank Kaushik Bhattacharya for suggestions that helped improve this paper.

OpenAI ChatGPT assisted in formulating proof strategies. OpenAI GPT-6.1-Sol suggested the bound improvement discussed in \cref{rmk:prony_construction_choices}. OpenAI GPT-5.3-codex in GitHub Copilot helped improve an earlier version of the viscoelasticity data-generation and training code published in \cite{bhattacharya2025learning} and implement the crystal plasticity data-generation pipeline. We used OpenAI Codex GPT-6 Astra to formalize our theoretical results in Lean 4.35.0, and we used Lean Comparator to verify the proofs against separately specified theorem statements. OpenAI GPT-6.1 Sol assisted in producing \texttt{TikZ} code for \cref{fig:graphical_abstract,fig:fnm_schematic,fig:local_chebyshev}. The authors take full responsibility for the manuscript's content.

\appendix

\makeatletter
\crefalias{section}{appendix}
\crefalias{subsection}{appendix} %
\@addtoreset{figure}{section}
\renewcommand{\thefigure}{\thesection.\arabic{figure}}
\makeatother

\section{Homogenized Tensors and Memory Kernel}

This appendix proves symmetry and definiteness of the homogenized tensors and memory kernel in \cref{app:homogenized_tensors_kernel_properties}, establishes the spectral representation of the kernel in \cref{app:spectral_rep}, and derives the homogenized model directly from the unit-cell problem in \cref{app:unit_cell_homogenization}.

\subsection{Symmetry and Definiteness of the Homogenized Tensors and Memory Kernel}\label{app:homogenized_tensors_kernel_properties}
We verify the symmetry and definiteness asserted in \cref{thm:viscoelastic_homogenization}. Let $\bV,\bE\in L^{\infty}_{\mathrm{per}}(\bT^d;\cV^4_{+})$, and let $\bV^{\mathrm{hom}}(\bV)$, $\bE^{\mathrm{hom}}(\bE)$, and $\bK(t;\bV,\bE)$ be defined by \cref{eq:homogenized_tensors,eq:memory_kernel}.

Since $\cV^{2}$ is finite-dimensional, it suffices to verify the symmetry and sign properties of the bilinear forms in \cref{eq:homogenized_tensors,eq:memory_kernel}. For $\bA\in\{\bV,\bE\}$, let $\cB_{\bA}(\overline e_1,\overline e_2)\coloneqq(\overline e_1,\cP_{\bA}^{\perp}\overline e_2)_{\bA}$. Because $\cP_{\bA}^{\perp}$ is the $\bA$-orthogonal projector onto the orthogonal complement of $\cX$, it is self-adjoint and idempotent with respect to $(\dummy,\dummy)_{\bA}$. Therefore $\cB_{\bA}(\overline e_1,\overline e_2)=(\cP_{\bA}^{\perp}\overline e_1,\cP_{\bA}^{\perp}\overline e_2)_{\bA}=\cB_{\bA}(\overline e_2,\overline e_1)$ and $\cB_{\bA}(\overline e,\overline e)=\|\cP_{\bA}^{\perp}\overline e\|_{\bA}^{2}\geq 0$. If $\cB_{\bA}(\overline e,\overline e)=0$, then $\cP_{\bA}^{\perp}\overline e=0$, so $\overline e=\cP_{\bA}\overline e\in \cX$. Since $\overline e$ is constant and $\cX$ contains no nonzero constant tensors, it follows that $\overline e=0$. Thus $\cB_{\bA}$ is positive definite, so $\bV^{\mathrm{hom}}(\bV),\bE^{\mathrm{hom}}(\bE)\in \cV^4_{+}$. For the kernel, set $\cA_t\coloneqq\cR_{\bV,\bE}\exp(-t\cR_{\bV,\bE})$ for $t\ge 0$. Since $\cR_{\bV,\bE}$ is positive and self-adjoint on $(\cX,(\dummy,\dummy)_{\bV})$, the spectral theorem implies that $\cA_t=g_t(\cR_{\bV,\bE})$, with $g_t(\lambda)\coloneqq\lambda \exp(-t\lambda)$, is also positive and self-adjoint. Hence $(\cA_t\cD_{\bV,\bE}\overline e_1,\cD_{\bV,\bE}\overline e_2)_{\bV}=(\cD_{\bV,\bE}\overline e_1,\cA_t\cD_{\bV,\bE}\overline e_2)_{\bV}$ and $(\cA_t\cD_{\bV,\bE}\overline e,\cD_{\bV,\bE}\overline e)_{\bV}\ge 0$, so $\bK(t;\bV,\bE)\in \cV^{4}$ is positive semidefinite for all $t\ge 0$.
\subsection{Proof of \cref{prop:spectral_rep}}\label{app:spectral_rep}

Write $\cR\coloneqq\cR_{\bV,\bE}$, $\cD\coloneqq\cD_{\bV,\bE}$, $\nu\coloneqq\nu_{\bV,\bE}$, and $\cI\coloneqq \cI_{\bV,\bE}$, and set $\lambda_{\min}\coloneqq \alpha_{\bE}/\beta_{\bV}$ and $\lambda_{\max}\coloneqq \beta_{\bE}/\alpha_{\bV}$, so $\cI=[\lambda_{\min},\lambda_{\max}]$. By the Riesz representation theorem, $\cR\in\cL(\cX)$ is the self-adjoint positive operator on $(\cX,(\dummy,\dummy)_{\bV})$. The coercivity and boundedness assumptions yield, for all $\eta\in\cX$, $\lambda_{\min}\|\eta\|_{\bV}^2
\leq (\cR\eta,\eta)_{\bV}
\leq \lambda_{\max}\|\eta\|_{\bV}^2$.
Thus, the spectrum of $\cR$ is contained in $\cI$. By the spectral theorem, there exists a projection-valued measure $\Pi\coloneqq\Pi_{\bV,\bE}$ supported on $\cI$ such that $\cR\exp(-t\cR)=\int_{\cI} \lambda \exp(-\lambda t)\,\mathrm{d}\Pi(\lambda)$ for all $t\ge 0$. For each Borel set $B\subset\cI$, define a bounded bilinear form on $\cV^{2}$ by
\begin{equation}\label{eq:nu_def}
(\overline{e}_1, \nu(B)\overline{e}_2)_F\coloneqq\bigl(\Pi(B)\cD \overline{e}_1,\ \cD \overline{e}_2\bigr)_{\bV}.
\end{equation}
Since $\Pi(B)$ is self-adjoint, $\nu(B)$ is symmetric. Since $\cV^{2}$ is finite-dimensional, \cref{eq:nu_def} defines a unique tensor in $\cV^4$, and the strong countable additivity of $\Pi$ implies that $\nu$ is countably additive. Moreover, for every Borel set $B\subset\cI$ and every $\overline e_1,\overline e_2\in\cV^{2}$, we have $|(\overline e_1, \nu(B)\overline e_2)_F|\le \|\Pi(B)\cD\overline e_1\|_{\bV}\,\|\cD\overline e_2\|_{\bV}\le \|\cD\overline e_1\|_{\bV}\,\|\cD\overline e_2\|_{\bV}$, so $\nu$ is finite. Combining \cref{eq:memory_kernel}, the spectral representation of $\cR\exp(-t\cR)$, and \cref{eq:nu_def} gives \cref{eq:K_laplace}.

It remains to prove \cref{eq:uniform_TV_bound_nu_explicit}. First, 
for any Borel set $B\subset\cI$, the operator $\Pi(B)$ is positive on $\cX$. Consequently, \cref{eq:nu_def} gives, for any $\overline{e} \in \cV^{2}$, $(\overline{e}, \nu(B)\overline{e})_F\geq 0$. Let $\{\overline e_\alpha\}_{\alpha=1}^{d_{\cV^2}}$ be an orthonormal basis of $\cV^2$ with respect to the Frobenius inner product. For any finite Borel partition $B_1,\dots,B_m$ of $\cI$, the inequality $\|\bA\|_F\le \operatorname{tr}(\bA)$ for positive semidefinite $\bA\in\cV^4$ gives $\sum_{i=1}^m \|\nu(B_i)\|_F
\le \sum_{i=1}^m \operatorname{tr}(\nu(B_i))
=\operatorname{tr}(\nu(\cI))
=\sum_{\alpha=1}^{d_{\cV^2}}\|\cD\overline e_\alpha\|_{\bV}^2$, where the last identity follows from \cref{eq:nu_def} and {$\Pi(\cI)=\mathrm{Id}_{\cX}$}. Also, for any $\overline e\in\cV^2$, since $\cP_{\bE}\overline e\in\cX$ and $\cP_{\bV}$ is the $\bV$-orthogonal projector onto $\cX$,
\begin{equation}\label{eq:mismatch_bound}
\|\cD\overline e\|_{\bV}^2
=\|\cP_{\bV}(\overline e-\cP_{\bE}\overline e)\|_{\bV}^2
\le \|\overline e-\cP_{\bE}\overline e\|_{\bV}^2
\le \frac{\beta_{\bV}}{\alpha_{\bE}}\|\overline e-\cP_{\bE}\overline e\|_{\bE}^2
\le \frac{\beta_{\bV}\beta_{\bE}}{\alpha_{\bE}}\|\overline e\|_F^2.
\end{equation}
Applying this bound to each basis vector $\overline e_\alpha$ and taking the supremum over all finite Borel partitions gives \cref{eq:uniform_TV_bound_nu_explicit}.

\subsection{Deriving the Homogenized Model from the Unit Cell Problem} \label{app:unit_cell_homogenization}

We verify the unit-cell representation \cref{eq:cell_constitutive_map} by deriving the homogenized model in \cref{eq:homogenized_model_identified} directly from the unit cell problem \cref{eqns:cell_problem}, complementing the Laplace-transform proof of \cref{thm:viscoelastic_homogenization}. {Let $T>0$ be arbitrary and finite, and write} $\cR\coloneqq\cR_{\bV,\bE}$ and $\cD\coloneqq\cD_{\bV,\bE}$.

Testing the balance law in \cref{eqns:cell_problem} against $\phi\in\mathring{H}^1_{\mathrm{per}}(\bT^d;\bR^d)$ and integrating by parts gives
\begin{equation}\label{eq:cell_weak}
(\dot{\overline{e}}(t)+\dot\eta(t),\psi)_{\bV}+(\overline{e}(t)+\eta(t),\psi)_{\bE}=0\quad \text{a.e. }t\in{[0,T]}
\end{equation}
where $\psi=\mathrm{sym}(\nabla_y\phi)\in\cX$ and $\eta(0)=-\cP_{\bE}\overline{e}(0)$. Since $\eta(t),\dot\eta(t)\in\cX$ for a.e.\ $t\in{[0,T]}$, the definitions of $\cP_{\bV}$, $\cP_{\bE}$, and $\cR$ give $(\dot{\overline{e}}(t),\psi)_{\bV}=(\cP_{\bV}\dot{\overline{e}}(t),\psi)_{\bV}$, $(\overline{e}(t),\psi)_{\bE}=(\cR\cP_{\bE}\overline{e}(t),\psi)_{\bV}$, and $(\eta(t),\psi)_{\bE}=(\cR\eta(t),\psi)_{\bV}$ for all $\psi\in\cX$ and a.e.\ $t\in{[0,T]}$. Hence \cref{eq:cell_weak} is equivalent to
\begin{equation}\label{eq:cell_ode}
\dot\eta(t)+\cR\eta(t)=-\cP_{\bV}\dot{\overline{e}}(t)-\cR\cP_{\bE}\overline{e}(t),\quad\eta(0)=-\cP_{\bE}\overline{e}(0).
\end{equation}
Set $\omega(t)\coloneqq\eta(t)+\cP_{\bE}\overline{e}(t)$ for a.e.\ $t\in{[0,T]}$. Then \cref{eq:cell_ode} and $\cD=\cP_{\bV}-\cP_{\bE}$ give $\dot\omega(t)+\cR\omega(t)=-\cD\dot{\overline{e}}(t)$ with $\omega(0)=0$. Integrating in time gives
\begin{equation}\label{eq:eta_sol}
\omega(t)=-\int_0^t\exp(-(t-s)\cR)\cD\dot{\overline{e}}(s)\,\mathrm{d}s.
\end{equation}
For any $\overline{e}_1\in\cV^2$, the macroscopic stress satisfies $(\overline{e}_1, \overline{\sigma}(t))_F=(\overline{e}_1,\dot{\overline{e}}(t)+\dot\eta(t))_{\bV}+(\overline{e}_1,\overline{e}(t)+\eta(t))_{\bE}$ for a.e.\ $t\in{[0,T]}$. Using $\eta(t)=\omega(t)-\cP_{\bE}\overline{e}(t)$ gives $\overline{e}(t)+\eta(t)=\cP_{\bE}^{\perp}\overline{e}(t)+\omega(t)$ and $\dot{\overline{e}}(t)+\dot\eta(t)=\cP_{\bE}^{\perp}\dot{\overline{e}}(t)+\dot\omega(t)$. Since $\cP_{\bE}^{\perp}=\cP_{\bV}^{\perp}+\cD$, this becomes
\begin{equation}\label{eq:stress_decomp}
(\overline{e}_1, \overline{\sigma}(t))_{F}
=\underbrace{(\overline{e}_1,\cP_{\bV}^{\perp}\dot{\overline{e}}(t))_{\bV}}_{\displaystyle (\overline{e}_1, \bV^{\mathrm{hom}}\dot{\overline{e}}(t))_F}
+\underbrace{(\overline{e}_1,\cP_{\bE}^{\perp}\overline{e}(t))_{\bE}}_{\displaystyle(\overline{e}_1,\bE^{\mathrm{hom}}\overline{e}(t))_F}
+\underbrace{(\overline{e}_1,\cD\dot{\overline{e}}(t)+\dot\omega(t))_{\bV}+(\overline{e}_1,\omega(t))_{\bE}}_{\displaystyle\eqqcolon R},
\end{equation}
where the first two terms are identified by \cref{eq:homogenized_tensors}.

Since $\dot\omega(t)+\cR\omega(t)=-\cD\dot{\overline{e}}(t)$, we have $\cD\dot{\overline{e}}(t)+\dot\omega(t)=-\cR\omega(t)$. The definitions of $\cP_{\bV}$, $\cP_{\bE}$, and $\cR$ imply $(\overline{e}_1,\cR\omega(t))_{\bV}=(\cR\cP_{\bV}\overline{e}_1,\omega(t))_{\bV}$ and $(\overline{e}_1,\omega(t))_{\bE}=(\cR\cP_{\bE}\overline{e}_1,\omega(t))_{\bV}$, so $R=-(\cR\cD\overline{e}_1,\omega(t))_{\bV}$. Substituting \cref{eq:eta_sol} and using the self-adjointness and commutativity of $\cR$ and $\exp(-t\cR)$ on $(\cX,(\dummy,\dummy)_{\bV})$ gives
\begin{equation*}
    R=\int_0^t\bigl(\cR\exp(-(t-s)\cR)\cD\overline{e}_1,\cD\dot{\overline{e}}(s)\bigr)_{\bV}\,\mathrm{d}s=\int_0^t\left(\overline{e}_1, \bK(t-s;\bV,\bE)\dot{\overline{e}}(s)\right)_F\,\mathrm{d}s
\end{equation*}
where the second equality is \cref{eq:memory_kernel}. Substituting into \cref{eq:stress_decomp} yields \cref{eq:homogenized_model_identified}.

\section{\texorpdfstring{{Supporting Results for \cref{thm:ua}}}{Supporting Results for the Universal Approximation Theorem}}\label{app:ua}

{Throughout this appendix, $\cM$ and $\cE$ are the admissible sets in \cref{thm:ua}, and {$\cT=[0,\infty)$}. We first establish the compactness and continuity properties used in \ref{item:ua_proof_s1}.}

\begin{lemma}[{Strong Convergence of Multiplication Operators}]\label{lem:strong_conv_tensor}
Let $\{\bA_n\}_{n\in\bN}\subset L^\infty_{\rm per}(\bT^d;\cV^4)$ satisfy
$\sup_{n\in\bN}\ \|\bA_n\|_{L^\infty_{\rm per}}\le \kappa_{\bA}<\infty$ and $\bA_n\to\bA$ in $L^1_{\rm per}(\bT^d;\cV^4)$. Then $\bA\in L^\infty_{\rm per}$ with $\|\bA\|_{L^\infty_{\rm per}}\le \kappa_{\bA}$ and $\| (\bA_n-\bA)e\|_{L^2_{\rm per}}\to 0$ for every $e\in L^2_{\rm per}(\bT^d;\cV^{2})$.
\end{lemma}

\begin{proof}[\bf Proof]
Since $\bA_n\to\bA$ in $L^1_{\rm per}(\bT^d;\cV^4)$ and $\bT^d$ has finite measure, a subsequence converges pointwise a.e.\ to $\bA$. Passing to the limit along this subsequence and using $\|\bA_n(y)\|_F\le \kappa_{\bA}$ a.e.\ gives $\|\bA(y)\|_F\le \kappa_{\bA}$ a.e. Hence $\bA\in L^\infty_{\rm per}$ and $\|\bA\|_{L^\infty_{\rm per}}\le \kappa_{\bA}$. Now take $e\in L^2_{\rm per}(\bT^d;\cV^2)$. By the subsequence criterion in $L^2_{\rm per}$, it suffices to show that every subsequence of $\{(\bA_n-\bA)e\}_{n\in\bN}$ has a further subsequence converging to $0$. Let $\bA_{n_j}$ be any subsequence. Since $\bA_{n_j}\to\bA$ in $L^1_{\rm per}$, after passing to a further subsequence we have $\bA_{n_j}(y)\to\bA(y)$ for a.e.\ $y\in\bT^d$. Therefore, $(\bA_{n_j}(y)-\bA(y))e(y)\to0$ for a.e.\ $y\in\bT^d$, and
$\|((\bA_{n_j}(y)-\bA(y))e(y))\|_F^2
\le \|\bA_{n_j}(y)-\bA(y)\|_F^2\|e(y)\|_F^2
\le (2\kappa_{\bA})^2\|e(y)\|_F^2$ for a.e.\ $y\in\bT^d$. The right-hand side belongs to $L^1_{\rm per}(\bT^d)$, so the dominated convergence theorem gives $\|(\bA_{n_j}-\bA)e\|_{L^2_{\rm per}}\to0$ along the further subsequence. The subsequence criterion therefore implies $\|(\bA_n-\bA)e\|_{L^2_{\rm per}}\to0$.
\end{proof}

\begin{lemma}[Compactness of the Admissible Set $\cM$ in $L^2_{\mathrm{per}}$] \label{lem:bv_compact}
The admissible set $\cM$ defined in \cref{eq:admissible_microstructure} is compact in $L^2_{\mathrm{per}}(\bT^d;\cV^{4})\times
L^2_{\mathrm{per}}(\bT^d;\cV^{4})$.
\end{lemma}

\begin{proof}[\bf Proof]
Let $\{(\bV_n,\bE_n)\}_{n\in\bN}\subset \cM$. Since $L^2_{\mathrm{per}}\times L^2_{\mathrm{per}}$ is a metric space, it suffices to prove sequential compactness. By \cref{eq:bv_bound_microstructure}, $\{\bV_n\}_{n\in\bN}$ and $\{\bE_n\}_{n\in\bN}$ are bounded in $\mathrm{BV}$. The compact embedding $\mathrm{BV}(\bT^d)\hookrightarrow L^1(\bT^d)$ gives a subsequence (not relabeled) and tensors $\bV,\bE\in L^1_{\mathrm{per}}(\bT^d;\cV^{4})$ such that
$ \bV_n \to \bV \quad\text{in }L^1_{\mathrm{per}}(\bT^d;\cV^{4})$ and $\bE_n \to \bE \quad\text{in }L^1_{\mathrm{per}}(\bT^d;\cV^{4})$. Applying \cref{lem:strong_conv_tensor} gives $\bV,\bE\in L^\infty_{\mathrm{per}}(\bT^d;\cV^4)$, $\|\bV\|_{L^\infty_{\mathrm{per}}}\le \sqrt{d_{\cV^2}}\,\beta_{\rm visc}$, and $\|\bE\|_{L^\infty_{\mathrm{per}}}\le \sqrt{d_{\cV^2}}\,\beta_{\rm elas}$. Since $\|\bV_n-\bV\|_{L^\infty_{\mathrm{per}}}\le 2\sqrt{d_{\cV^2}}\,\beta_{\rm visc}$ and $\|\bE_n-\bE\|_{L^\infty_{\mathrm{per}}}\le 2\sqrt{d_{\cV^2}}\,\beta_{\rm elas}$, H\"older's inequality gives
\begin{equation*}
    \|\bV_n-\bV\|_{L^2_{\mathrm{per}}}^2
\leq \|\bV_n-\bV\|_{L^\infty_{\mathrm{per}}}\|\bV_n-\bV\|_{L^1_{\mathrm{per}}}
    \le 2\sqrt{d_{\cV^2}}\,\beta_{\rm visc}\|\bV_n-\bV\|_{L^1_{\mathrm{per}}}\to 0.
\end{equation*}
Therefore, $\bV_n\to \bV$ strongly in $L^2_{\mathrm{per}}(\bT^d;\cV^{4})$. The same argument gives $\bE_n\to \bE$ strongly in $L^2_{\mathrm{per}}(\bT^d;\cV^{4})$. Lower semicontinuity of total variation under $L^1$ convergence implies $\|\bV\|_{\mathrm{BV}}\le \liminf_{n\to\infty}\|\bV_n\|_{\mathrm{BV}}\leq \gamma_{\rm visc}$ and $\|\bE\|_{\mathrm{BV}}\le \liminf_{n\to\infty}\|\bE_n\|_{\mathrm{BV}}\leq \gamma_{\rm elas}$, so the BV bounds are preserved. For the uniform ellipticity bounds, \cref{lem:strong_conv_tensor} and Cauchy--Schwarz give, for every $e\in L^2_{\mathrm{per}}(\bT^d;\cV^2)$,
\begin{equation*}
    \left|\int_{\bT^d} \left(e(y),(\bV_n(y)-\bV(y))e(y)\right)_F\,\mathrm{d}y\right| \le \|e\|_{L^2_{\mathrm{per}}}\|(\bV_n-\bV)e\|_{L^2_{\mathrm{per}}} \to 0.
\end{equation*}
Therefore, $\|e\|_{\bV_n}\to \|e\|_{\bV}$, and passing to the limit in the uniform ellipticity bounds for $\bV_n$ leads to the same bounds for $\bV$. The same argument applies to $\bE$. Therefore $(\bV,\bE)\in \cM$.
\end{proof}

\begin{proposition}[{Strong Convergence of the Unit-Cell Operators}]\label{prop:strong_conv}
Let $\{(\bV_n,\bE_n)\}_{n\in\bN}\subset\cM$ and assume that for some pair $(\bV,\bE)$, $\bV_n\to\bV$ and $\bE_n\to\bE$ in $L^1_{\rm per}(\bT^d;\cV^4)$. We have
\begin{gather}
    \|(\cP_{\bV_n}-\cP_{\bV})e\|_{L^2_{\rm per}}\to0,\quad \|(\cP_{\bE_n}-\cP_{\bE})e\|_{L^2_{\rm per}}\to0\quad\text{for every }e\in L^2_{\rm per}, \label{eq:P_conv}\\
    \|(\cD_{\bV_n,\bE_n}-\cD_{\bV,\bE})e\|_{L^2_{\rm per}}\to0\quad\text{for every }e\in L^2_{\rm per}, \label{eq:D_conv}\\
    \|(\cR_{\bV_n,\bE_n}-\cR_{\bV,\bE})\eta\|_{L^2_{\rm per}}\to0\quad\text{for every }\eta\in\cX, \label{eq:R_conv}
\end{gather}
Furthermore, for all $k\in\bN$, we have
\begin{gather}
\|(\cR^k_{\bV_n,\bE_n}-\cR^k_{\bV,\bE})\eta\|_{L^2_{\rm per}}\to0\quad\text{for every }\eta\in\cX. \label{eq:R_k_conv}
\end{gather}
\end{proposition}

\begin{proof}[\bf Proof]
{By \cref{lem:bv_compact}, the limit pair $(\bV,\bE)$ belongs to $\cM$. Let $e\in L^2_{\rm per}(\bT^d;\cV^2)$ and set $\eta_n\coloneqq\cP_{\bV_n}e$, $\eta\coloneqq\cP_{\bV}e$. By \cref{eqns:operator_P}, $(\eta_n-e,\psi)_{\bV_n}=0$ and $(\eta-e,\psi)_{\bV}=0$ for all $\psi\in\cX$. Choosing $\psi=\eta_n-\eta$, subtracting the two identities, and applying coercivity and Cauchy--Schwarz give
\begin{equation*}
    \alpha_{\rm visc}\|\eta_n-\eta\|_{L^2_{\rm per}}^2
\le \|(\bV_n-\bV)(e-\eta)\|_{L^2_{\rm per}}\|\eta_n-\eta\|_{L^2_{\rm per}}.
\end{equation*}
If $\|\eta_n-\eta\|_{L^2_{\rm per}}>0$, we divide by this norm, and otherwise there is nothing to prove. Thus $\alpha_{\rm visc}\|\eta_n-\eta\|_{L^2_{\rm per}}
\le \|(\bV_n-\bV)(e-\eta)\|_{L^2_{\rm per}}\to 0$ by \cref{lem:strong_conv_tensor}. Hence $\cP_{\bV_n}e\to\cP_{\bV}e$ in $L^2_{\rm per}$. The same argument with $\bE_n\to\bE$ gives $\cP_{\bE_n}e\to\cP_{\bE}e$, and therefore \cref{eq:P_conv,eq:D_conv} follow.}

Now set $\cR_n\coloneqq\cR_{\bV_n,\bE_n}$ and $\cR\coloneqq\cR_{\bV,\bE}$. For $\eta\in\cX$, let $\varrho_n\coloneqq\cR_n\eta$ and $\varrho\coloneqq\cR\eta$. {Using \cref{eq:operator_L} for $\cR_n$ and $\cR$, choosing $\psi=\varrho_n-\varrho$, and arguing as above, we obtain
\begin{equation}\label{eq:R_strong}
\alpha_{\rm visc}\|\varrho_n-\varrho\|_{L^2_{\rm per}}
\le \|(\bE_n-\bE)\eta\|_{L^2_{\rm per}}+\|(\bV_n-\bV)\varrho\|_{L^2_{\rm per}}\to 0,
\end{equation}
again by \cref{lem:strong_conv_tensor}. Thus \cref{eq:R_conv} holds.

Testing the identity $(\cR_n\eta,\psi)_{\bV_n}=(\eta,\psi)_{\bE_n}$ from \cref{eq:operator_L} with $\psi=\cR_n\eta$ and using the coercivity of $\bV_n$ together with the upper bound on $\bE_n$ gives
\begin{equation}\label{eq:R_uniform_bound}
\sup_{n\in\bN}\|\cR_n\|_{\cL(\cX)} \le \beta_{\rm elas}/\alpha_{\rm visc}.
\end{equation}
The convergence of $\cR_n^k$ follows by induction on $k$. The case $k=1$ is \cref{eq:R_strong}. If the claim holds for $k-1$, then $\cR_n^k\eta-\cR^k\eta=\cR_n\bigl(\cR_n^{k-1}\eta-\cR^{k-1}\eta\bigr)+(\cR_n-\cR)\cR^{k-1}\eta$. The first term tends to $0$ by the inductive hypothesis and \cref{eq:R_uniform_bound}, while the second tends to $0$ by \cref{eq:R_strong}. Therefore, \cref{eq:R_k_conv} holds.}
\end{proof}

\begin{corollary}[Continuity of the Homogenized Tensors]\label{cor:continuity_tensors}
Let $\{(\bV_n,\bE_n)\}_{n\in\bN}\subset\cM$ and assume that for some pair $(\bV,\bE)$, $\bV_n\to\bV$ and $\bE_n\to\bE$ in $L^1_{\rm per}(\bT^d;\cV^4)${.} Then $(\bV,\bE)\in\cM$ and we have $\bV^{\mathrm{hom}}(\bV_n)\to \bV^{\mathrm{hom}}(\bV)$ and $\bE^{\mathrm{hom}}(\bE_n)\to \bE^{\mathrm{hom}}(\bE)$.
\end{corollary}

\begin{proof}[\bf Proof]
By \cref{lem:bv_compact}, the limit pair $(\bV,\bE)$ belongs to $\cM$. Let $\overline e_1,\overline e_2\in\cV^2$. Using \cref{eq:homogenized_tensors} and adding and subtracting $(\overline{e}_1,\cP_{\bV}^{\perp}\overline e_2)_{\bV_n}$, we obtain
\begin{align*}
    \bigl| \left(\overline e_1, \bV^{\mathrm{hom}}(\bV_n)\overline{e}_2\right)_F-\left(\overline e_1,\bV^{\mathrm{hom}}(\bV)\overline e_2\right)_F\bigr|\le &\bigl|(\overline e_1,(\cP_{\bV_n}^{\perp}-\cP_{\bV}^{\perp})\overline e_2)_{\bV_n}\bigr|\\
    &+\left|\int_{\bT^d}\left(\overline e_1,(\bV_n-\bV)[\cP_{\bV}^{\perp}\overline e_2](y)\right)_F \,\mathrm{d}y\right|.
\end{align*}
The first term is bounded by $\|\bV_n\|_{L^\infty_{\rm per}}\|\overline e_1\|_F\|(\cP_{\bV_n}-\cP_{\bV})\overline e_2\|_{L^2_{\rm per}}$ and therefore converges to $0$ by the uniform $L^\infty_{\rm per}$ bound on $\bV_n$ and \cref{prop:strong_conv}. The second term is bounded by $\|(\bV_n-\bV)\overline e_1\|_{L^2_{\rm per}}\|\cP_{\bV}^{\perp}\overline e_2\|_{L^2_{\rm per}}$ and also converges to $0$ by \cref{lem:strong_conv_tensor}. Therefore, $\left(\overline e_1,\bV^{\mathrm{hom}}(\bV_n)\overline e_2\right)_F\to {\left(\overline e_1,\bV^{\mathrm{hom}}(\bV)\overline e_2\right)_F}$ for all $\overline e_1,\overline e_2\in\cV^2$. Since $\cV^2$ is finite-dimensional, this is equivalent to $\bV^{\mathrm{hom}}(\bV_n)\to \bV^{\mathrm{hom}}(\bV)$ in $\cV^4$. The argument for $\bE^{\mathrm{hom}}$ is identical.
\end{proof}

\begin{lemma}[Weak Convergence of the Spectral Measure]
\label{lem:weak_convergence}
Suppose $\{(\bV_n,\bE_n)\}_{n\in\bN}\subset\cM$ satisfies $\bV_n\to\bV$ and $\bE_n\to\bE$ in $L^1_{\mathrm{per}}(\bT^d;\cV^{4})$. Then the sequence of measures $\nu_{\bV_n,\bE_n}$, as defined by \cref{eq:nu_def}, converges weakly to $\nu_{\bV,\bE}$. That is, for any continuous function $f \in C(\cI)$ with $\cI\coloneqq[\alpha_{\rm elas}/\beta_{\rm visc}, \beta_{\rm elas}/\alpha_{\rm visc}]$ and any $\overline{e}_1, \overline{e}_2 \in \cV^{2}$,
\begin{equation*}
    \lim_{n \to \infty} \int_{\cI} f(\lambda) \, \mathrm{d}(\overline{e}_1,\nu_{\bV_n,\bE_n}(\lambda)\overline{e}_2)_F = \int_{\cI} f(\lambda) \, \mathrm{d}(\overline{e}_1,\nu_{\bV,\bE}(\lambda)\overline{e}_2)_F.
\end{equation*}
\end{lemma}

\begin{proof}[\bf Proof]
By \cref{lem:bv_compact}, the limit pair $(\bV,\bE)$ belongs to $\cM$. Write $\nu_n\coloneqq \nu_{\bV_n,\bE_n}$, $\nu\coloneqq \nu_{\bV,\bE}$, $\cR_n\coloneqq \cR_{\bV_n,\bE_n}$, $\cR\coloneqq \cR_{\bV,\bE}$, $\cD_n\coloneqq \cD_{\bV_n,\bE_n}$, and $\cD\coloneqq \cD_{\bV,\bE}$. For $g\in C(\cI)$, let $g(\cR_n)\coloneqq \int_{\cI} g(\lambda)\,\mathrm{d}\Pi_{\bV_n,\bE_n}(\lambda)$, and define $g(\cR)$ analogously. Then \cref{eq:nu_def} gives $\int_{\cI} f(\lambda)\,\mathrm{d}(\overline e_1,\nu_n\overline e_2)_F
=\Bigl(f(\cR_n)\cD_n\overline e_1,\cD_n\overline e_2\Bigr)_{\bV_n}$,
and the same identity holds with $n$ replaced by the limit pair $(\bV,\bE)$. Therefore, it suffices to prove
\begin{equation}\label{eq:weak_converge_g}
    \Bigl(f(\cR_n)\cD_n\overline e_1,\cD_n\overline e_2\Bigr)_{\bV_n}\to \Bigl(f(\cR)\cD\overline e_1,\cD\overline e_2\Bigr)_{\bV}.
\end{equation}

{For any $\varepsilon_p>0$, the Weierstrass approximation theorem gives a polynomial $p$ such that
\begin{equation*}
\sup_{\lambda\in\cI}|f(\lambda)-p(\lambda)|\le\varepsilon_p.
\end{equation*}
}
Adding and subtracting the polynomial terms gives
\begin{align}
&\left|\Bigl(f(\cR_n)\cD_n\overline e_1,\cD_n\overline e_2\Bigr)_{\bV_n}-\Bigl(f(\cR)\cD\overline e_1,\cD\overline e_2\Bigr)_{\bV}\right|\nonumber\\
&\le \left|\Bigl((f-p)(\cR_n)\cD_n\overline e_1,\cD_n\overline e_2\Bigr)_{\bV_n}\right|
+\left|\Bigl(p(\cR_n)\cD_n\overline e_1,\cD_n\overline e_2\Bigr)_{\bV_n}-\Bigl(p(\cR)\cD\overline e_1,\cD\overline e_2\Bigr)_{\bV}\right|\nonumber\\
&\qquad +\left|\Bigl((p-f)(\cR)\cD\overline e_1,\cD\overline e_2\Bigr)_{\bV}\right|.\label{eq:weak_converge_decomposed}
\end{align}
For any $\eta\in\cX$,

\begin{equation*}
    \|(f-p)(\cR_n)\eta\|_{\bV_n}^2
=\int_{\cI}|f(\lambda)-p(\lambda)|^2\,\mathrm{d}(\Pi_{\bV_n,\bE_n}(\lambda)\eta,\eta)_{\bV_n}
\le \varepsilon_p^2\|\eta\|_{\bV_n}^2,
\end{equation*}
and the same estimate holds with $(\cR_n,\bV_n,\Pi_{\bV_n,\bE_n})$ replaced by $(\cR,\bV,\Pi_{\bV,\bE})$. {Cauchy--Schwarz and \cref{eq:mismatch_bound} therefore show} that the first and third terms in \cref{eq:weak_converge_decomposed} are each bounded by $\varepsilon_p\beta_{\rm visc}\beta_{\rm elas}\|\overline e_1\|_F\|\overline e_2\|_F/\alpha_{\rm elas}$.

For the polynomial term, set $\phi_n\coloneqq p(\cR_n)\cD_n\overline e_1$, $\phi\coloneqq p(\cR)\cD\overline e_1$, $\psi_n\coloneqq \cD_n\overline e_2$, and $\psi\coloneqq \cD\overline e_2$. Then the middle term in \cref{eq:weak_converge_decomposed} is $| (\phi_n,\psi_n)_{\bV_n}-(\phi,\psi)_{\bV}|$. By \cref{eq:D_conv}, $\psi_n\to\psi$ strongly in $L^2_{\rm per}$. If $p(\lambda)=\sum_{k=0}^K c_k\lambda^k$, then $\phi_n-\phi=\sum_{k=0}^K c_k\bigl(\cR_n^k\cD_n\overline e_1-\cR^k\cD\overline e_1\bigr)$. For $k=0,\dots,K$,
\begin{equation*}
\cR_n^k\cD_n\overline e_1-\cR^k\cD\overline e_1
=\cR_n^k(\cD_n-\cD)\overline e_1+(\cR_n^k-\cR^k)\cD\overline e_1.
\end{equation*}
The first term tends to $0$ because \cref{eq:D_conv} gives $\cD_n\overline e_1\to \cD\overline e_1$ in $L^2_{\rm per}$ and \cref{eq:R_uniform_bound} implies that $\{\cR_n^k\}_n$ is uniformly bounded on $\cX$ for fixed $k$. The second term tends to $0$ by \cref{eq:R_k_conv} with $\eta=\cD\overline e_1\in\cX$. Hence $\phi_n\to\phi$ strongly in $L^2_{\rm per}$. {We now write}
\begin{align*}
    (\phi_n,\psi_n)_{\bV_n}-(\phi,\psi)_{\bV}
    &=((\phi_n-\phi),\psi_n)_{\bV_n}+\int_{\bT^d}\phi(y):(\bV_n(y)-\bV(y)):\psi_n(y)\,\mathrm{d}y\\
    &\qquad +(\phi,(\psi_n-\psi))_{\bV}.
\end{align*}
The first and third terms tend to $0$ by Cauchy--Schwarz, the uniform $L^\infty_{\rm per}$ bounds on $\bV_n$ and $\bV$, and the strong $L^2_{\rm per}$ convergence of $\phi_n$ and $\psi_n$. The integral term tends to $0$ by \cref{lem:strong_conv_tensor}, since $\bV_n\to\bV$ in $L^1_{\rm per}$ and $(\psi_n)$ is bounded in $L^2_{\rm per}$. Consequently, the polynomial term in \cref{eq:weak_converge_decomposed} converges to $0$. Taking $\limsup_{n\to\infty}$ in \cref{eq:weak_converge_decomposed} and then letting $\varepsilon_p\to0$ gives \cref{eq:weak_converge_g}.
\end{proof}

{
We next construct the Prony-series approximation used in \ref{item:ua_proof_s2}.

\begin{lemma}[Weighted Spectral Measure]\label{lem:weighted_spectral_measure}
Set $\lambda_{\min}\coloneqq\alpha_{\rm elas}/\beta_{\rm visc}$, $\lambda_{\max}\coloneqq\beta_{\rm elas}/\alpha_{\rm visc}$, and $\cI\coloneqq[\lambda_{\min},\lambda_{\max}]$. For $(\bV,\bE)\in\cM$, define the positive semidefinite tensor-valued measure
\begin{equation}\label{eq:weighted_measure_def}
\mathrm d\mu_{\bV,\bE}(\lambda)\coloneqq\frac{\lambda}{\beta_{\rm elas}/\beta_{\rm visc}+\lambda}\,\mathrm d\nu_{\bV,\bE}(\lambda).
\end{equation}
Then
\begin{equation}\label{eq:weighted_measure_mass}
\sup_{(\bV,\bE)\in\cM}\|\mu_{\bV,\bE}\|_{\mathrm{TV}(\cI)}<2d_{\cV^2}\beta_{\rm visc}.
\end{equation}
Moreover, $(\bV,\bE)\mapsto\mu_{\bV,\bE}$ is weakly continuous under $L^1_{\rm per}$ convergence in $\cM$.
\end{lemma}

\begin{proof}[\bf Proof]
By \cref{eq:nu_def}, $\nu_{\bV,\bE}(B)$ is positive semidefinite for every Borel set $B\subset\cI$, so the same holds for $\mu_{\bV,\bE}(B)$. For $s>0$, let $\bA_s=s^{-1}\bE+\bV$, as in the proof of \cref{thm:viscoelastic_homogenization}. The Laplace-transform identity there and \cref{eq:K_laplace} give
\begin{equation*}
\int_{\cI}\frac{\lambda}{s+\lambda}\,\mathrm d\nu_{\bV,\bE}(\lambda)
=\hat{\bK}(s)=\bA_s^{\rm hom}-s^{-1}\bE^{\rm hom}-\bV^{\rm hom}.
\end{equation*}
The variational characterization of static homogenization, with a constant trial strain, gives $(\overline e,\bA_s^{\rm hom}\overline e)_F\le(s^{-1}\beta_{\rm elas}+\beta_{\rm visc})\|\overline e\|_F^2$. Since periodic strain correctors have zero mean, the same characterization gives $(\overline e,\bV^{\rm hom}\overline e)_F\ge\alpha_{\rm visc}\|\overline e\|_F^2$ and $(\overline e,\bE^{\rm hom}\overline e)_F\ge\alpha_{\rm elas}\|\overline e\|_F^2$. At $s=\beta_{\rm elas}/\beta_{\rm visc}$, these bounds imply
\begin{equation*}
(\overline e,\mu_{\bV,\bE}(\cI)\overline e)_F\le\left(2\beta_{\rm visc}-\alpha_{\rm visc}-\frac{\beta_{\rm visc}\alpha_{\rm elas}}{\beta_{\rm elas}}\right)\|\overline e\|_F^2.
\end{equation*}
The inequality $\|\bA\|_F\le\operatorname{tr}\bA$ for positive semidefinite $\bA\in\cV^4$ gives $\|\mu_{\bV,\bE}\|_{\mathrm{TV}(\cI)}\le\operatorname{tr}\mu_{\bV,\bE}(\cI)$, as in the proof of \cref{prop:spectral_rep}. Taking the trace in the preceding bound and using $\alpha_{\rm visc},\alpha_{\rm elas}>0$ proves \cref{eq:weighted_measure_mass}. Finally, for every $f\in C(\cI)$, the function $\lambda\mapsto f(\lambda)\lambda/(\beta_{\rm elas}/\beta_{\rm visc}+\lambda)$ is continuous. Weak continuity follows from \cref{lem:weak_convergence}.
\end{proof}

\begin{lemma}[Weighted Relaxation Response]\label{lem:weighted_response}
For $\overline e\in\cE$, $t\ge0$, and $\lambda\in\bC$, define
\begin{equation}\label{eq:weighted_response_def}
\cW_{\overline e,t}(\lambda)\coloneqq\left(\frac{\beta_{\rm elas}}{\beta_{\rm visc}}+\lambda\right)\int_0^t\exp(-\lambda(t-s))\dot{\overline e}(s)\,\mathrm ds.
\end{equation}
For each fixed $\overline e$ and $t$, every tensor component of $\cW_{\overline e,t}(\lambda)$ is an entire function of $\lambda$. If $\operatorname{Re}\lambda>0$, then
\begin{equation}\label{eq:weighted_response_bound}
\sup_{\overline e\in\cE}\sup_{t\ge0}\|\cW_{\overline e,t}(\lambda)\|_F\le\left(1+\frac{|\lambda|}{\operatorname{Re}\lambda}\right)\frac{\delta_{\rm stress}}{\beta_{\rm visc}}.
\end{equation}
Here the Frobenius norm for complex tensors is $\|Z\|_F=(\sum_{i,j=1}^d|Z_{ij}|^2)^{1/2}$.
\end{lemma}

\begin{proof}[\bf Proof]
For each $t\ge0$, differentiation under the integral on the finite interval $[0,t]$ proves analyticity. Write $r=\operatorname{Re}\lambda>0$, $c=|\lambda|/r\ge1$, and $z(t)=\int_0^t\exp(-\lambda(t-s))\dot{\overline e}(s)\,\mathrm ds$. Integration by parts gives
\begin{equation*}
z(t)=\overline e(t)-\exp(-\lambda t)\overline e(0)-\lambda\int_0^t\exp(-\lambda(t-s))\overline e(s)\,\mathrm ds.
\end{equation*}
Consequently, $\|z(t)\|_F\le[1+\exp(-rt)+c(1-\exp(-rt))]\|\overline e\|_{L^\infty}\le(1+c)\|\overline e\|_{L^\infty}$. By the definition of $z$, we also have $|\lambda|\|z(t)\|_F\le c\|\dot{\overline e}\|_{L^\infty}$. Hence
\begin{align*}
\|\cW_{\overline e,t}(\lambda)\|_F
\le \frac{\beta_{\rm elas}}{\beta_{\rm visc}}(1+c)\|\overline e\|_{L^\infty}+c\|\dot{\overline e}\|_{L^\infty}\le\frac{1+c}{\beta_{\rm visc}}\left(\beta_{\rm elas}\|\overline e\|_{L^\infty}+\beta_{\rm visc}\|\dot{\overline e}\|_{L^\infty}\right).
\end{align*}
The definition of $\cE$ gives \cref{eq:weighted_response_bound}.
\end{proof}

\begin{definition}[Prony Series by Local Chebyshev Interpolation]\label{def:prony_series_approximation}
A Prony-series approximation of $\bK(\dummy;\bV,\bE)$ is a kernel
\begin{equation}\label{eq:Kn_def}
\bK^{(n)}(t;\bV,\bE)\coloneqq\sum_{k=1}^n\bK_k(\bV,\bE)\exp(-\lambda_k t),
\end{equation}
with microstructure-independent rates $\lambda_k\in\cI$. If $\lambda_{\min}=\lambda_{\max}$, set $n=1$, $\lambda_1=\lambda_{\min}$, and
\begin{equation}\label{eq:Ki_def_degenerate}
\bK_1(\bV,\bE)\coloneqq\int_{\cI}\lambda\,\mathrm d\nu_{\bV,\bE}(\lambda).
\end{equation}
Otherwise, set $\mathsf R=\ln(\lambda_{\max}/\lambda_{\min})$, $L=\lceil\mathsf R\rceil$, $h=\mathsf R/L\le1$, and $J=L+1$. For $j=1,\dots,J$, define
\begin{gather*}
\tau_j=(j-1)h,\qquad a_j=\max\{0,\tau_j-h\},\qquad b_j=\min\{\mathsf R,\tau_j+h\},\\
\theta_j(\lambda)=\max\left\{1-\frac{|\ln(\lambda/\lambda_{\min})-\tau_j|}{h},0\right\},\qquad\lambda\in\cI.
\end{gather*}
The functions $\theta_j$ are continuous and nonnegative, satisfy $\sum_{j=1}^J\theta_j=1$ on $\cI$, and are supported on $\cI_j=[\lambda_{\min}\exp(a_j),\lambda_{\min}\exp(b_j)]$. For the reference coordinate $x\in[-1,1]$, define the local changes of variables by
\begin{equation}\label{eq:log_spectral_coordinate}
g_j(x)=\lambda_{\min}\exp\left(\frac{a_j+b_j}{2}+\frac{b_j-a_j}{2}x\right),\qquad
\psi_j(\lambda)=\frac{2\ln(\lambda/\lambda_{\min})-a_j-b_j}{b_j-a_j}.
\end{equation}
For $m\in\bN$, let $x_i=\cos((2i-1)\pi/(2m))$, $i=1,\dots,m$, be the Chebyshev nodes, and let $\ell_i$ be the associated Lagrange basis polynomials. Set
\begin{equation}\label{eq:Ki_def}
\lambda_{ji}=g_j(x_i),\qquad
\bK_{ji}(\bV,\bE)=\left(\frac{\beta_{\rm elas}}{\beta_{\rm visc}}+\lambda_{ji}\right)\int_{\cI}\theta_j(\lambda)\ell_i(\psi_j(\lambda))\,\mathrm d\mu_{\bV,\bE}(\lambda).
\end{equation}
Use the single index $k=(j-1)m+i$, so that $n=Jm$, $\lambda_k=\lambda_{ji}$, and $\bK_k=\bK_{ji}$.
\end{definition}

\begin{lemma}[Local Approximation of the Weighted Response]\label{lem:uniform_prony_discretization_beta}
Suppose $\lambda_{\min}<\lambda_{\max}$, and use the response in \cref{eq:weighted_response_def} and the maps and nodes in \cref{def:prony_series_approximation}. For $\overline e\in\cE$, $t\ge0$, and $j=1,\dots,J$, set $f_j(x)=\cW_{\overline e,t}(g_j(x))$. The interpolant $p_j(x)=\sum_{i=1}^m\ell_i(x)f_j(x_i)$ has degree at most $m-1$ and satisfies
\begin{equation}\label{eq:chebyshev_interp_bound_iv}
\sup_{\overline e\in\cE}\sup_{t\ge0}\max_{1\le j\le J}\sup_{x\in[-1,1]}\|f_j(x)-p_j(x)\|_F
\le\frac{4\delta_{\rm stress}}{\beta_{\rm visc}}\left(1+\sec\frac34\right)2^{-(m-1)}.
\end{equation}
\end{lemma}

\begin{proof}[\bf Proof]
Let $\cB_2$ be the closed Bernstein ellipse with parameter $2$. Its imaginary semiaxis is $3/4$. Since $(b_j-a_j)/2\le h\le1$, \cref{eq:log_spectral_coordinate} gives $|\arg g_j(z)|\le3/4$ on $\cB_2$. In particular, $\operatorname{Re}g_j(z)>0$ and $|g_j(z)|/\operatorname{Re}g_j(z)\le\sec(3/4)$. By \cref{lem:weighted_response}, $f_j$ is analytic on a neighborhood of $\cB_2$ and
\begin{equation*}
\sup_{z\in\cB_2}\|f_j(z)\|_F\le\frac{\delta_{\rm stress}}{\beta_{\rm visc}}\left(1+\sec\frac34\right).
\end{equation*}
Denote this bound by $M_f$. For any complex linear functional $\varphi$ of norm one, the Chebyshev coefficients $c_k$ of $\varphi\circ f_j$ satisfy $|c_k|\le2M_f2^{-k}$ \cite[Theorem 8.1]{trefethen2019approximation}. Let $I_m$ denote interpolation at the roots of $T_m$. Each Chebyshev polynomial $T_k$ interpolates to zero or to $\pm T_r$ for some $r<m$, so $\|I_mT_k\|_{L^\infty(-1,1)}\le1$. Since interpolation is exact for $k<m$,
\begin{equation*}
\|\varphi\circ f_j-I_m(\varphi\circ f_j)\|_{L^\infty(-1,1)}
\le2\sum_{k=m}^{\infty}|c_k|
\le4M_f2^{-(m-1)}.
\end{equation*}
Taking the supremum over $\varphi$ proves \cref{eq:chebyshev_interp_bound_iv}.
\end{proof}

\begin{proposition}[Uniform Approximation by the Prony-Series Internal-Variable Model]\label{prop:uniform_prony_homogenized_model}
Let the rates and coefficients be given by \cref{def:prony_series_approximation}. Each coefficient $\bK_k$ depends continuously on $(\bV,\bE)\in\cM$ in the $L^1_{\rm per}$ topology. For $(\bV,\bE)\in\cM$ and $\overline e\in\cE$, define $q_k\in W^{1,\infty}(\cT;\cV^2)$ by
\begin{equation}\label{eq:q_i_ode}
\dot q_k(t)+\lambda_k q_k(t)=\bK_k(\bV,\bE)\dot{\overline e}(t),\qquad q_k(0)=0,\qquad k=1,\dots,n,
\end{equation}
and set
\begin{equation}\label{eq:prony_homogenized_model}
\overline\sigma_{\rm IV}(t)=\bV^{\rm hom}(\bV)\dot{\overline e}(t)+\bE^{\rm hom}(\bE)\overline e(t)+\sum_{k=1}^n q_k(t).
\end{equation}
Denote the resulting constitutive model by $\cH_{\rm IV}$. This is the homogenized constitutive law with $\bK$ replaced by $\bK^{(n)}$. If $\lambda_{\min}=\lambda_{\max}$, the approximation is exact with $n=1$. Otherwise, with $n=Jm$,
\begin{align}\label{eq:prony_homogenized_error_uniform}
\sup_{(\bV,\bE)\in\cM}\sup_{\overline e\in\cE}\esssup_{t\in\cT}\|\cH(\{\overline e(s)\}_{s\in[0,t]};\bV,\bE)-\overline\sigma_{\rm IV}(t)\|_F
&<\kappa_{\rm IV}2^{-(m-1)},\\[-1ex]
\kappa_{\rm IV}&\coloneqq8d_{\cV^2}\left(1+\sec\frac34\right)\delta_{\rm stress}.\nonumber
\end{align}
\end{proposition}

\begin{proof}[\bf Proof]
For a singleton interval, coefficient continuity follows from \cref{eq:Ki_def_degenerate,lem:weak_convergence}. Otherwise, each function $\lambda\mapsto\theta_j(\lambda)\ell_i(\psi_j(\lambda))$ is continuous on $\cI$, and the rates are independent of $(\bV,\bE)$. Coefficient continuity then follows from \cref{eq:Ki_def,lem:weighted_spectral_measure}.

Duhamel's formula gives
\begin{equation}\label{eq:q_i_duhamel}
q_k(t)=\int_0^t\exp(-\lambda_k(t-s))\bK_k(\bV,\bE)\dot{\overline e}(s)\,\mathrm ds.
\end{equation}
Since $\lambda_k>0$ and $\dot{\overline e}\in L^\infty$, this solution is bounded, and \cref{eq:q_i_ode} bounds its derivative. Summing \cref{eq:q_i_duhamel} and using \cref{eq:Kn_def} proves the asserted constitutive representation. If $\cI$ is a singleton, \cref{eq:Ki_def_degenerate,eq:K_laplace} give exactness.

For a nondegenerate interval, abbreviate $\mu=\mu_{\bV,\bE}$. By \cref{eq:weighted_measure_def,eq:K_laplace}, the memory contribution to the stress equals $\int_{\cI}\mathrm d\mu(\lambda)\cW_{\overline e,t}(\lambda)$, where the tensor measure acts on the vector-valued integrand. The coefficients in \cref{eq:Ki_def} give the corresponding contribution in the Prony-series model, $\sum_{j=1}^J\int_{\cI}\theta_j(\lambda)\,\mathrm d\mu(\lambda)p_j(\psi_j(\lambda))$. Hence their difference is
\begin{equation*}
\sum_{j=1}^J\int_{\cI}\theta_j(\lambda)\,\mathrm d\mu(\lambda)\bigl[f_j(\psi_j(\lambda))-p_j(\psi_j(\lambda))\bigr].
\end{equation*}
The interpolation bound applies on the support of $\theta_j$, where $\psi_j(\lambda)\in[-1,1]$. If $|\mu|$ denotes the total variation measure, the norm of this difference is at most
\begin{align*}
\frac{4\delta_{\rm stress}}{\beta_{\rm visc}}\left(1+\sec\frac34\right)2^{-(m-1)}\sum_{j=1}^J\int_{\cI}\theta_j(\lambda)\,\mathrm d|\mu|(\lambda)
&=\frac{4\delta_{\rm stress}}{\beta_{\rm visc}}\left(1+\sec\frac34\right)2^{-(m-1)}\|\mu\|_{\mathrm{TV}(\cI)}\\
&<\kappa_{\rm IV}2^{-(m-1)},
\end{align*}
by the partition of unity and the uniform mass bound \cref{eq:weighted_measure_mass}. The instantaneous terms coincide, so this proves \cref{eq:prony_homogenized_error_uniform}.
\end{proof}

}

{
We now approximate the stress and kinetic responses of the Prony-series model and bound the resulting errors in the internal variables and stress, as required in \ref{item:ua_proof_s3}.

\begin{lemma}[Uniform Bound for the Internal Variable]\label{lem:ua_step3_targets}
{Choose an approximation from \cref{def:prony_series_approximation} with $n$ terms}, and let $\xi_{\rm IV}\coloneqq\mathsf{vec}_2^{(n)}(q_1,\dots,q_n)\in\bR^{nd_{\cV^2}}$, where the $q_k$ solve \cref{eq:q_i_ode} and $\mathsf{vec}_2^{(n)}$ concatenates their Mandel vectors. There is a constant $\delta_{\xi_{\rm IV}}>0$ such that
\begin{equation*}
\sup_{(\bV,\bE)\in\cM}\sup_{\overline e\in\cE}\|\xi_{\rm IV}\|_{L^\infty}\le\delta_{\xi_{\rm IV}}.
\end{equation*}
\end{lemma}

\begin{proof}[\bf Proof]
By \cref{eq:q_i_duhamel}, $\lambda_k\ge\lambda_{\min}>0$, and \cref{eq:admissible_strain},
\begin{equation*}
\|q_k(t)\|_F\le\frac{\delta_{\rm stress}}{\beta_{\rm visc}\lambda_{\min}}\|\bK_k(\bV,\bE)\|_F,\qquad k=1,\dots,n.
\end{equation*}
{Continuity of each coefficient and compactness of $\cM$ (\cref{prop:uniform_prony_homogenized_model,lem:bv_compact}) give $C_k\coloneqq\sup_{(\bV,\bE)\in\cM}\|\bK_k(\bV,\bE)\|_F<\infty$. Since $\mathsf{vec}_2^{(n)}$ is an isometry, one may take
\begin{equation*}
\delta_{\xi_{\rm IV}}=1+\frac{\delta_{\rm stress}}{\beta_{\rm visc}\lambda_{\min}}\left(\sum_{k=1}^n C_k^2\right)^{1/2}.
\end{equation*}}
\end{proof}
}

\begin{lemma}[Uniform Approximation by Neural Operator Response Mappings]\label{lem:ua_step3_no}
{{Consider an approximation from \cref{def:prony_series_approximation} with $n$ terms}. For $\eta_{\xi}>0$, let} $\delta_{\xi_{\mathrm{IV}}}>0$ be given by \cref{lem:ua_step3_targets}. For $\xi\in\bR^{nd_{\cV^2}}$, write $\xi=\mathsf{vec}_2^{(n)}(q_1,\dots,q_n)$ with $q_1,\dots,q_n\in\cV^2$, and define the stress and kinetic responses associated with the Prony-series internal-variable model in \cref{prop:uniform_prony_homogenized_model} by
\begin{align*}
\cF_{\mathrm{IV}}(\overline{e},\dot{\overline{e}},\xi;\bV,\bE)
&\coloneqq
\bV^{\rm hom}(\bV)\dot{\overline{e}}
+\bE^{\rm hom}(\bE)\overline{e}
+\sum_{i=1}^{n} q_i,\\
\cG_{\mathrm{IV}}(\dot{\overline{e}},\xi;\bV,\bE)
&\coloneqq
\mathsf{vec}_2^{(n)}\big(\bK_1(\bV,\bE)\dot{\overline{e}}-\lambda_1 q_1,\dots,
\bK_n(\bV,\bE)\dot{\overline{e}}-\lambda_n q_n\big).
\end{align*}
{Define the compact input set
\begin{equation*}
\begin{aligned}
\cK_{\rm vec}\coloneqq\bigl\{(\overline e,\dot{\overline e},\xi)\in\cV^2\times\cV^2\times\bR^{nd_{\cV^2}}\ \big|\ &\beta_{\rm elas}\|\overline e\|_F+\beta_{\rm visc}\|\dot{\overline e}\|_F\le\delta_{\rm stress},\|\xi\|\le\delta_{\xi_{\rm IV}}+\eta_\xi\bigr\}.
\end{aligned}
\end{equation*}
}
Then, for any $\varepsilon_{\cF},\varepsilon_{\cG}>0$, there exist response mappings $(\cF_{\mathrm{NO}},\cG_{\mathrm{NO}})\in \mathfrak{F}(nd_{\cV^2})\times \mathfrak{G}(nd_{\cV^2})$ such that
    \begin{align}
    \sup_{\substack{(\overline e,\dot{\overline e},\xi)\in \cK_{\mathrm{vec}}\\(\bV,\bE)\in \cM}}
    \left\|\cF_{\mathrm{NO}}(\overline e,\dot{\overline e},\xi;\bV,\bE)-\cF_{\mathrm{IV}}(\overline e,\dot{\overline e},\xi;\bV,\bE)\right\|_F
    &\leq \varepsilon_{\cF},\label{eq:ua_full_F_approx}\\
    \sup_{\substack{(\overline e,\dot{\overline e},\xi)\in \cK_{\mathrm{vec}}\\(\bV,\bE)\in \cM}}
    \left\|\cG_{\mathrm{NO}}(\dot{\overline e},\xi;\bV,\bE)-\cG_{\mathrm{IV}}(\dot{\overline e},\xi;\bV,\bE)\right\|
    &\leq \varepsilon_{\cG}.\label{eq:ua_full_G_approx}
    \end{align}
Moreover, for every $(\bV,\bE)\in\cM$, the map $(\dot{\overline e},\xi)\mapsto \cG_{\mathrm{NO}}(\dot{\overline e},\xi;\bV,\bE)$ is locally Lipschitz.
\end{lemma}

\begin{proof}[\bf Proof]
The set $\cK_{\mathrm{vec}}$ is compact by definition, and \cref{lem:bv_compact} gives the compactness of $\cM$. By \cref{cor:continuity_tensors,prop:uniform_prony_homogenized_model}, the maps $(\bV,\bE)\mapsto \bV^{\rm hom}(\bV)$, $(\bV,\bE)\mapsto \bE^{\rm hom}(\bE)$, and $(\bV,\bE)\mapsto \bK_i(\bV,\bE)$ are continuous on $\cM$, so $\cF_{\mathrm{IV}}$ and $\cG_{\mathrm{IV}}$ are continuous on $\cK_{\mathrm{vec}}\times\cM$.

Given $\varepsilon_{\cF},\varepsilon_{\cG}>0$, define $\widehat{\cK}_{\cF}\coloneqq
\{(\mathsf{vec}_2\overline e,\mathsf{vec}_2\dot{\overline e},\xi)\mid (\overline e,\dot{\overline e},\xi)\in \cK_{\mathrm{vec}}\}$,
$\widehat{\cK}_{\cG}
\coloneqq
\{(\mathsf{vec}_2\dot{\overline e},\xi)\mid (\overline e,\dot{\overline e},\xi)\in \cK_{\mathrm{vec}}\}$, and $\widehat{\cM}\coloneqq
\{(\mathsf{vec}_4\circ\bV,\mathsf{vec}_4\circ\bE)\mid (\bV,\bE)\in\cM\}$.
These sets are compact as continuous images of $\cK_{\mathrm{vec}}$ and $\cM$ under the relevant vectorization maps. Define
\begin{align*}
g_{\cF}\Bigl((\mathsf{vec}_2\overline e,\mathsf{vec}_2\dot{\overline e},\xi),(\mathsf{vec}_4\circ\bV,\mathsf{vec}_4\circ\bE)\Bigr)
&\coloneqq
\mathsf{vec}_2\bigl(\cF_{\mathrm{IV}}(\overline e,\dot{\overline e},\xi;\bV,\bE)\bigr),\\
g_{\cG}\Bigl((\mathsf{vec}_2\dot{\overline e},\xi),(\mathsf{vec}_4\circ\bV,\mathsf{vec}_4\circ\bE)\Bigr)
&\coloneqq
\cG_{\mathrm{IV}}(\dot{\overline e},\xi;\bV,\bE).
\end{align*}
{Then $g_{\cF}$ and $g_{\cG}$ are continuous on $\widehat{\cK}_{\cF}\times\widehat{\cM}$ and $\widehat{\cK}_{\cG}\times\widehat{\cM}$, respectively. Apply~\ref{no:ua} to $g_{\cF}$ and~\ref{no:lip} to $g_{\cG}$ to obtain maps $\Phi_{\cF}$ and $\Phi_{\cG}$ satisfying
\begin{equation*}
\sup_{\widehat{\cK}_{\cF}\times\widehat{\cM}}\|\Phi_{\cF}-g_{\cF}\|\le\varepsilon_{\cF},\qquad
\sup_{\widehat{\cK}_{\cG}\times\widehat{\cM}}\|\Phi_{\cG}-g_{\cG}\|\le\varepsilon_{\cG},
\end{equation*}
with $\Phi_{\cG}$ locally Lipschitz in its vector input for each microstructure. Let $(\cF_{\rm NO},\cG_{\rm NO})\in\mathfrak F(nd_{\cV^2})\times\mathfrak G(nd_{\cV^2})$ be the induced response mappings. The isometry of $\mathsf{vec}_2$ gives \cref{eq:ua_full_F_approx,eq:ua_full_G_approx} and preserves the local Lipschitz property of $\cG_{\rm NO}$.}
\end{proof}

{\begin{proposition}[Global Existence and Uniform Error Bounds]\label{prop:ua_step3_stability}
{Consider an approximation from \cref{def:prony_series_approximation} with $n$ terms}. Let $\varepsilon_{\cF},\eta_\xi>0$, let $0<\varepsilon_{\cG}\le\lambda_{\min}\eta_\xi$, and choose $(\cF_{\rm NO},\cG_{\rm NO})$ as in \cref{lem:ua_step3_no}. For every $(\bV,\bE)\in\cM$ and $\overline e\in\cE$, the initial value problem
\begin{equation*}
\dot\xi_{\rm NO}(t)=\cG_{\rm NO}(\dot{\overline e}(t),\xi_{\rm NO}(t);\bV,\bE)\quad\text{for a.e.\ }t\in\cT,\qquad\xi_{\rm NO}(0)=0,
\end{equation*}
{admits a unique solution $\xi_{\rm NO}\in W^{1,\infty}(\cT;\bR^{nd_{\cV^2}})$}. Let $\overline\sigma_{\rm NO}$ be its stress response in \cref{eq:learned_model_full}. Then
\begin{equation}\label{eq:ua_err}
\begin{aligned}
\sup_{(\bV,\bE)\in\cM}\sup_{\overline e\in\cE}\|\xi_{\rm NO}-\xi_{\rm IV}\|_{L^\infty}&\le\frac{\varepsilon_{\cG}}{\lambda_{\min}},\\
\sup_{(\bV,\bE)\in\cM}\sup_{\overline e\in\cE}\|\overline\sigma_{\rm NO}-\overline\sigma_{\rm IV}\|_{L^\infty}&\le\varepsilon_{\cF}+\frac{\sqrt n}{\lambda_{\min}}\varepsilon_{\cG},
\end{aligned}
\end{equation}
where $\xi_{\rm IV}$ is the internal variable defined in \cref{lem:ua_step3_targets} and $\overline\sigma_{\rm IV}$ is given by \cref{eq:prony_homogenized_model}.
\end{proposition}}

\begin{proof}[\bf Proof]
Let $(\bV,\bE)\in\cM$ and $\overline e\in\cE$. For $j=1,2$, write $\xi_j=\mathsf{vec}_2^{(n)}(q_{j,1},\dots,q_{j,n})$. The definition of $\cG_{\mathrm{IV}}$ in \cref{lem:ua_step3_no} gives
\begin{equation*}
\cG_{\mathrm{IV}}(\dot{\overline e},\xi_1;\bV,\bE)-\cG_{\mathrm{IV}}(\dot{\overline e},\xi_2;\bV,\bE)
=\mathsf{vec}_2^{(n)}\big(-\lambda_1(q_{1,1}-q_{2,1}),\dots,-\lambda_n(q_{1,n}-q_{2,n})\big).
\end{equation*}
Since $\mathsf{vec}_2^{(n)}$ is an isometry,
\begin{equation}\label{eq:G_IV_dissipative}
\big(\xi_1-\xi_2\big)^{\top}\Big(\cG_{\mathrm{IV}}(\dot{\overline e},\xi_1;\bV,\bE)-\cG_{\mathrm{IV}}(\dot{\overline e},\xi_2;\bV,\bE)\Big)
\le -\lambda_{\min}\|\xi_1-\xi_2\|^2.
\end{equation}

By \cref{lem:ua_step3_no}, the map $(\dot{\overline e},\xi)\mapsto \cG_{\mathrm{NO}}(\dot{\overline e},\xi;\bV,\bE)$ is locally Lipschitz. Since $\dot{\overline e}\in L^\infty(\cT;\cV^2)$, Carath\'eodory's theorem gives a unique maximal {locally} absolutely continuous solution $\xi_{\mathrm{NO}}$ on $[0,T_{\max})$, where $T_{\max}$ denotes its maximal existence time. Set $\zeta(t)\coloneqq \xi_{\mathrm{NO}}(t)-\xi_{\mathrm{IV}}(t)$ and define
\begin{equation*}
t_*\coloneqq \sup\left\{{0\le\tau<T_{\max}}\ \middle|\ \|\zeta(t)\|\le \eta_\xi\ \text{for all }t\in[0,\tau]\right\}.
\end{equation*}
Since $\zeta(0)=0$, we have $t_*>0$. On {$[0,t_*)$}, \cref{lem:ua_step3_targets} and the definition of $\cE$ give $(\overline e(t),\dot{\overline e}(t),\xi_{\mathrm{NO}}(t))\in\cK_{\mathrm{vec}}$ for a.e.\ $t$. Hence \cref{eq:ua_full_G_approx} applies there. For a.e.\ {$t\in[0,t_*)$},
\begin{align*}
\dot\zeta(t)
&=
\Big(\cG_{\mathrm{NO}}(\dot{\overline e}(t),\xi_{\mathrm{NO}}(t);\bV,\bE)
-\cG_{\mathrm{IV}}(\dot{\overline e}(t),\xi_{\mathrm{NO}}(t);\bV,\bE)\Big)\\
&\quad+
\Big(\cG_{\mathrm{IV}}(\dot{\overline e}(t),\xi_{\mathrm{NO}}(t);\bV,\bE)
-\cG_{\mathrm{IV}}(\dot{\overline e}(t),\xi_{\mathrm{IV}}(t);\bV,\bE)\Big).
\end{align*}
Taking the Euclidean inner product with $\zeta(t)$ and using \cref{eq:ua_full_G_approx,eq:G_IV_dissipative} gives
\begin{equation*}
\frac{1}{2}\frac{\mathrm{d}}{\mathrm{d}t}\|\zeta(t)\|^2
\le \varepsilon_{\cG}\|\zeta(t)\| - \lambda_{\min}\|\zeta(t)\|^2.
\end{equation*}
Therefore, for a.e.\ {$t\in[0,t_*)$} such that $\|\zeta(t)\|>0$, $\frac{\mathrm{d}}{\mathrm{d}t}\|\zeta(t)\|\le \varepsilon_{\cG}-\lambda_{\min}\|\zeta(t)\|$. Since $\zeta(0)=0$, comparison yields
\begin{equation}\label{eq:zeta_uniform_bound}
\|\zeta(t)\|\le \lambda_{\min}^{-1}\varepsilon_{\cG}\big(1-\exp(-\lambda_{\min}t)\big)\le \lambda_{\min}^{-1}\varepsilon_{\cG}
\end{equation}
on {$[0,t_*)$. If $t_*<T_{\max}$, continuity and $\varepsilon_{\cG}\le\lambda_{\min}\eta_\xi$ give $\|\zeta(t_*)\|<\eta_\xi$, contradicting the definition of $t_*$. Thus $t_*=T_{\max}$. Since $\xi_{\mathrm{NO}}$ remains bounded by $\delta_{\xi_{\mathrm{IV}}}+\eta_\xi$ and $\dot{\overline e}$ is bounded, the continuation criterion for Carath\'eodory ODEs gives $T_{\max}=\infty$. Moreover, \cref{eq:ua_full_G_approx} and the boundedness of $\cG_{\mathrm{IV}}$ on $\cK_{\mathrm{vec}}\times\cM$ imply $\dot\xi_{\mathrm{NO}}\in L^\infty$, so $\xi_{\mathrm{NO}}\in W^{1,\infty}(\cT;\bR^{nd_{\cV^2}})$.} Hence \cref{eq:zeta_uniform_bound} holds on $\cT$, and $\|\xi_{\mathrm{NO}}-\xi_{\mathrm{IV}}\|_{L^\infty}\le \lambda_{\min}^{-1}\varepsilon_{\cG}$.

With $(\overline e(t),\dot{\overline e}(t),\xi_{\mathrm{NO}}(t))\in\cK_{\mathrm{vec}}$ for a.e.\ $t\in\cT$, we apply \cref{eq:ua_full_F_approx}. By the definition of $\cF_{\mathrm{IV}}$ in \cref{lem:ua_step3_no},
\begin{gather*}
{\overline\sigma_{\rm IV}}(t)=\cF_{\mathrm{IV}}(\overline e(t),\dot{\overline e}(t),\xi_{\mathrm{IV}}(t);\bV,\bE)
\quad \text{for a.e.\ } t\in\cT.\\
\|\overline{\sigma}_{\mathrm{NO}}(t)-{\overline\sigma_{\rm IV}}(t)\|_F
\leq
\varepsilon_{\cF}+
\left\|\cF_{\mathrm{IV}}(\overline e(t),\dot{\overline e}(t),\xi_{\mathrm{NO}}(t);\bV,\bE)-\cF_{\mathrm{IV}}(\overline e(t),\dot{\overline e}(t),\xi_{\mathrm{IV}}(t);\bV,\bE)\right\|_F.
\end{gather*}
Since $\cF_{\mathrm{IV}}$ is affine in $\xi$ and sums the $n$ tensorial blocks, Cauchy--Schwarz gives
\begin{equation*}
\|\cF_{\mathrm{IV}}(\overline e(t),\dot{\overline e}(t),\xi_{\mathrm{NO}}(t);\bV,\bE)-\cF_{\mathrm{IV}}(\overline e(t),\dot{\overline e}(t),\xi_{\mathrm{IV}}(t);\bV,\bE)\|_F
\le \sqrt{n} \|\xi_{\mathrm{NO}}(t)-\xi_{\mathrm{IV}}(t)\|.
\end{equation*}
Therefore, $\|\overline{\sigma}_{\mathrm{NO}}(t)-{\overline\sigma_{\rm IV}}(t)\|_F\leq \varepsilon_{\cF} + \sqrt{n}\lambda_{\min}^{-1}\varepsilon_{\cG}$ for a.e.\ $t\in\cT$. Taking suprema over $(\bV,\bE)\in\cM$ and $\overline e\in\cE$ gives \cref{eq:ua_err}.
\end{proof}

\begin{proof}[\bf Completion of the Proof of \cref{thm:ua}]
Let $\varepsilon>0$ and choose $\varepsilon_{\rm IV}\in(0,\varepsilon)$. By \cref{prop:uniform_prony_homogenized_model}, choose a Prony-series model with $n$ terms whose uniform stress error is at most $\varepsilon_{\rm IV}$. Let $\eta_\xi>0$ and set
\begin{equation}\label{eq:ua_error_tolerances}
\varepsilon_{\cF}\coloneqq\frac{\varepsilon-\varepsilon_{\rm IV}}{2},\qquad
\varepsilon_{\cG}\coloneqq\min\left\{\frac{\lambda_{\min}\eta_\xi}{2},\frac{\lambda_{\min}(\varepsilon-\varepsilon_{\rm IV})}{2\sqrt n}\right\}.
\end{equation}
Choose the response mappings $(\cF_{\rm NO},\cG_{\rm NO})$ by \cref{lem:ua_step3_no}. These tolerances satisfy the hypothesis of \cref{prop:ua_step3_stability} and $\varepsilon_{\cF}+\sqrt n\varepsilon_{\cG}/\lambda_{\min}\le\varepsilon-\varepsilon_{\rm IV}$, so the neural operator model approximates the Prony-series model with uniform stress error at most $\varepsilon-\varepsilon_{\rm IV}$. The triangle inequality gives \cref{eq:ua} with $n_\xi=nd_{\cV^2}$.
\end{proof}

\section{Implementation Details for Viscoelasticity}\label{app:kv2d_details}

This appendix provides additional information on the numerical results for learning viscoelasticity.

\subsection{Numerical Solutions of the Unit Cell Problem}\label{app:kv2d_cell_solver}
The cell problem is discretized using a $128\times128$ uniform triangular mesh with $\bP_1$ finite elements and periodic boundary conditions. The homogenization solver uses a mixed formulation for the displacement corrector together with a Lagrange multiplier enforcing a zero-mean constraint. We use an explicit fourth-order Runge--Kutta scheme with $1024$ time steps on {$[0,1]$}. The displacement field at each time step is obtained via a direct solver. The solver is implemented in FEniCS \cite{logg2012automated, alnaes2015fenics}.

\subsection{Trajectory Sampling}\label{app:kv2d_strain_sampling}
For each independent component $a\in\{11,22,12\}$ of the macroscopic strain, we first sample the total number of control times $N_c$ uniformly from $\{3,\dots,21\}$. Let $0=t_0<t_1<\cdots<t_{N_c-1}=1$ be the control times, with the interior points sampled uniformly in $[0,1]$ and then sorted. The strain is initialized by $\overline e_a(t_0)=0$. For each $k=1,\dots,N_c-1$, we sample an independent sign $v_{a,k}\in\{-1,1\}$ and set $
    \overline e_a(t_k)
    =
    \overline e_a(t_{k-1})
    +
    \left(e_{\mathrm{max}}v_{a,k}-\overline e_a(t_{k-1})\right)
    \sqrt{t_k-t_{k-1}},$ with $e_{\mathrm{max}}=0.2$.
The values between control times are obtained via piecewise cubic Hermite interpolating polynomials (PCHIP). This procedure is applied independently to the three strain components.

\subsection{Architecture and Training}\label{app:kv2d_training}
For both the stress and kinetic response mappings in \cref{eq:F_and_G}, each Fourier neural mapping uses GELU activation, $4$ retained Fourier modes, a width of $32$, a projection width of $64$, and three layers. The temporal grids are downsampled to 256 steps for the training. The baseline dataset contains $4096$ training, $1024$ validation, and $2048$ test samples. We use the explicit Euler method to evolve the dynamics of the internal variable. Training uses the AdamW optimizer with a batch size of $32$, an initial learning rate of $10^{-3}$, and $400$ epochs. A cosine learning-rate scheduler decreases the learning rate to $10^{-5}$ by the final epoch. We use a penalty weight of $\alpha=1$ for $n_{\xi}=1$ and $\alpha=0.1$ for $n_{\xi}\in\{5, 10, 15\}$ and $\alpha=0.01$ for $n_{\xi}\in\{20,25\}$. The architecture and training are implemented in PyTorch \cite{paszke2019pytorch}.

\subsection{Voronoi Dataset}\label{app:kv2d_voronoi_dataset}
The baseline microstructures are piecewise constant on periodic Voronoi tessellations. The number of seeds is sampled uniformly from $\{5,\dots,15\}$, and the seed locations are sampled uniformly in the unit cell. Within each Voronoi region, the elastic and viscous Poisson ratios are sampled independently from $[0.05,0.45]$. The elastic and viscous Young's moduli are sampled independently from log-uniform distributions, with baseline ranges $[0.5,2.0]$ and $[0.02,0.2]$, respectively.

\subsection{Datasets with Varied Modulus Ranges}\label{app:kv2d_modulus_ranges}
For each material contrast shown in the middle panel of \cref{fig:ve_complexity}, we first generate a separate Voronoi dataset. To vary the elastic contrast range, we either set the elastic modulus to $1$ or sample it from the ranges $[0.71,1.41]$, $[0.35,2.83]$, or $[0.25,4]$. Note that even when the elastic modulus is constant, the elastic tensor remains heterogeneous because the Poisson ratio varies across each Voronoi region. To vary the viscous modulus range, we sample the viscous modulus from $[0.04,0.1]$, $[0.028,0.14]$, $[0.014,0.28]$, or $[0.01,0.4]$. Since the same Poisson ratio range is used throughout and the plane-stress isotropic tensor has smallest and largest eigenvalues $E/(1+\nu)$ and $E/(1-\nu)$, the material contrast in this numerical study is $\beta/\alpha=(E_{\max}/E_{\min})(1+\nu_{\max})/(1-\nu_{\max})$, where $\nu_{\max}=0.45$.

\subsection{Fitting a Linear Model}\label{app:kv2d_linear_fit}
For each microstructure, we evaluate the stress and kinetic responses of the trained neural operator model on the sampled trajectories and fit them by affine functions of the strain, strain rate, and internal state. After subtracting the fitted offsets, the linear part has the form
\begin{align*}
    \mathsf{vec}_2(\overline{\sigma}_{\rm NO})
    =A_{\dot e}\mathsf{vec}_2(\dot{\overline e})+A_e\mathsf{vec}_2(\overline e)+B\xi,\qquad\dot\xi=C\mathsf{vec}_2(\dot{\overline e})+D\xi,
\end{align*}
where $\xi\in\bR^{n_\xi}$ is the internal variable of the trained neural operator model. The tensorial form in \cref{eqns:learned_fitted_dynamics} follows from a change of coordinates. For $n_\xi=nd_{\cV^2}$ with $n\in\bN$, introduce tensors $q_{\rm NO}^{(i)}\in\cV^2$, $i=1,\dots,n$, and seek an invertible matrix $S$ such that
\begin{equation*}
    \xi=S\mathsf{vec}_2^{(n)}(q_{\rm NO}^{(1)},\dots,q_{\rm NO}^{(n)}),\qquad BS=[I\;\cdots\;I],
\end{equation*}
where each identity block has size $d_{\cV^2}\times d_{\cV^2}$. If $B$ has full row rank, such an $S$ exists, and the internal-variable contribution to the stress is $\sum_i q_{\rm NO}^{(i)}$. The matrices $A_{\dot e}$ and $A_e$ are the Mandel representations of $\bV^{\rm hom}_{\rm NO}$ and $\bE^{\rm hom}_{\rm NO}$. Partitioning $S^{-1}DS$ and $S^{-1}C$ into $d_{\cV^2}\times d_{\cV^2}$ blocks gives the Mandel representations of $\bL_{\rm NO}^{(ij)}$ and $\bK_{\rm NO}^{(i)}$, respectively. For $n_\xi=15$ in two dimensions, $d_{\cV^2}=3$ and hence $n=5$. Since $S^{-1}DS$ is similar to $D$, the eigenvalues plotted in \cref{fig:kv2d_iv_diagnostics} are invariants of the fitted linear dynamics.

\subsection{Particulate Dataset}\label{app:kv2d_particulate_dataset}

The particulate dataset consists of two-phase periodic microstructures with circular inclusions. In each sample, the elastic and viscous modulus and Poisson-ratio fields share the same inclusion geometry. The number of inclusions is drawn uniformly from $\{1,\dots,20\}$, and their sizes are chosen to give a particle area fraction in $[0.03,0.40]$ subject to non-overlap in the periodic cell. The material fields are constant on each phase. We sample elastic moduli log-uniformly from $[0.5,1.0]$ in the matrix and from $[1.0,2.0]$ in the inclusions. Viscous moduli use the reverse ordering on the logarithmic scale over $[0.02,0.2]$, so that the matrix phase is more viscous. We sample elastic and viscous Poisson ratios independently from $[0.05,0.45]$ in both phases.

\subsection{Median-Error Prediction Visualization}\label{app:kv2d_median_examples}
\begin{figure}[htp]
    \centering
\addtolength{\tabcolsep}{-5 pt}
    {\bf The median-error test sample: Voronoi dataset}\\
    
    \scalebox{0.85}{
    \begin{tabular}{>{\centering\arraybackslash} m{0.05\linewidth} >{\centering\arraybackslash} m{0.15\linewidth} >{\centering\arraybackslash} m{0.15\linewidth}}
       & Elasticity & Viscosity  \\
       \raisebox{5pt}{\rotatebox[origin=c]{90}{Modulus}}& \includegraphics[width=\linewidth]{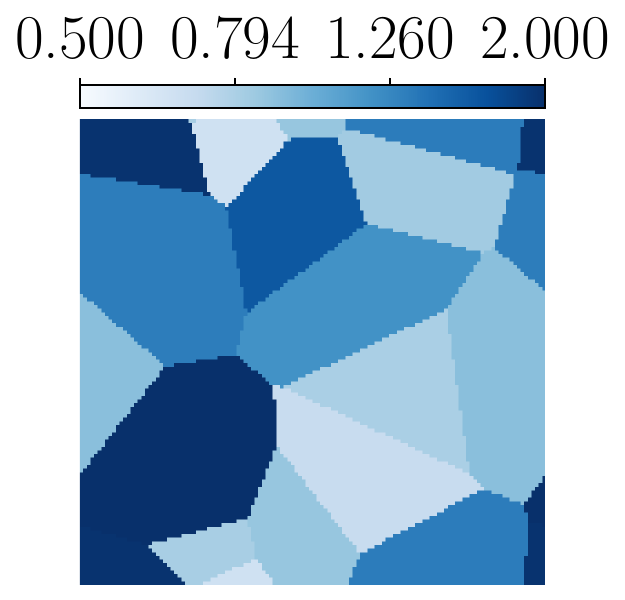}  & \includegraphics[width=\linewidth]{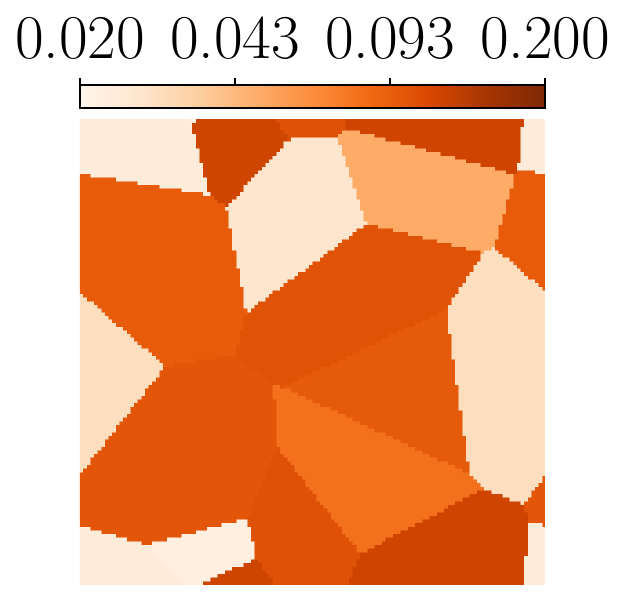} \\
       \raisebox{5pt}{\rotatebox[origin=c]{90}{\makecell{Poisson\\ratio}}}& \includegraphics[width=\linewidth]{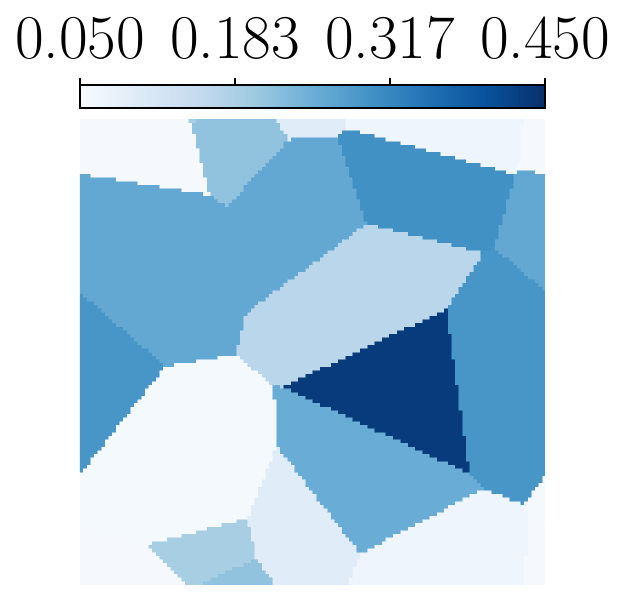} & \includegraphics[width=\linewidth]{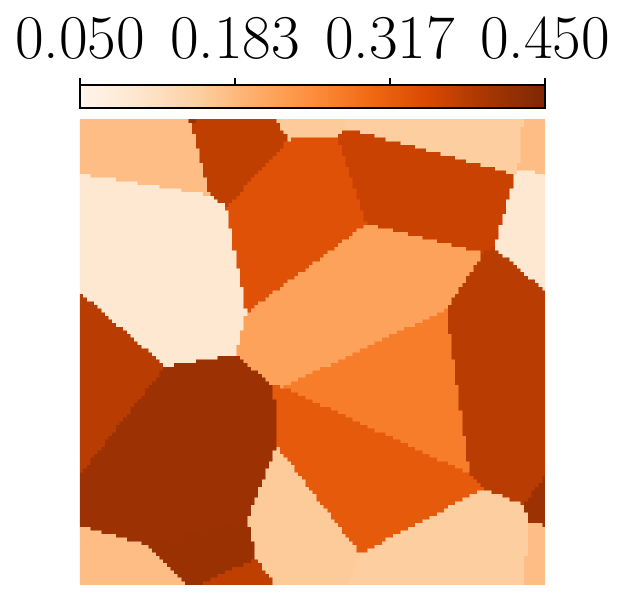}
    \end{tabular}
    }
    \begin{tabular}{c c}
        \includegraphics[width=0.245\linewidth, trim={0 0.9in 0 0},clip]{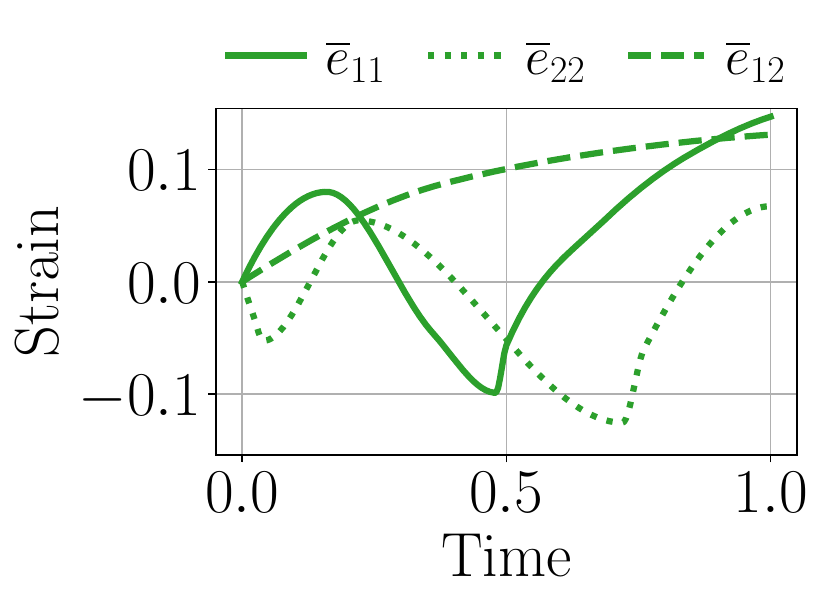} & \includegraphics[width=0.25\linewidth, trim={0 0.9in 0 0},clip]{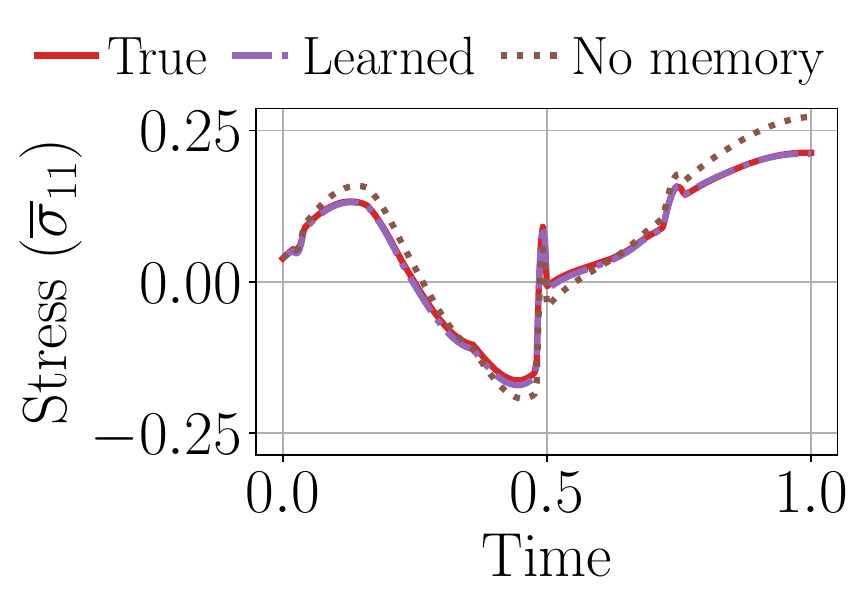}  \\
        \includegraphics[width=0.25\linewidth, trim={0 0 0 0.6in},clip]{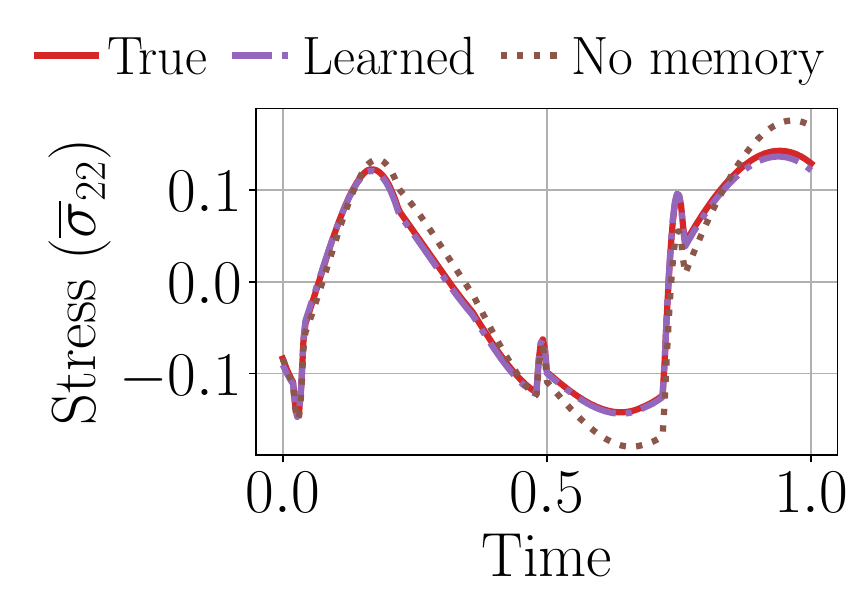} & \includegraphics[width=0.25\linewidth, trim={0 0 0 0.6in},clip]{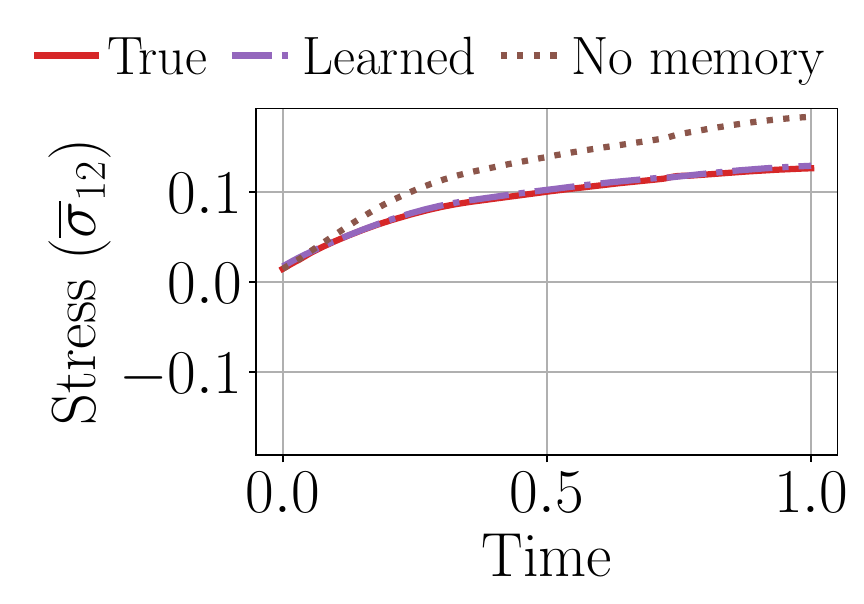}\\
        \multicolumn{2}{c}{\includegraphics[width=0.27\linewidth]{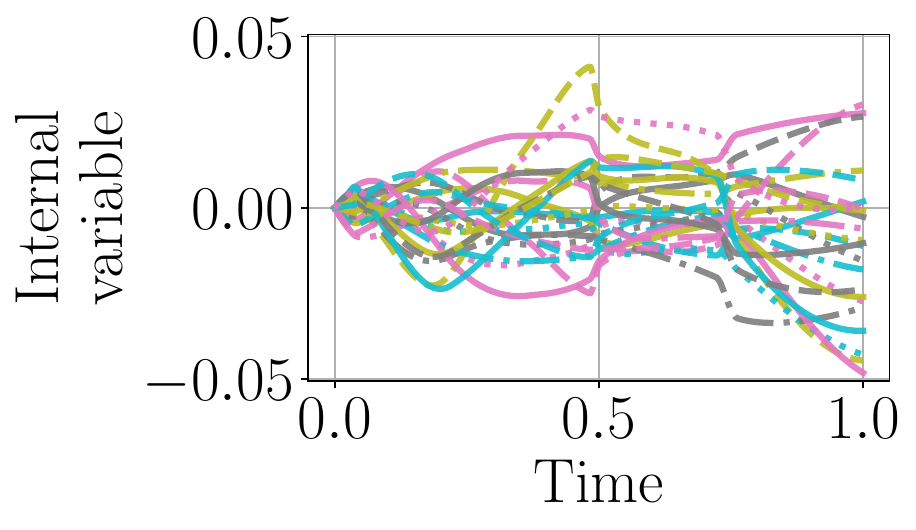}}
    \end{tabular}\\
\addtolength{\tabcolsep}{5 pt}
    \caption{Prediction by the trained neural operator model for a viscoelasticity sample with median relative stress error, using internal-variable dimension $n_{\xi}=25$. The left panels show the microstructural property fields, while the right panels show the prescribed strain trajectory, the predicted and reference stress trajectory, and the predicted internal-variable trajectory. The memoryless prediction ($\bK=0$) is included for comparison.}
    \label{fig:kv2d_prediction_median}
\end{figure}

\begin{figure}[htp]
    \centering
\addtolength{\tabcolsep}{-5 pt}
    {\bf The median-error test sample: Particulate dataset}\\
    
    \scalebox{0.85}{
    \begin{tabular}{>{\centering\arraybackslash} m{0.05\linewidth} >{\centering\arraybackslash} m{0.15\linewidth} >{\centering\arraybackslash} m{0.15\linewidth}}
       & Elasticity & Viscosity  \\
       \raisebox{5pt}{\rotatebox[origin=c]{90}{Modulus}}& \includegraphics[width=\linewidth]{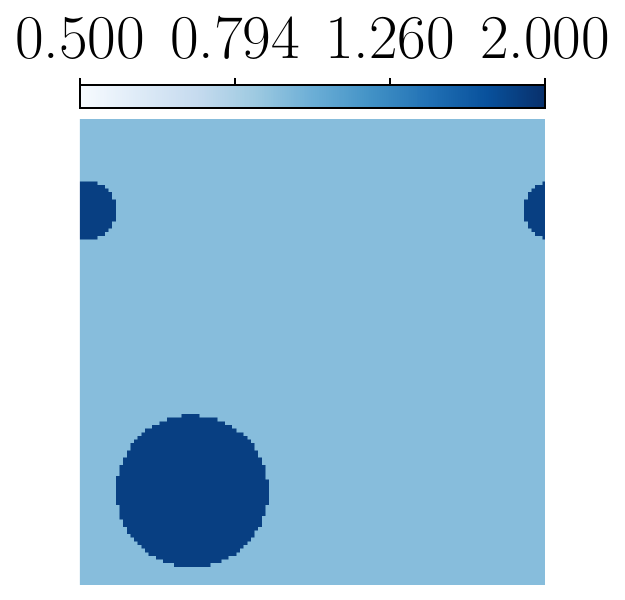}  & \includegraphics[width=\linewidth]{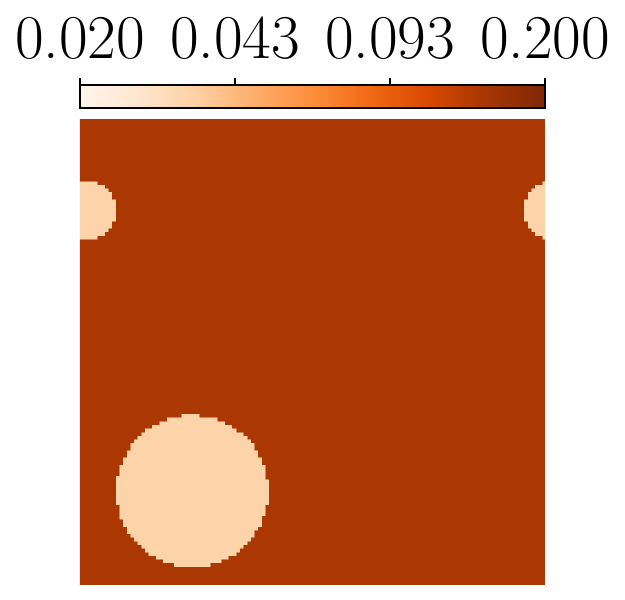} \\
       \raisebox{5pt}{\rotatebox[origin=c]{90}{\makecell{Poisson\\ratio}}}& \includegraphics[width=\linewidth]{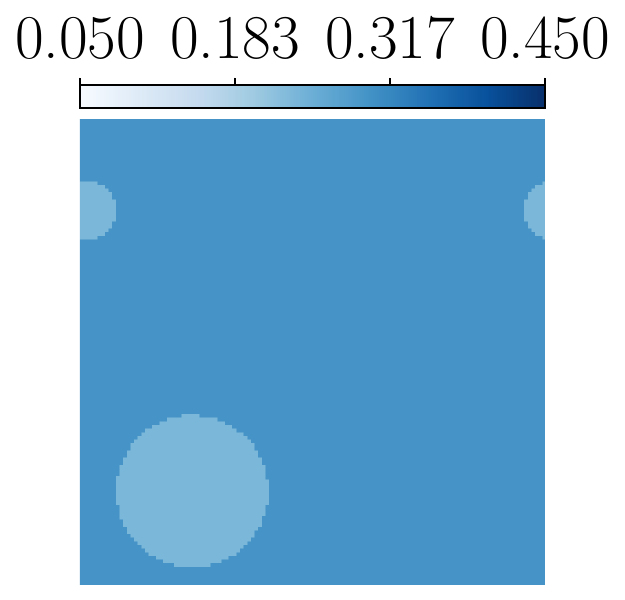} & \includegraphics[width=\linewidth]{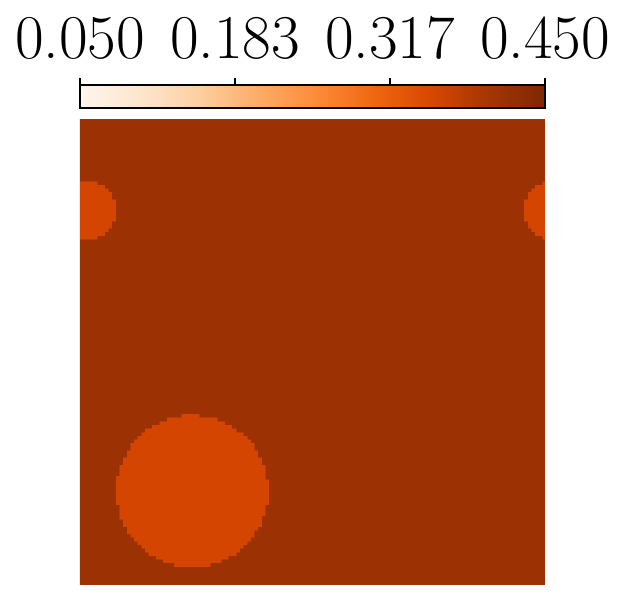}
    \end{tabular}
    }
    \begin{tabular}{c c}
        \includegraphics[width=0.245\linewidth, trim={0 0.9in 0 0},clip]{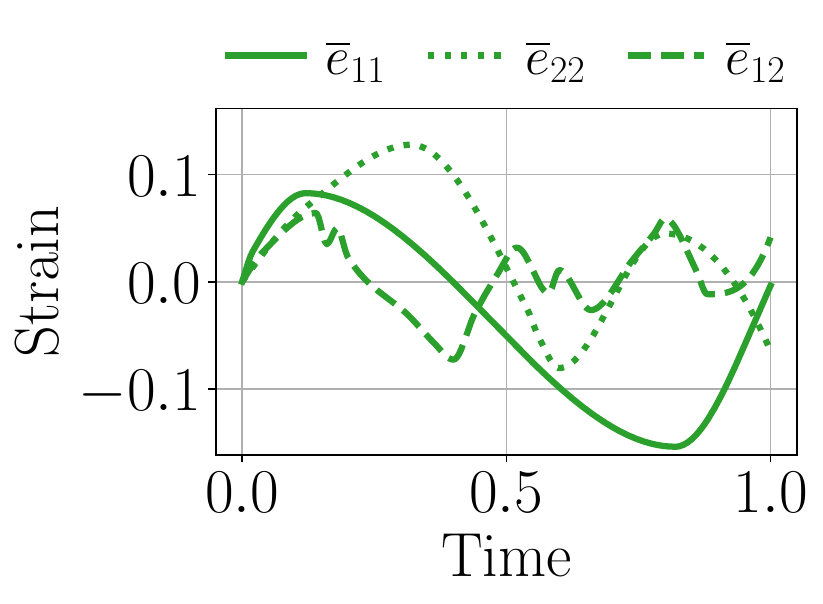} & \includegraphics[width=0.25\linewidth, trim={0 0.9in 0 0},clip]{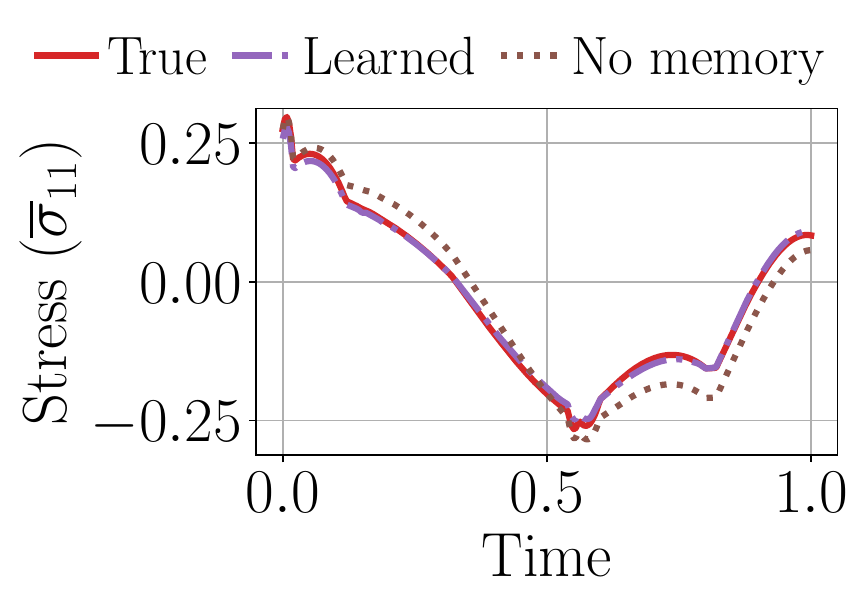}  \\
        \includegraphics[width=0.25\linewidth, trim={0 0 0 0.6in},clip]{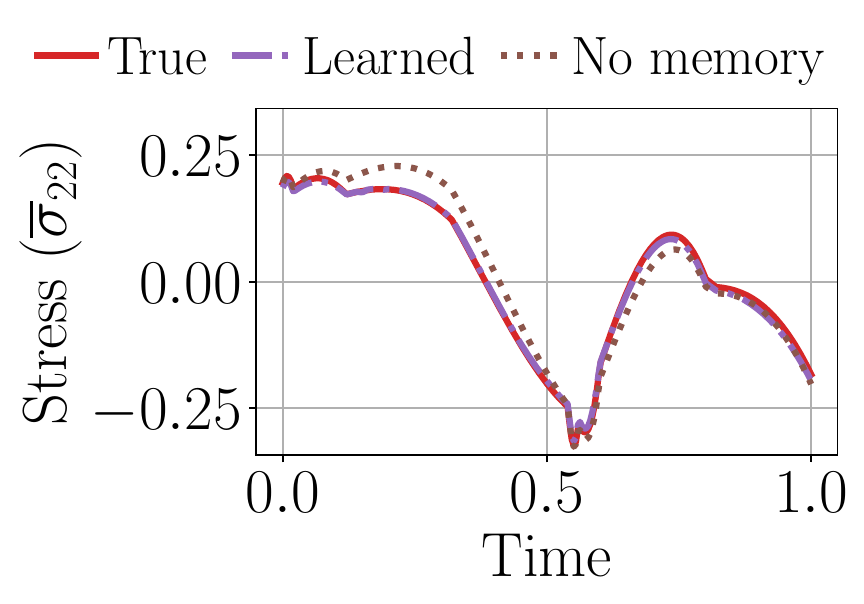} & \includegraphics[width=0.25\linewidth, trim={0 0 0 0.6in},clip]{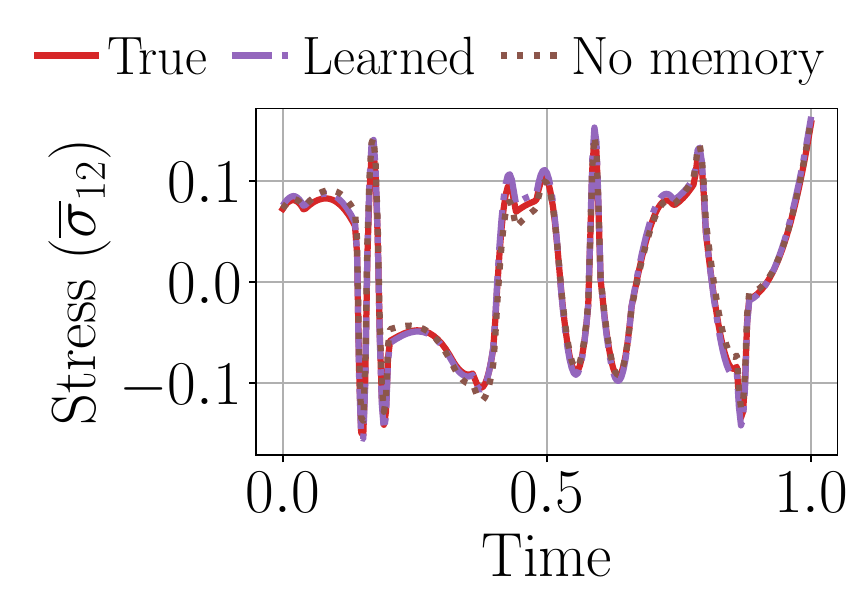}\\
        \multicolumn{2}{c}{\includegraphics[width=0.27\linewidth]{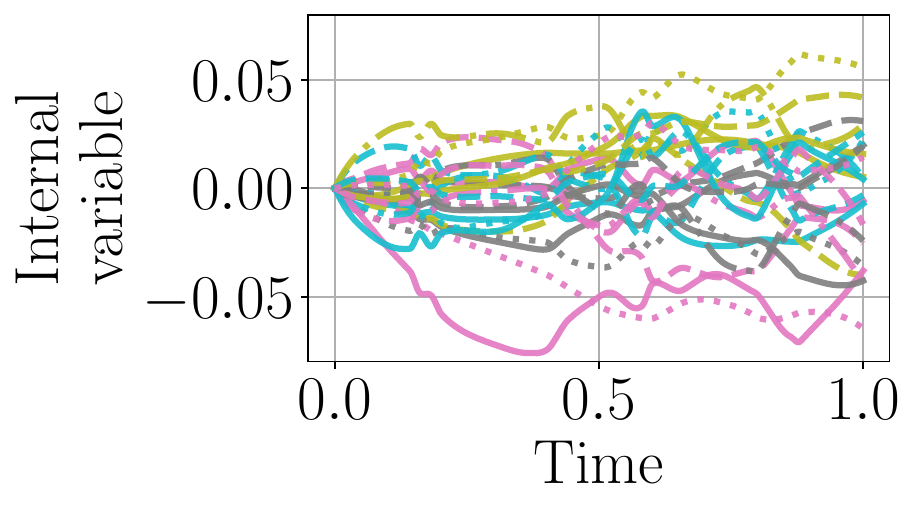}}
    \end{tabular}\\
\addtolength{\tabcolsep}{5 pt}
    \caption{Prediction by the trained neural operator model for a particulate viscoelasticity sample with median relative stress error, using internal-variable dimension $n_{\xi}=25$. The left panels show the microstructural property fields, while the right panels show the prescribed strain trajectory, the predicted and reference stress trajectory, and the predicted internal-variable trajectory. The memoryless prediction ($\bK=0$) is included for comparison.}
    \label{fig:kv2d_particulate_median}
\end{figure}

\Cref{fig:kv2d_prediction_median,fig:kv2d_particulate_median} show test samples whose relative stress errors are closest to the median for the baseline Voronoi and particulate datasets, respectively, using the trained neural operator model with $n_{\xi}=25$. These examples display the four scalar material fields, the imposed strain path, the predicted and reference stress components, and the predicted internal-variable trajectory. The memoryless response indicates the part of the response not captured without internal variables.

\section{Implementation Details for Crystal Plasticity}\label{app:cp3d_details}

{This appendix specifies the local constitutive equations, data-generation procedure, and training setup for the crystal-plasticity experiments in Section~\ref{subsec:3d_crystal_plasticity}.}

{
\subsection{Constitutive Equations}\label{app:cp3d_constitutive}
The cell problem \cref{eq:cp_cell_problem} is closed by the following local constitutive equations. We use the multiplicative decomposition into elastic and plastic distortions, $F=F_eF_p$, with elastic Green--Lagrange strain and stress
\begin{align*}
    E_e&=\tfrac12(F_e^{\top}F_e-I),
    &S_e&=\bC_{\rm crys}E_e,
    &P&=F_eS_eF_p^{-\top}.
\end{align*}
Here $\bC_{\rm crys}$ is the cubic single-crystal stiffness and $S_e$ is the elastic second Piola--Kirchhoff stress. The stiffness and slip systems are expressed in crystal coordinates. Plastic flow on the twelve FCC $\{111\}\langle110\rangle$ systems satisfies
\begin{align*}
    \dot F_p=L_pF_p,\qquad L_p=\sum_{\alpha=1}^{12}\dot\gamma^\alpha s^\alpha\otimes m^\alpha,\qquad \dot\gamma^\alpha=\dot\gamma_0\left|\frac{\tau^\alpha}{r^\alpha}\right|^n\operatorname{sgn}(\tau^\alpha),\qquad \tau^\alpha=S_e(s^\alpha\otimes m^\alpha),
\end{align*}
where $\dot\gamma^\alpha$ is the shear rate on system $\alpha$, with unit slip direction $s^\alpha$, unit plane normal $m^\alpha$, and slip resistance $r^\alpha$. The resolved shear stress follows DAMASK's small-elastic-strain approximation $F_e^{\top}F_e\simeq I$ in the plastic driving stress \cite[Appendix C]{roters2019damask}. The slip resistances evolve by self and latent hardening according to
\begin{equation*}
    \dot r^\alpha=h_0\operatorname{sgn}\!\left(1-\frac{r^\alpha}{r_\infty}\right)\left|1-\frac{r^\alpha}{r_\infty}\right|^a\sum_{\beta=1}^{12}Q_{\alpha\beta}|\dot\gamma^\beta|.
\end{equation*}
Here $h_0(y)$ is the hardening scale, $r_\infty(y)$ is the saturation resistance, and $Q_{\alpha\beta}$ is the FCC slip-interaction matrix. The initial orientation $O_0(y)$ maps specimen coordinates to crystal coordinates and is incorporated into the initial plastic distortion. With initial slip resistance $r_0(y)$, the initially stress-free state is specified by
\begin{equation*}
    \chi(y,0)=0,\qquad F_p(y,0)=O_0(y),\qquad r^\alpha(y,0)=r_0(y).
\end{equation*}
Since $\overline F(0)=I$, these conditions give $F_e(y,0)=O_0(y)^{\top}$ and $S_e(y,0)=0$.

The elastic constants are $(C_{11},C_{12},C_{44})=(168.9,121.8,75.8)\,\mathrm{GPa}$, the reference shear rate is $\dot\gamma_0=3\times10^{-3}\,\mathrm{s}^{-1}$, the rate exponent is $n=20$, and the saturation exponent is $a=0.6$. The seven coefficients defining $Q$ are $(1,1,5.123,0.574,1.123,1.123,1)$ for self, coplanar, collinear, Hirth, glissile I, glissile II, and Lomer interactions, respectively.} The grainwise fields $r_0$, $r_\infty$, and $h_0$ are sampled as follows.

\subsection{Microstructure Sampling}\label{app:cp3d_microstructure_sampling}
{Let $\{G_j\}_{j=1}^{N_g}$ denote the grains of the periodic Voronoi tessellation. The microstructural input is piecewise constant on grains, with $M(y)=(q_j,r_{0,j},r_{\infty,j},h_{0,j})\in\bR^7$ for $y\in G_j$. Here $q_j\in\bR^4$ is the orientation quaternion representing $O_0$ on $G_j$. For the single- and two-mode textures, the spread angle is sampled over the range $5^{\circ}$ to $35^{\circ}$.}

{
We sample the three slip-law parameters around the copper reference values using correlated grainwise scale factors. For each sample, we draw $\delta\sim\mathrm{Unif}[0.05,0.20]$, a resistance scale $\bar s_r$ log-uniformly from $[0.80,1.35]$, and a hardening scale $\bar s_h$ log-uniformly from $[0.80,1.50]$. For each grain $G_j$, we draw a standard Gaussian pair $(Z_{r,j},Z_{h,j})$ with correlation $0.65$ and set
\begin{equation*}
    s_{r,j}=\bar s_r\exp(\delta Z_{r,j}),\qquad s_{h,j}=\bar s_h\exp(\delta Z_{h,j}).
\end{equation*}
We set
\begin{equation*}
    r_{0,j}=31\,\mathrm{MPa}\,s_{r,j},\qquad r_{\infty,j}=96.4\,\mathrm{MPa}\,s_{r,j},\qquad h_{0,j}=0.35\,\mathrm{GPa}\,s_{h,j}.
\end{equation*}
Thus $r_{0,j}$ and $r_{\infty,j}$ vary together with a constant ratio, while $h_{0,j}$ varies through a distinct but positively correlated factor.
}

\subsection{Trajectory Sampling}\label{app:cp3d_loading_sampling}
We parameterize each macroscopic loading path by separating the rotation, volumetric stretch, and deviatoric stretch components. Specifically, we write
\begin{equation*}
    \overline F(t)=R(t)\exp(E(t)),
    \quad
    R(t)=\exp(\operatorname{skew}(\omega(t))),
\end{equation*}
where $E(t)$ is the logarithmic strain, $\omega(t)\in\bR^3$ is the rotation vector, and {$\operatorname{skew}(\omega(t))$} is the associated skew symmetric matrix, i.e., $\operatorname{skew}(w)v=w\times v$ for $v,w\in\bR^3$. The logarithmic strain is decomposed as
\begin{equation*}
    E(t)=\frac{\eta_{\rm vol}(t)}{3}I+E_{\rm dev}(t),
    \quad
    \eta_{\rm vol}(t)=\log\det\overline F(t),
    \quad
    \operatorname{tr}E_{\rm dev}(t)=0.
\end{equation*}
After choosing a basis $\{B_a\}_{a=1}^{5}$ for trace-free symmetric $3\times3$ tensors, we have $E_{\rm dev}(t)=\sum_{a=1}^{5}e_a(t)B_a$.
Thus, each loading path is sampled through the nine scalar coordinates $\left(\eta_{\rm vol}(t),e_1(t),\dots,e_5(t), \omega_1(t),\omega_2(t),\omega_3(t)\right)$
independently according to the procedure in \cref{app:kv2d_strain_sampling}. The initial scalar values are set so that $\overline F(0)=I$. Once the scalar paths have been sampled and interpolated, the deformation gradient is reconstructed on the common output time grid. This reconstruction enforces $\det\overline F(t)=\exp(\eta_{\rm vol}(t))>0$.

\subsection{Architecture and Training}\label{app:cp_training}
The architecture and training are similar to \cref{app:kv2d_training}, except that we use 8 Fourier modes, a width of 48, and a projection width of 96. We use a penalty weight of $\alpha=0.001$ across all internal-variable dimensions. We train on trajectory data with a larger step size (1/64) rather than the data-generation step size (1/128), and only on [0, 0.875] rather than the full unit time interval. This reduction lets us make the most of our computing resources; we use 4 Nvidia Tesla P100 GPUs for training. The testing is performed on the full trajectory data.

\subsection{Median-Error Prediction Visualization}

\begin{figure}[!htp]
    \centering
    {\bf The median error test sample: Crystal plasticity dataset}\\[5pt]
    \scalebox{0.7}{
\addtolength{\tabcolsep}{-5 pt}
    \begin{tabular}{>{\centering\arraybackslash} m{0.18\linewidth} >{\centering\arraybackslash} m{0.06\linewidth}
    >{\centering\arraybackslash} m{0.16\linewidth} 
    >{\centering\arraybackslash} m{0.16\linewidth}
    >{\centering\arraybackslash} m{0.16\linewidth}}
        \raisebox{-40 pt}{\multirow{4}{*}{\includegraphics[width=\linewidth]{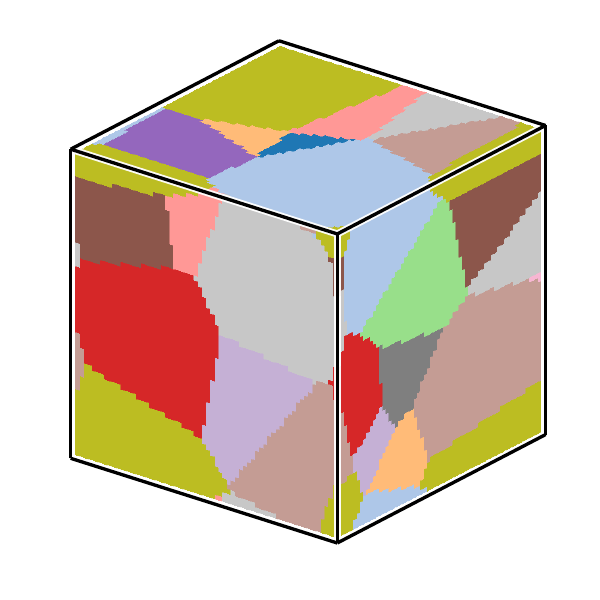}}} & \raisebox{45 pt}{\multirow{2}{*}{\rotatebox[origin=c]{90}{\makecell{Slip-law parameters \\ (MPa)}}}} &
        \includegraphics[width=\linewidth]{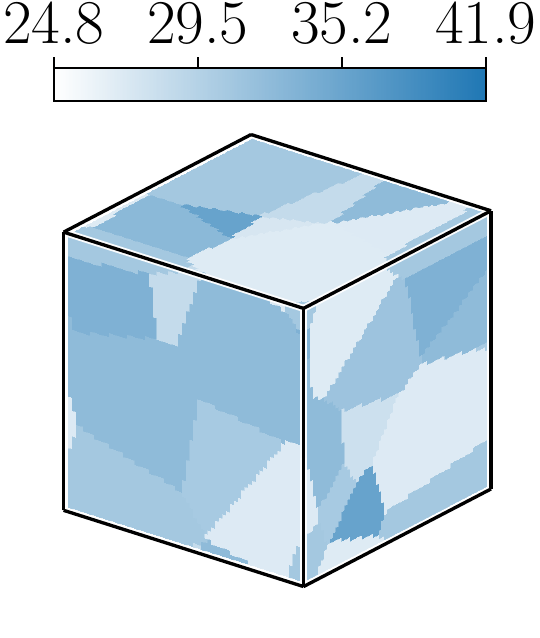}  & \includegraphics[width=\linewidth]{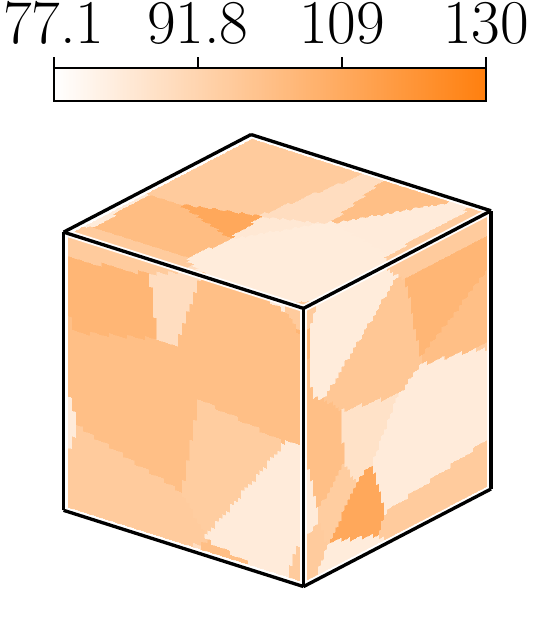} & \includegraphics[width=\linewidth]{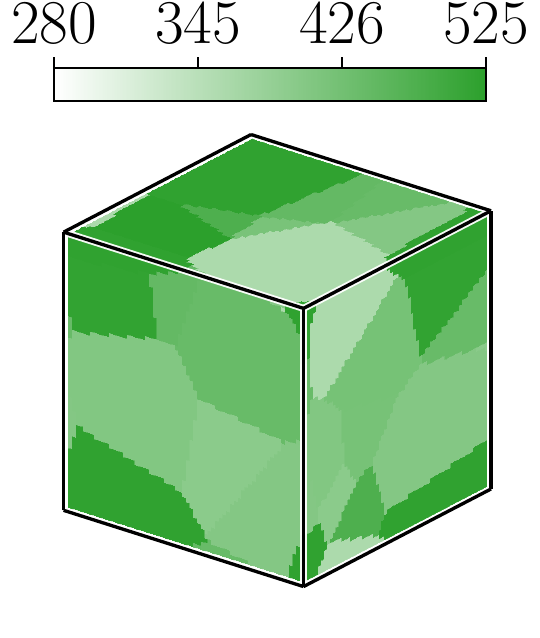} \\
         & & Initial slip resistance &Saturation slip resistance & Hardening scale \\
            & \raisebox{40 pt}{\rotatebox[origin=b]{90}{Pole figures}} & \multicolumn{3}{c}{\parbox[c]{0.48\linewidth}{\centering\includegraphics[width=0.875\linewidth]{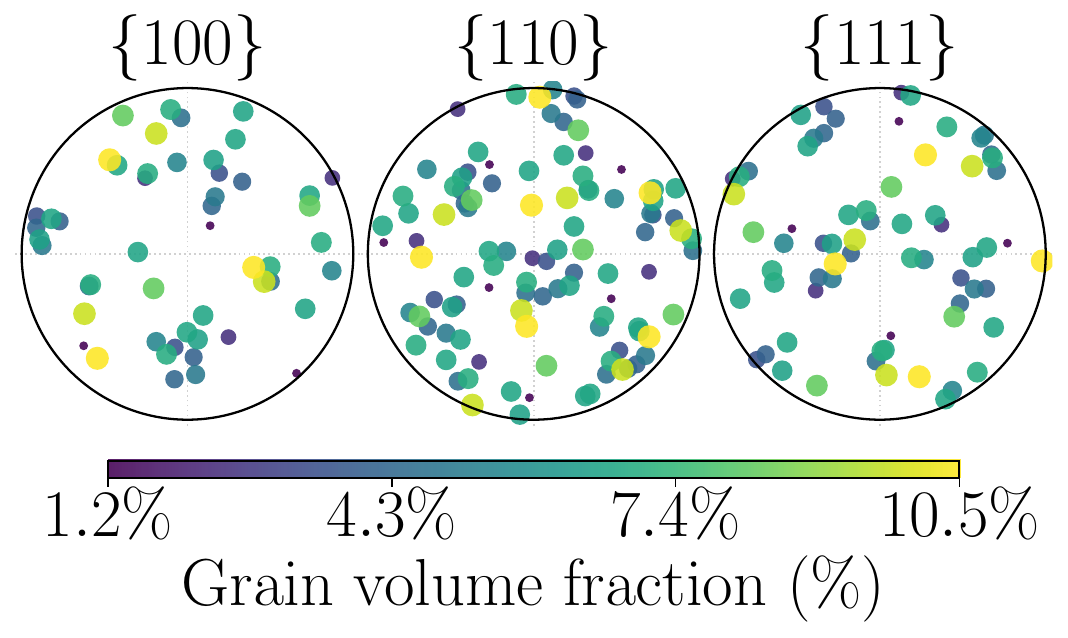}}}  \\
    \end{tabular}
    \begin{tabular}{c}
         \includegraphics[width=0.4\linewidth, trim={0 1.0in 0 0},clip]{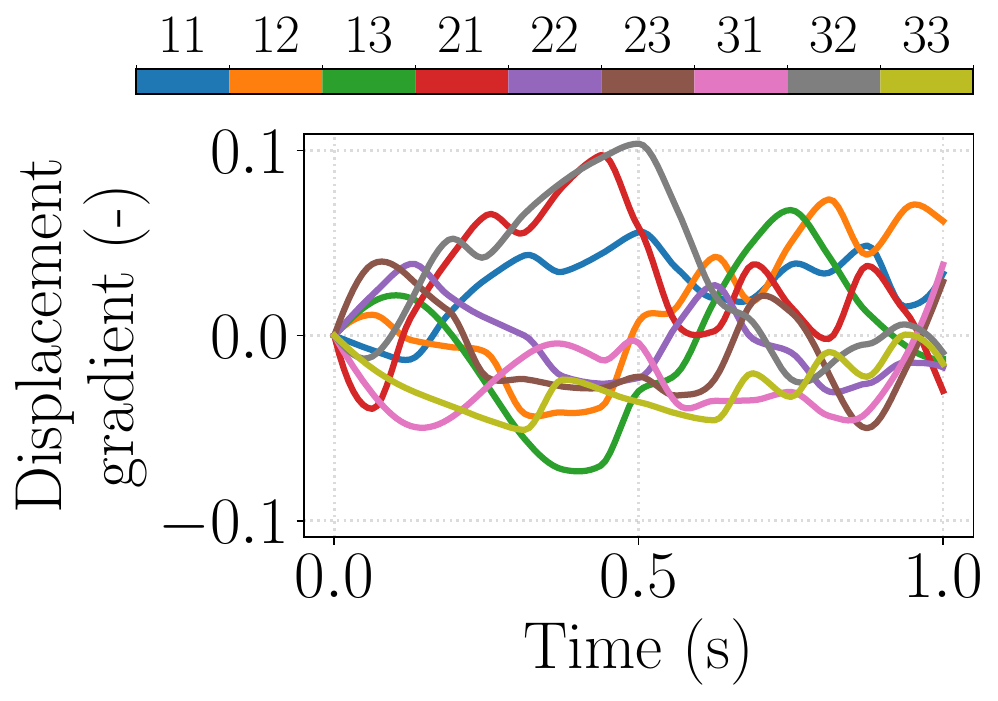}  \\
        \hspace{0.012\linewidth}\includegraphics[width=0.39\linewidth, trim={0 0 0 0.8in},clip]{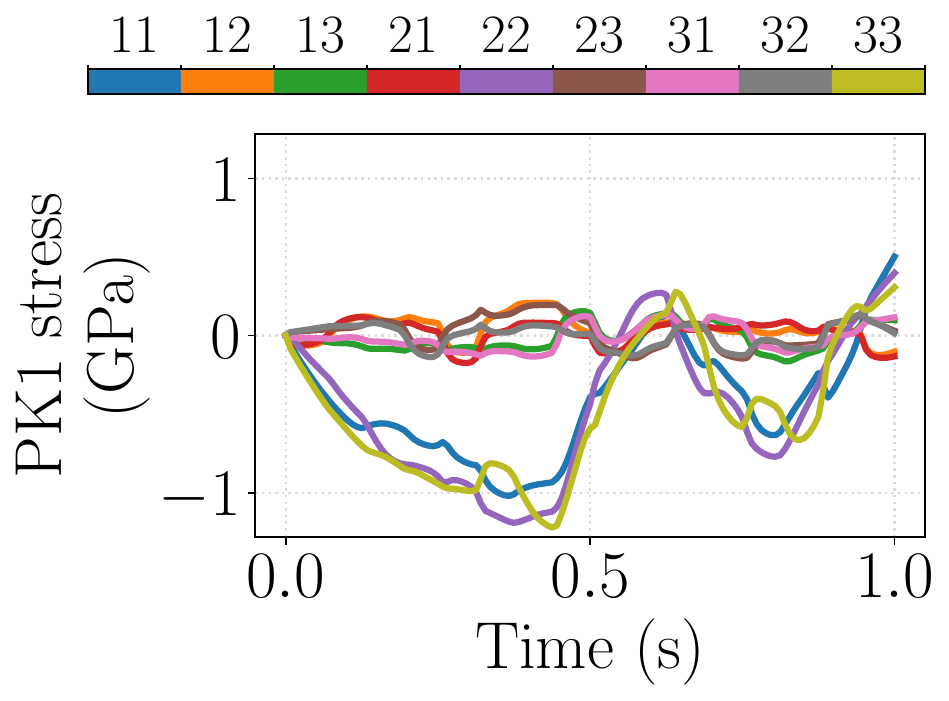}
    \end{tabular}
\addtolength{\tabcolsep}{5 pt}
    }
\scalebox{0.5}{
\addtolength{\tabcolsep}{-5 pt}
\centering
        
        \begin{tabular}{c c c c}
        & 11 & 12 & 13 \\
        \multirow{3}{*}{\rotatebox[origin=c]{90}{PK1 Stress (GPa)}}&\includegraphics[width=0.25\linewidth]{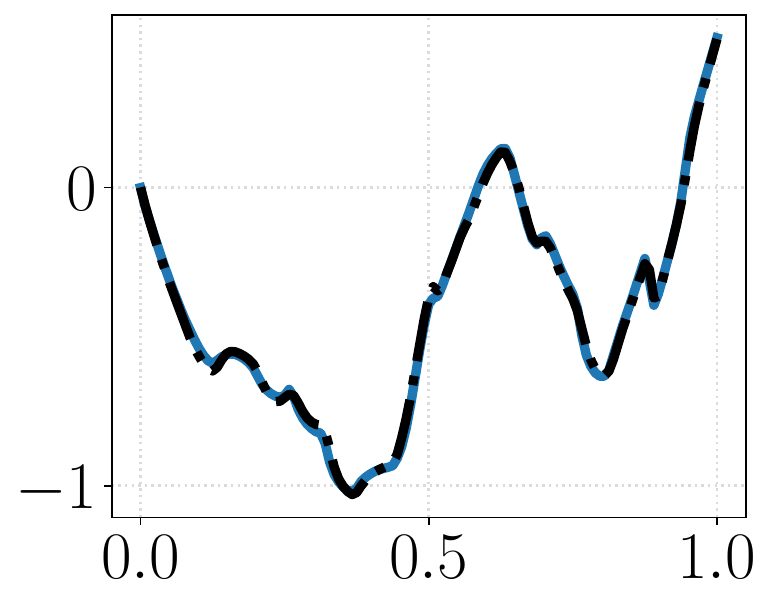} &\includegraphics[width=0.25\linewidth]{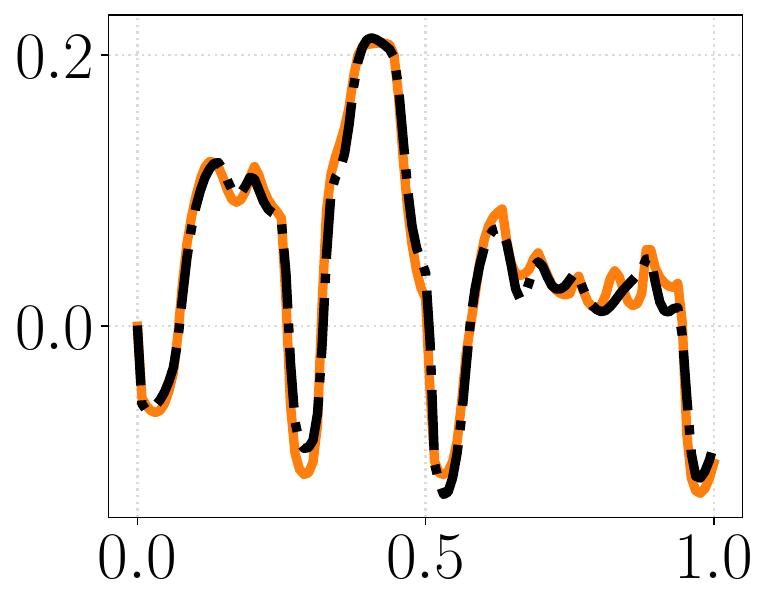}  
         & \includegraphics[width=0.265\linewidth]{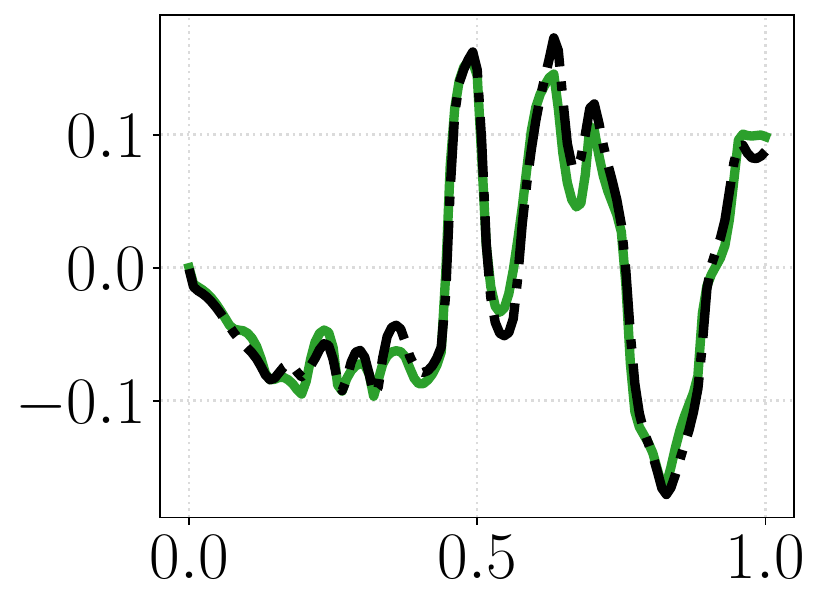} \\
        & 21 & 22 & 23 \\
         &\includegraphics[width=0.26\linewidth]{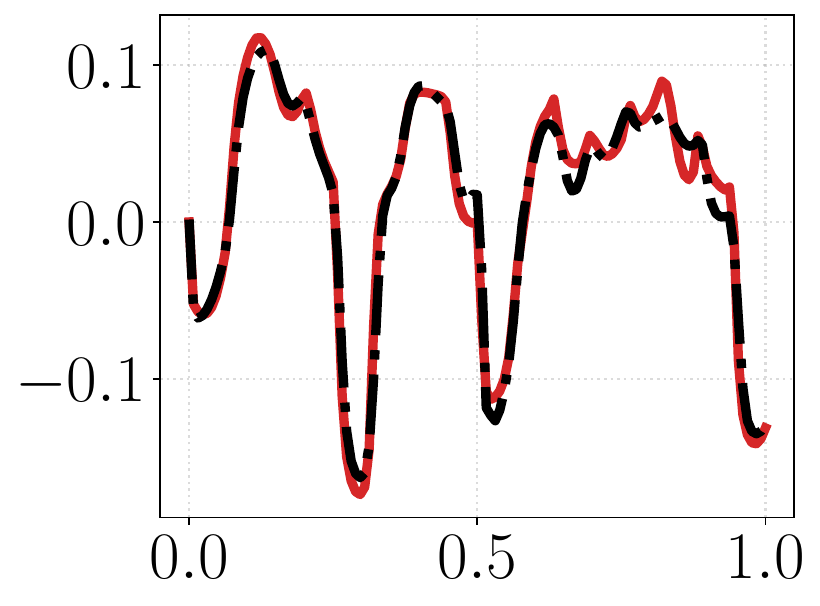} &\includegraphics[width=0.245\linewidth]{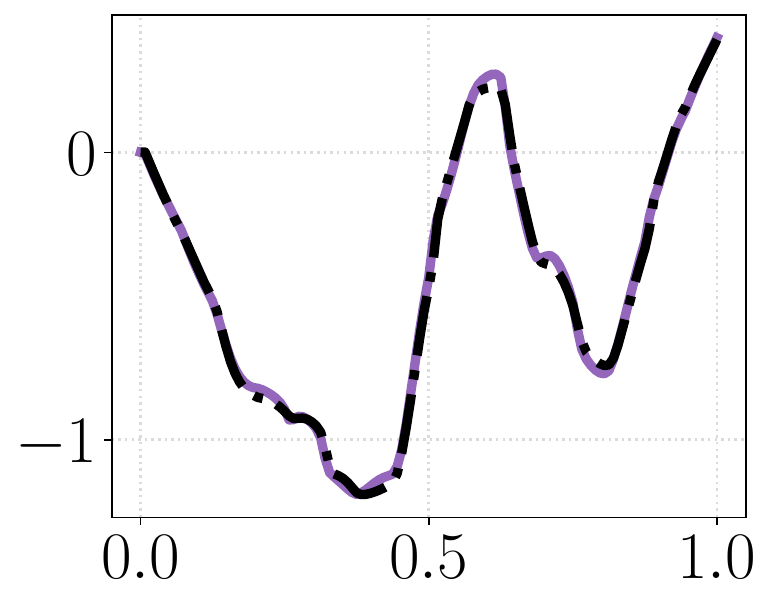}  
         & \includegraphics[width=0.245\linewidth]{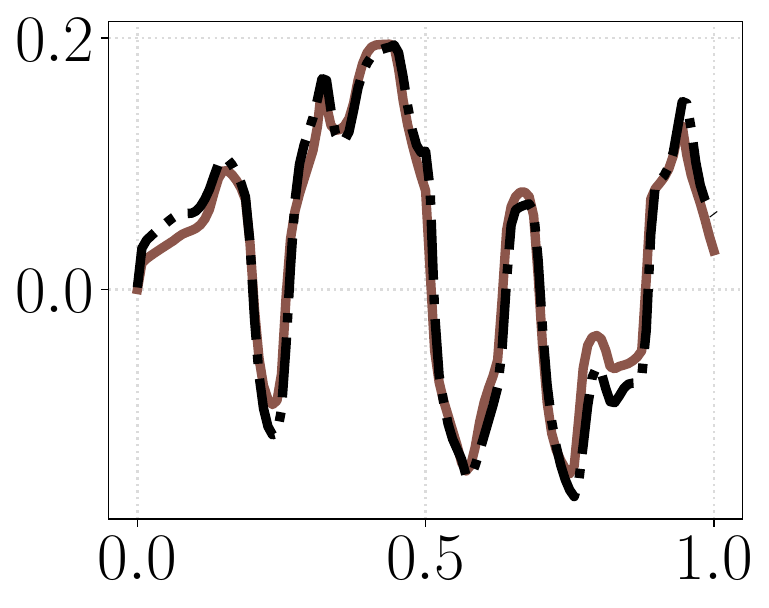} \\
        & 31 & 32 & 33 \\
        &\includegraphics[width=0.26\linewidth]{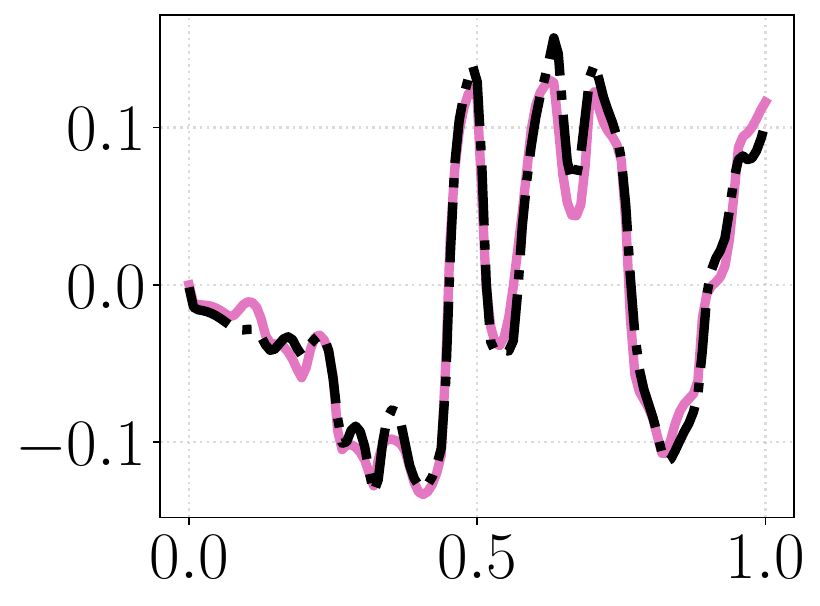} &\includegraphics[width=0.26\linewidth]{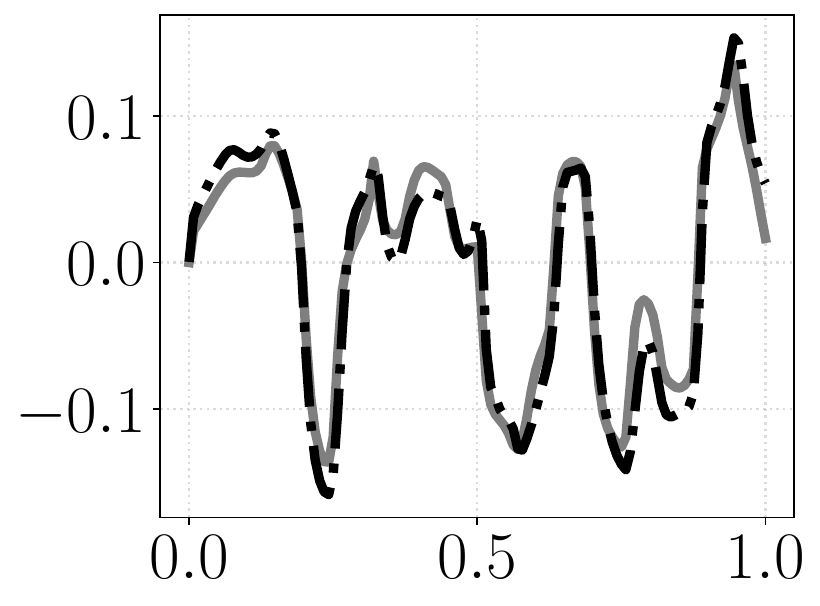}  
         & \includegraphics[width=0.245\linewidth]{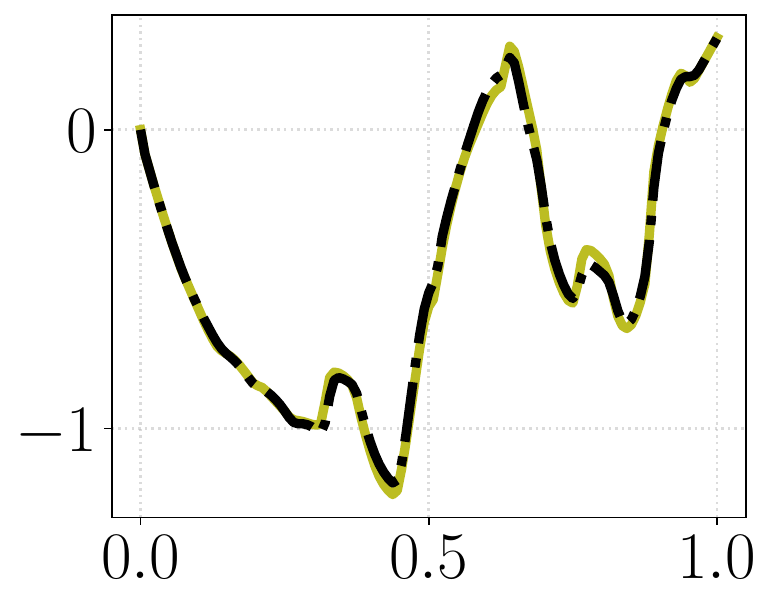}\\
         & \multicolumn{3}{c}{Time (s)}
    \end{tabular}
    \includegraphics[width=0.6\linewidth]{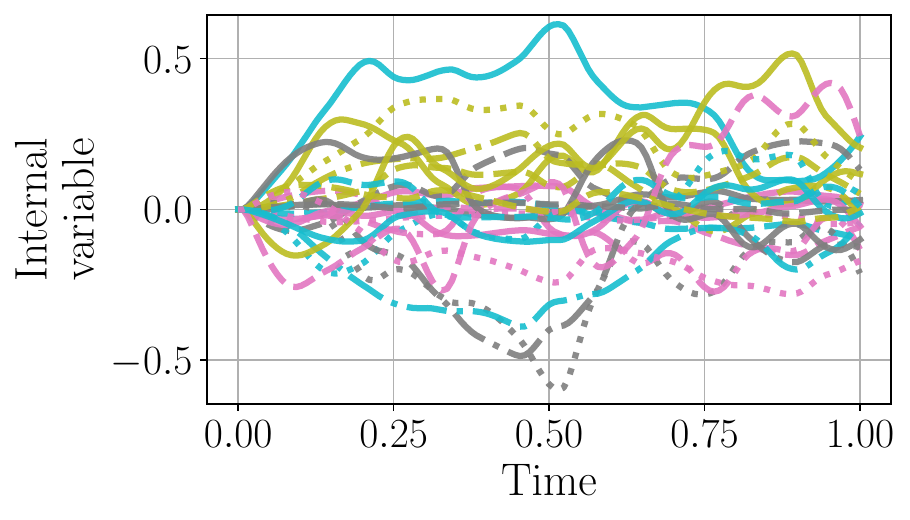}
    \addtolength{\tabcolsep}{5 pt}
    }
    
    \caption{(\emph{Top}) The test sample of the crystal plasticity dataset with a median relative stress prediction error (7.5\%) for the learned model. (\emph{Bottom left}) The stress-prediction components, where solid lines denote the reference homogenized response and dashed-dotted lines denote the stress response predicted by the learned model. (\emph{Bottom right}) The components of the internal variable predicted by the learned model.}
    \label{fig:cp_median_prediction}
\end{figure}

\Cref{fig:cp_median_prediction} shows the median-error test sample for the crystal-plasticity dataset, including the sampled microstructure, loading path, stress response, and learned internal variables.

\begin{thebibliography}{75}
\expandafter\ifx\csname natexlab\endcsname\relax\def\natexlab#1{#1}\fi
\providecommand{\url}[1]{\texttt{#1}}
\providecommand{\href}[2]{#2}
\providecommand{\path}[1]{#1}
\providecommand{\DOIprefix}{doi:}
\providecommand{\ArXivprefix}{arXiv:}
\providecommand{\URLprefix}{URL: }
\providecommand{\Pubmedprefix}{pmid:}
\providecommand{\doi}[1]{\href{http://dx.doi.org/#1}{\path{#1}}}
\providecommand{\Pubmed}[1]{\href{pmid:#1}{\path{#1}}}
\providecommand{\bibinfo}[2]{#2}
\ifx\xfnm\relax \def\xfnm[#1]{\unskip,\space#1}\fi
\bibitem[{Bensoussan et~al.(1978)Bensoussan, Lions, and
  Papanicolaou}]{bensoussan1978asymptotic}
\bibinfo{author}{A.~Bensoussan}, \bibinfo{author}{J.-L. Lions},
  \bibinfo{author}{G.~Papanicolaou}, \bibinfo{title}{Asymptotic Analysis for
  Periodic Structures}, volume~\bibinfo{volume}{5} of
  \textit{\bibinfo{series}{Studies in Mathematics and Its Applications}},
  \bibinfo{publisher}{North-Holland}, \bibinfo{year}{1978}.
  \DOIprefix\doi{10.1016/S0168-2024(08)70192-9}.
\bibitem[{S{\'a}nchez-Palencia(1980)}]{sanchezpalencia1980nonhomogeneous}
\bibinfo{author}{E.~S{\'a}nchez-Palencia}, \bibinfo{title}{Non-Homogeneous
  Media and Vibration Theory}, volume \bibinfo{volume}{127} of
  \textit{\bibinfo{series}{Lecture Notes in Physics}},
  \bibinfo{publisher}{Springer}, \bibinfo{year}{1980}.
  \DOIprefix\doi{10.1007/3-540-10000-8}.
\bibitem[{Allaire(1992)}]{allaire1992homogenization}
\bibinfo{author}{G.~Allaire},
\newblock \bibinfo{title}{Homogenization and two-scale convergence},
\newblock \bibinfo{journal}{SIAM Journal on Mathematical Analysis}
  \bibinfo{volume}{23} (\bibinfo{year}{1992}) \bibinfo{pages}{1482--1518}.
  \DOIprefix\doi{10.1137/0523084}.
\bibitem[{Jikov et~al.(1994)Jikov, Kozlov, and
  Oleinik}]{jikov1994homogenization}
\bibinfo{author}{V.~V. Jikov}, \bibinfo{author}{S.~M. Kozlov},
  \bibinfo{author}{O.~A. Oleinik}, \bibinfo{title}{Homogenization of
  Differential Operators and Integral Functionals},
  \bibinfo{publisher}{Springer}, \bibinfo{year}{1994}.
  \DOIprefix\doi{10.1007/978-3-642-84659-5}.
\bibitem[{Pavliotis and Stuart(2008)}]{pavliotis2008multiscale}
\bibinfo{author}{G.~A. Pavliotis}, \bibinfo{author}{A.~M. Stuart},
  \bibinfo{title}{Multiscale Methods: Averaging and Homogenization},
  volume~\bibinfo{volume}{53} of \textit{\bibinfo{series}{Texts in Applied
  Mathematics}}, \bibinfo{publisher}{Springer}, \bibinfo{address}{New York,
  NY}, \bibinfo{year}{2008}. \DOIprefix\doi{10.1007/978-0-387-73829-1}.
\bibitem[{Feyel and Chaboche(2000)}]{feyel2000fe2}
\bibinfo{author}{F.~Feyel}, \bibinfo{author}{J.-L. Chaboche},
\newblock \bibinfo{title}{{FE}$^2$ multiscale approach for modelling the
  elastoviscoplastic behaviour of long fibre {SiC/Ti} composite materials},
\newblock \bibinfo{journal}{Computer Methods in Applied Mechanics and
  Engineering} \bibinfo{volume}{183} (\bibinfo{year}{2000})
  \bibinfo{pages}{309--330}. \DOIprefix\doi{10.1016/S0045-7825(99)00224-8}.
\bibitem[{Miehe et~al.(1999)Miehe, Schr{\"o}der, and
  Schotte}]{miehe1999computational}
\bibinfo{author}{C.~Miehe}, \bibinfo{author}{J.~Schr{\"o}der},
  \bibinfo{author}{J.~Schotte},
\newblock \bibinfo{title}{Computational homogenization analysis in finite
  plasticity: simulation of texture development in polycrystalline materials},
\newblock \bibinfo{journal}{Computer Methods in Applied Mechanics and
  Engineering} \bibinfo{volume}{171} (\bibinfo{year}{1999})
  \bibinfo{pages}{387--418}. \DOIprefix\doi{10.1016/S0045-7825(98)00218-7}.
\bibitem[{Geers et~al.(2010)Geers, Kouznetsova, and
  Brekelmans}]{geers2010multiscale}
\bibinfo{author}{M.~G.~D. Geers}, \bibinfo{author}{V.~G. Kouznetsova},
  \bibinfo{author}{W.~A.~M. Brekelmans},
\newblock \bibinfo{title}{Multi-scale computational homogenization: trends and
  challenges},
\newblock \bibinfo{journal}{Journal of Computational and Applied Mathematics}
  \bibinfo{volume}{234} (\bibinfo{year}{2010}) \bibinfo{pages}{2175--2182}.
  \DOIprefix\doi{10.1016/j.cam.2009.08.077}.
\bibitem[{Ghavamian and Simone(2019)}]{ghavamian2019accelerating}
\bibinfo{author}{F.~Ghavamian}, \bibinfo{author}{A.~Simone},
\newblock \bibinfo{title}{Accelerating multiscale finite element simulations of
  history-dependent materials using a recurrent neural network},
\newblock \bibinfo{journal}{Computer Methods in Applied Mechanics and
  Engineering} \bibinfo{volume}{357} (\bibinfo{year}{2019})
  \bibinfo{pages}{112594}. \DOIprefix\doi{10.1016/j.cma.2019.112594}.
\bibitem[{Logarzo et~al.(2021)Logarzo, Capuano, and Rimoli}]{logarzo2021smart}
\bibinfo{author}{H.~J. Logarzo}, \bibinfo{author}{G.~Capuano},
  \bibinfo{author}{J.~J. Rimoli},
\newblock \bibinfo{title}{Smart constitutive laws: Inelastic homogenization
  through machine learning},
\newblock \bibinfo{journal}{Computer Methods in Applied Mechanics and
  Engineering} \bibinfo{volume}{373} (\bibinfo{year}{2021})
  \bibinfo{pages}{113482}. \DOIprefix\doi{10.1016/j.cma.2020.113482}.
\bibitem[{Liu et~al.(2022)Liu, Kovachki, Li, Azizzadenesheli, Anandkumar,
  Stuart, and Bhattacharya}]{liu2022learning}
\bibinfo{author}{B.~Liu}, \bibinfo{author}{N.~Kovachki},
  \bibinfo{author}{Z.~Li}, \bibinfo{author}{K.~Azizzadenesheli},
  \bibinfo{author}{A.~Anandkumar}, \bibinfo{author}{A.~M. Stuart},
  \bibinfo{author}{K.~Bhattacharya},
\newblock \bibinfo{title}{A learning-based multiscale method and its
  application to inelastic impact problems},
\newblock \bibinfo{journal}{Journal of the Mechanics and Physics of Solids}
  \bibinfo{volume}{158} (\bibinfo{year}{2022}) \bibinfo{pages}{104668}.
  \DOIprefix\doi{10.1016/j.jmps.2021.104668}.
\bibitem[{Liu et~al.(2023)Liu, Ocegueda, Trautner, Stuart, and
  Bhattacharya}]{liu2023learning}
\bibinfo{author}{B.~Liu}, \bibinfo{author}{E.~Ocegueda},
  \bibinfo{author}{M.~Trautner}, \bibinfo{author}{A.~M. Stuart},
  \bibinfo{author}{K.~Bhattacharya},
\newblock \bibinfo{title}{Learning macroscopic internal variables and history
  dependence from microscopic models},
\newblock \bibinfo{journal}{Journal of the Mechanics and Physics of Solids}
  \bibinfo{volume}{178} (\bibinfo{year}{2023}) \bibinfo{pages}{105329}.
  \DOIprefix\doi{10.1016/j.jmps.2023.105329}.
\bibitem[{Bhattacharya et~al.(2023)Bhattacharya, Liu, Stuart, and
  Trautner}]{bhattacharya2023learning}
\bibinfo{author}{K.~Bhattacharya}, \bibinfo{author}{B.~Liu},
  \bibinfo{author}{A.~M. Stuart}, \bibinfo{author}{M.~Trautner},
\newblock \bibinfo{title}{Learning markovian homogenized models in
  viscoelasticity},
\newblock \bibinfo{journal}{Multiscale Modeling \& Simulation}
  \bibinfo{volume}{21} (\bibinfo{year}{2023}) \bibinfo{pages}{641--679}.
  \DOIprefix\doi{10.1137/22M1499200}.
\bibitem[{Karimi and Bhattacharya(2024)}]{karimi2024learning}
\bibinfo{author}{M.~Karimi}, \bibinfo{author}{K.~Bhattacharya},
\newblock \bibinfo{title}{A learning-based multiscale model for reactive flow
  in porous media},
\newblock \bibinfo{journal}{Water Resources Research} \bibinfo{volume}{60}
  (\bibinfo{year}{2024}) \bibinfo{pages}{e2023WR036303}.
  \DOIprefix\doi{10.1029/2023WR036303}.
\bibitem[{Raj et~al.(2026)Raj, Cao, Stuart, and Bhattacharya}]{raj2026neural}
\bibinfo{author}{M.~Raj}, \bibinfo{author}{L.~Cao},
  \bibinfo{author}{A.~Stuart}, \bibinfo{author}{K.~Bhattacharya},
\newblock \bibinfo{title}{A neural-network framework to learn history-dependent
  constitutive laws and identifiability of internal variables},
\newblock \bibinfo{journal}{Journal of the Mechanics and Physics of Solids}
  \bibinfo{volume}{217} (\bibinfo{year}{2026}) \bibinfo{pages}{106807}.
  \DOIprefix\doi{10.1016/j.jmps.2026.106807}.
\bibitem[{Mozaffar et~al.(2019)Mozaffar, Bostanabad, Chen, Ehmann, Cao, and
  Bessa}]{mozaffar2019deep}
\bibinfo{author}{M.~Mozaffar}, \bibinfo{author}{R.~Bostanabad},
  \bibinfo{author}{W.~Chen}, \bibinfo{author}{K.~Ehmann},
  \bibinfo{author}{J.~Cao}, \bibinfo{author}{M.~A. Bessa},
\newblock \bibinfo{title}{Deep learning predicts path-dependent plasticity},
\newblock \bibinfo{journal}{Proceedings of the National Academy of Sciences}
  \bibinfo{volume}{116} (\bibinfo{year}{2019}) \bibinfo{pages}{26414--26420}.
  \DOIprefix\doi{10.1073/pnas.1911815116}.
\bibitem[{Liu et~al.(2019)Liu, Wu, and Koishi}]{liu2019a}
\bibinfo{author}{Z.~Liu}, \bibinfo{author}{C.~T. Wu},
  \bibinfo{author}{M.~Koishi},
\newblock \bibinfo{title}{A deep material network for multiscale topology
  learning and accelerated nonlinear modeling of heterogeneous materials},
\newblock \bibinfo{journal}{Computer Methods in Applied Mechanics and
  Engineering} \bibinfo{volume}{345} (\bibinfo{year}{2019})
  \bibinfo{pages}{1138--1168}. \DOIprefix\doi{10.1016/j.cma.2018.09.020}.
\bibitem[{Gorji et~al.(2020)Gorji, Mozaffar, Heidenreich, Cao, and
  Mohr}]{gorji2020on}
\bibinfo{author}{M.~B. Gorji}, \bibinfo{author}{M.~Mozaffar},
  \bibinfo{author}{J.~N. Heidenreich}, \bibinfo{author}{J.~Cao},
  \bibinfo{author}{D.~Mohr},
\newblock \bibinfo{title}{On the potential of recurrent neural networks for
  modeling path dependent plasticity},
\newblock \bibinfo{journal}{Journal of the Mechanics and Physics of Solids}
  \bibinfo{volume}{143} (\bibinfo{year}{2020}) \bibinfo{pages}{103972}.
  \DOIprefix\doi{10.1016/j.jmps.2020.103972}.
\bibitem[{Rao and Liu(2020)}]{rao2020threedimensional}
\bibinfo{author}{C.~Rao}, \bibinfo{author}{Y.~Liu},
\newblock \bibinfo{title}{Three-dimensional convolutional neural network
  ({3D-CNN}) for heterogeneous material homogenization},
\newblock \bibinfo{journal}{Computational Materials Science}
  \bibinfo{volume}{184} (\bibinfo{year}{2020}) \bibinfo{pages}{109850}.
  \DOIprefix\doi{10.1016/j.commatsci.2020.109850}.
\bibitem[{Kovachki et~al.(2023)Kovachki, Li, Liu, Azizzadenesheli,
  Bhattacharya, Stuart, and Anandkumar}]{kovachki2023neural}
\bibinfo{author}{N.~Kovachki}, \bibinfo{author}{Z.~Li},
  \bibinfo{author}{B.~Liu}, \bibinfo{author}{K.~Azizzadenesheli},
  \bibinfo{author}{K.~Bhattacharya}, \bibinfo{author}{A.~M. Stuart},
  \bibinfo{author}{A.~Anandkumar},
\newblock \bibinfo{title}{Neural operator: Learning maps between function
  spaces with applications to {PDEs}},
\newblock \bibinfo{journal}{Journal of Machine Learning Research}
  \bibinfo{volume}{24} (\bibinfo{year}{2023}) \bibinfo{pages}{1--97}.
\bibitem[{Bhattacharya et~al.(2024)Bhattacharya, Kovachki, Rajan, Stuart, and
  Trautner}]{bhattacharya2024learning}
\bibinfo{author}{K.~Bhattacharya}, \bibinfo{author}{N.~B. Kovachki},
  \bibinfo{author}{A.~Rajan}, \bibinfo{author}{A.~M. Stuart},
  \bibinfo{author}{M.~Trautner},
\newblock \bibinfo{title}{Learning homogenization for elliptic operators},
\newblock \bibinfo{journal}{SIAM Journal on Numerical Analysis}
  \bibinfo{volume}{62} (\bibinfo{year}{2024}) \bibinfo{pages}{1844--1873}.
  \DOIprefix\doi{10.1137/23M1585015}.
\bibitem[{Bhattacharya(2026)}]{bhattacharya2026multiscale}
\bibinfo{author}{K.~Bhattacharya},
\newblock \bibinfo{title}{Multiscale modeling of materials and neural
  operators},
\newblock \bibinfo{journal}{MRS Bulletin} \bibinfo{volume}{51}
  (\bibinfo{year}{2026}) \bibinfo{pages}{913--923}.
  \DOIprefix\doi{10.1557/s43577-026-01117-8}.
\bibitem[{Bhattacharya et~al.(2026)Bhattacharya, Cao, Stepaniants, Stuart, and
  Trautner}]{bhattacharya2025learning}
\bibinfo{author}{K.~Bhattacharya}, \bibinfo{author}{L.~Cao},
  \bibinfo{author}{G.~Stepaniants}, \bibinfo{author}{A.~M. Stuart},
  \bibinfo{author}{M.~Trautner},
\newblock \bibinfo{title}{Learning memory and material dependent constitutive
  laws},
\newblock \bibinfo{journal}{The SMAI Journal of Computational Mathematics}
  \bibinfo{volume}{12} (\bibinfo{year}{2026}) \bibinfo{pages}{219--267}.
  \DOIprefix\doi{10.5802/smai-jcm.148}.
\bibitem[{Coleman and Gurtin(1967)}]{coleman1967thermodynamics}
\bibinfo{author}{B.~D. Coleman}, \bibinfo{author}{M.~E. Gurtin},
\newblock \bibinfo{title}{Thermodynamics with internal state variables},
\newblock \bibinfo{journal}{The Journal of Chemical Physics}
  \bibinfo{volume}{47} (\bibinfo{year}{1967}) \bibinfo{pages}{597--613}.
  \DOIprefix\doi{10.1063/1.1711937}.
\bibitem[{Rice(1971)}]{rice1971inelastic}
\bibinfo{author}{J.~Rice},
\newblock \bibinfo{title}{Inelastic constitutive relations for solids: An
  internal-variable theory and its application to metal plasticity},
\newblock \bibinfo{journal}{Journal of the Mechanics and Physics of Solids}
  \bibinfo{volume}{19} (\bibinfo{year}{1971}) \bibinfo{pages}{433--455}.
  \DOIprefix\doi{10.1016/0022-5096(71)90010-X}.
\bibitem[{Halphen and Son(1975)}]{halphen1975sur}
\bibinfo{author}{B.~Halphen}, \bibinfo{author}{N.~Q. Son},
\newblock \bibinfo{title}{Sur les mat{\'e}riaux standard
  g{\'e}n{\'e}ralis{\'e}s},
\newblock \bibinfo{journal}{Journal de M{\'e}canique} \bibinfo{volume}{14}
  (\bibinfo{year}{1975}) \bibinfo{pages}{39--63}.
\bibitem[{He and Chen(2022)}]{he2022thermodynamically}
\bibinfo{author}{X.~He}, \bibinfo{author}{J.-S. Chen},
\newblock \bibinfo{title}{Thermodynamically consistent machine-learned internal
  state variable approach for data-driven modeling of path-dependent
  materials},
\newblock \bibinfo{journal}{Computer Methods in Applied Mechanics and
  Engineering} \bibinfo{volume}{402} (\bibinfo{year}{2022})
  \bibinfo{pages}{115348}. \DOIprefix\doi{10.1016/j.cma.2022.115348}.
\bibitem[{Jones et~al.(2022)Jones, Frankel, and Johnson}]{jones2022neural}
\bibinfo{author}{R.~E. Jones}, \bibinfo{author}{A.~L. Frankel},
  \bibinfo{author}{K.~L. Johnson},
\newblock \bibinfo{title}{A neural ordinary differential equation framework for
  modeling inelastic stress response via internal state variables},
\newblock \bibinfo{journal}{Journal of Machine Learning for Modeling and
  Computing} \bibinfo{volume}{3} (\bibinfo{year}{2022}) \bibinfo{pages}{1--35}.
  \DOIprefix\doi{10.1615/JMachLearnModelComput.2022042917}.
\bibitem[{Rosenkranz et~al.(2024)Rosenkranz, Kalina, Brummund, Sun, and
  K{\"a}stner}]{rosenkranz2024viscoelasticity}
\bibinfo{author}{M.~Rosenkranz}, \bibinfo{author}{K.~A. Kalina},
  \bibinfo{author}{J.~Brummund}, \bibinfo{author}{W.~Sun},
  \bibinfo{author}{M.~K{\"a}stner},
\newblock \bibinfo{title}{Viscoelasticty with physics-augmented neural
  networks: model formulation and training methods without prescribed internal
  variables},
\newblock \bibinfo{journal}{Computational Mechanics} \bibinfo{volume}{74}
  (\bibinfo{year}{2024}) \bibinfo{pages}{1279--1301}.
  \DOIprefix\doi{10.1007/s00466-024-02477-1}.
\bibitem[{Flaschel et~al.(2025)Flaschel, Steinmann, {De Lorenzis}, and
  Kuhl}]{flaschel2025convex}
\bibinfo{author}{M.~Flaschel}, \bibinfo{author}{P.~Steinmann},
  \bibinfo{author}{L.~{De Lorenzis}}, \bibinfo{author}{E.~Kuhl},
\newblock \bibinfo{title}{Convex neural networks learn generalized standard
  material models},
\newblock \bibinfo{journal}{Journal of the Mechanics and Physics of Solids}
  \bibinfo{volume}{200} (\bibinfo{year}{2025}) \bibinfo{pages}{106103}.
  \DOIprefix\doi{10.1016/j.jmps.2025.106103}.
\bibitem[{Jones and Fuhg(2025)}]{jones2025attention}
\bibinfo{author}{R.~E. Jones}, \bibinfo{author}{J.~N. Fuhg},
\newblock \bibinfo{title}{An attention-based neural ordinary differential
  equation framework for modeling inelastic processes},
\newblock \bibinfo{journal}{Computer Methods in Applied Mechanics and
  Engineering} \bibinfo{volume}{445} (\bibinfo{year}{2025})
  \bibinfo{pages}{118186}. \DOIprefix\doi{10.1016/j.cma.2025.118186}.
\bibitem[{Qiu et~al.(2025)Qiu, Huang, and Reina}]{qiu2025bridging}
\bibinfo{author}{W.~Qiu}, \bibinfo{author}{S.~Huang},
  \bibinfo{author}{C.~Reina},
\newblock \bibinfo{title}{Bridging statistical mechanics and thermodynamics
  away from equilibrium: A data-driven approach for learning internal variables
  and their dynamics},
\newblock \bibinfo{journal}{Journal of the Mechanics and Physics of Solids}
  \bibinfo{volume}{203} (\bibinfo{year}{2025}) \bibinfo{pages}{106211}.
  \DOIprefix\doi{10.1016/j.jmps.2025.106211}.
\bibitem[{Lu et~al.(2021)Lu, Jin, Pang, Zhang, and
  Karniadakis}]{lu2021learning}
\bibinfo{author}{L.~Lu}, \bibinfo{author}{P.~Jin}, \bibinfo{author}{G.~Pang},
  \bibinfo{author}{Z.~Zhang}, \bibinfo{author}{G.~E. Karniadakis},
\newblock \bibinfo{title}{Learning nonlinear operators via {DeepONet} based on
  the universal approximation theorem of operators},
\newblock \bibinfo{journal}{Nature Machine Intelligence} \bibinfo{volume}{3}
  (\bibinfo{year}{2021}) \bibinfo{pages}{218--229}.
  \DOIprefix\doi{10.1038/s42256-021-00302-5}.
\bibitem[{Lanthaler et~al.(2025)Lanthaler, Li, and
  Stuart}]{lanthaler2025nonlocality}
\bibinfo{author}{S.~Lanthaler}, \bibinfo{author}{Z.~Li}, \bibinfo{author}{A.~M.
  Stuart},
\newblock \bibinfo{title}{Nonlocality and nonlinearity implies universality in
  operator learning},
\newblock \bibinfo{journal}{Constructive Approximation} \bibinfo{volume}{62}
  (\bibinfo{year}{2025}) \bibinfo{pages}{261--303}.
  \DOIprefix\doi{10.1007/s00365-025-09718-3}.
\bibitem[{Coleman and Noll(1960)}]{coleman1960an}
\bibinfo{author}{B.~D. Coleman}, \bibinfo{author}{W.~Noll},
\newblock \bibinfo{title}{An approximation theorem for functionals, with
  applications in continuum mechanics},
\newblock \bibinfo{journal}{Archive for Rational Mechanics and Analysis}
  \bibinfo{volume}{6} (\bibinfo{year}{1960}) \bibinfo{pages}{355--370}.
  \DOIprefix\doi{10.1007/BF00276168}.
\bibitem[{Coleman and Noll(1961)}]{coleman1961foundations}
\bibinfo{author}{B.~D. Coleman}, \bibinfo{author}{W.~Noll},
\newblock \bibinfo{title}{Foundations of linear viscoelasticity},
\newblock \bibinfo{journal}{Reviews of Modern Physics} \bibinfo{volume}{33}
  (\bibinfo{year}{1961}) \bibinfo{pages}{239--249}.
  \DOIprefix\doi{10.1103/RevModPhys.33.239}.
\bibitem[{Tartar(1990)}]{tartar1990memory}
\bibinfo{author}{L.~Tartar},
\newblock \bibinfo{title}{Memory effects and homogenization},
\newblock \bibinfo{journal}{Archive for Rational Mechanics and Analysis}
  \bibinfo{volume}{111} (\bibinfo{year}{1990}) \bibinfo{pages}{121--133}.
  \DOIprefix\doi{10.1007/BF00375404}.
\bibitem[{Ortiz(2026)}]{ortiz2026linear}
\bibinfo{author}{M.~Ortiz},
\newblock \bibinfo{title}{Linear viscoelasticity: Mechanics, analysis and
  approximation},
\newblock \bibinfo{journal}{Archives of Computational Methods in Engineering}
  \bibinfo{volume}{33} (\bibinfo{year}{2026}) \bibinfo{pages}{4133--4153}.
  \DOIprefix\doi{10.1007/s11831-025-10426-1}.
\bibitem[{Francfort and Suquet(1986)}]{francfort1986homogenization}
\bibinfo{author}{G.~A. Francfort}, \bibinfo{author}{P.~M. Suquet},
\newblock \bibinfo{title}{Homogenization and mechanical dissipation in
  thermoviscoelasticity},
\newblock \bibinfo{journal}{Archive for Rational Mechanics and Analysis}
  \bibinfo{volume}{96} (\bibinfo{year}{1986}) \bibinfo{pages}{265--293}.
\bibitem[{de~Moura and Ullrich(2021)}]{demoura2021lean4}
\bibinfo{author}{L.~de~Moura}, \bibinfo{author}{S.~Ullrich},
\newblock \bibinfo{title}{The {Lean} 4 theorem prover and programming
  language},
\newblock in: \bibinfo{booktitle}{Automated Deduction---CADE 28}, volume
  \bibinfo{volume}{12699} of \textit{\bibinfo{series}{Lecture Notes in Computer
  Science}}, \bibinfo{publisher}{Springer}, \bibinfo{year}{2021}, pp.
  \bibinfo{pages}{625--635}. \DOIprefix\doi{10.1007/978-3-030-79876-5_37}.
\bibitem[{{The mathlib Community}(2020)}]{mathlib2020}
\bibinfo{author}{{The mathlib Community}},
\newblock \bibinfo{title}{The {Lean} mathematical library},
\newblock in: \bibinfo{booktitle}{Proceedings of the 9th ACM SIGPLAN
  International Conference on Certified Programs and Proofs}, CPP 2020,
  \bibinfo{publisher}{ACM}, \bibinfo{year}{2020}, pp.
  \bibinfo{pages}{367--381}. \DOIprefix\doi{10.1145/3372885.3373824}.
\bibitem[{Baumgaertel and Winter(1989)}]{baumgaertel1989determination}
\bibinfo{author}{M.~Baumgaertel}, \bibinfo{author}{H.~H. Winter},
\newblock \bibinfo{title}{Determination of discrete relaxation and retardation
  time spectra from dynamic mechanical data},
\newblock \bibinfo{journal}{Rheologica Acta} \bibinfo{volume}{28}
  (\bibinfo{year}{1989}) \bibinfo{pages}{511--519}.
  \DOIprefix\doi{10.1007/BF01332922}.
\bibitem[{Tschoegl(1989)}]{tschoegl1989phenomenological}
\bibinfo{author}{N.~W. Tschoegl}, \bibinfo{title}{The Phenomenological Theory
  of Linear Viscoelastic Behavior: An Introduction}, \bibinfo{edition}{1} ed.,
  \bibinfo{publisher}{Springer Berlin, Heidelberg}, \bibinfo{year}{1989}.
  \URLprefix \url{https://doi.org/10.1007/978-3-642-73602-5}.
  \DOIprefix\doi{10.1007/978-3-642-73602-5}.
\bibitem[{Baumgaertel and Winter(1992)}]{baumgaertel1992interrelation}
\bibinfo{author}{M.~Baumgaertel}, \bibinfo{author}{H.~Winter},
\newblock \bibinfo{title}{Interrelation between continuous and discrete
  relaxation time spectra},
\newblock \bibinfo{journal}{Journal of Non-Newtonian Fluid Mechanics}
  \bibinfo{volume}{44} (\bibinfo{year}{1992}) \bibinfo{pages}{15--36}.
  \DOIprefix\doi{10.1016/0377-0257(92)80043-W}.
\bibitem[{Park and Schapery(1999)}]{park1999methods}
\bibinfo{author}{S.~Park}, \bibinfo{author}{R.~Schapery},
\newblock \bibinfo{title}{Methods of interconversion between linear
  viscoelastic material functions. part i—a numerical method based on prony
  series},
\newblock \bibinfo{journal}{International Journal of Solids and Structures}
  \bibinfo{volume}{36} (\bibinfo{year}{1999}) \bibinfo{pages}{1653--1675}.
  \DOIprefix\doi{10.1016/S0020-7683(98)00055-9}.
\bibitem[{Beylkin and Monzón(2005)}]{beylkin2005approximation}
\bibinfo{author}{G.~Beylkin}, \bibinfo{author}{L.~Monzón},
\newblock \bibinfo{title}{On approximation of functions by exponential sums},
\newblock \bibinfo{journal}{Applied and Computational Harmonic Analysis}
  \bibinfo{volume}{19} (\bibinfo{year}{2005}) \bibinfo{pages}{17--48}.
  \DOIprefix\doi{10.1016/j.acha.2005.01.003}.
\bibitem[{Huang et~al.(2025)Huang, Nelsen, and Trautner}]{huang2025operator}
\bibinfo{author}{D.~Z. Huang}, \bibinfo{author}{N.~H. Nelsen},
  \bibinfo{author}{M.~Trautner},
\newblock \bibinfo{title}{An operator learning perspective on
  parameter-to-observable maps},
\newblock \bibinfo{journal}{Foundations of Data Science} \bibinfo{volume}{7}
  (\bibinfo{year}{2025}) \bibinfo{pages}{163--225}.
  \DOIprefix\doi{10.3934/fods.2024037}.
\bibitem[{Helnwein(2001)}]{helnwein2001some}
\bibinfo{author}{P.~Helnwein},
\newblock \bibinfo{title}{Some remarks on the compressed matrix representation
  of symmetric second-order and fourth-order tensors},
\newblock \bibinfo{journal}{Computer Methods in Applied Mechanics and
  Engineering} \bibinfo{volume}{190} (\bibinfo{year}{2001})
  \bibinfo{pages}{2753--2770}. \DOIprefix\doi{10.1016/S0045-7825(00)00263-2}.
\bibitem[{Ciarlet(2021)}]{ciarlet2021mathematical}
\bibinfo{author}{P.~G. Ciarlet}, \bibinfo{title}{Mathematical Elasticity:
  Three-Dimensional Elasticity}, Classics in Applied Mathematics,
  \bibinfo{publisher}{Society for Industrial and Applied Mathematics},
  \bibinfo{address}{Philadelphia, PA}, \bibinfo{year}{2021}.
  \DOIprefix\doi{10.1137/1.9781611976786}.
\bibitem[{Ambrosio et~al.(2000)Ambrosio, Fusco, and
  Pallara}]{ambrosio2000functions}
\bibinfo{author}{L.~Ambrosio}, \bibinfo{author}{N.~Fusco},
  \bibinfo{author}{D.~Pallara}, \bibinfo{title}{Functions of Bounded Variation
  and Free Discontinuity Problems}, \bibinfo{publisher}{Oxford University
  Press}, \bibinfo{year}{2000}.
  \DOIprefix\doi{10.1093/oso/9780198502456.001.0001}.
\bibitem[{Chen and Chen(1995)}]{chen1995universal}
\bibinfo{author}{T.~Chen}, \bibinfo{author}{H.~Chen},
\newblock \bibinfo{title}{Universal approximation to nonlinear operators by
  neural networks with arbitrary activation functions and its application to
  dynamical systems},
\newblock \bibinfo{journal}{IEEE Transactions on Neural Networks}
  \bibinfo{volume}{6} (\bibinfo{year}{1995}) \bibinfo{pages}{911--917}.
  \DOIprefix\doi{10.1109/72.392253}.
\bibitem[{Bhattacharya et~al.(2021)Bhattacharya, Hosseini, Kovachki, and
  Stuart}]{bhattacharya2021pca}
\bibinfo{author}{K.~Bhattacharya}, \bibinfo{author}{B.~Hosseini},
  \bibinfo{author}{N.~B. Kovachki}, \bibinfo{author}{A.~M. Stuart},
\newblock \bibinfo{title}{Model reduction and neural networks for parametric
  {PDEs}},
\newblock \bibinfo{journal}{The SMAI Journal of Computational Mathematics}
  \bibinfo{volume}{7} (\bibinfo{year}{2021}) \bibinfo{pages}{121--157}.
  \DOIprefix\doi{10.5802/smai-jcm.74}.
\bibitem[{Kovachki et~al.(2021)Kovachki, Lanthaler, and
  Mishra}]{kovachki2021universal}
\bibinfo{author}{N.~Kovachki}, \bibinfo{author}{S.~Lanthaler},
  \bibinfo{author}{S.~Mishra},
\newblock \bibinfo{title}{On universal approximation and error bounds for
  fourier neural operators},
\newblock \bibinfo{journal}{Journal of Machine Learning Research}
  \bibinfo{volume}{22} (\bibinfo{year}{2021}) \bibinfo{pages}{1--76}.
  \URLprefix \url{http://jmlr.org/papers/v22/21-0806.html}.
\bibitem[{Nelsen and Stuart(2021)}]{nelsen2021random}
\bibinfo{author}{N.~H. Nelsen}, \bibinfo{author}{A.~M. Stuart},
\newblock \bibinfo{title}{The random feature model for input-output maps
  between banach spaces},
\newblock \bibinfo{journal}{SIAM Journal on Scientific Computing}
  \bibinfo{volume}{43} (\bibinfo{year}{2021}) \bibinfo{pages}{A3212--A3243}.
  \DOIprefix\doi{10.1137/20M133957X}.
\bibitem[{Luo et~al.(2026)Luo, O'Leary-Roseberry, Chen, and
  Ghattas}]{luo2025dimension}
\bibinfo{author}{D.~Luo}, \bibinfo{author}{T.~O'Leary-Roseberry},
  \bibinfo{author}{P.~Chen}, \bibinfo{author}{O.~Ghattas},
\newblock \bibinfo{title}{Dimension reduction for derivative-informed operator
  learning: An analysis of approximation errors},
\newblock \bibinfo{journal}{Journal of Machine Learning Research}
  \bibinfo{volume}{27} (\bibinfo{year}{2026}) \bibinfo{pages}{1--94}.
  \URLprefix \url{https://www.jmlr.org/papers/v27/25-1861.html}.
\bibitem[{Trefethen(2019)}]{trefethen2019approximation}
\bibinfo{author}{L.~N. Trefethen}, \bibinfo{title}{Approximation Theory and
  Approximation Practice}, \bibinfo{edition}{extended edition} ed.,
  \bibinfo{publisher}{Society for Industrial and Applied Mathematics},
  \bibinfo{address}{Philadelphia, PA}, \bibinfo{year}{2019}.
  \DOIprefix\doi{10.1137/1.9781611975949}.
\bibitem[{Asaro(1983)}]{asaro1983crystal}
\bibinfo{author}{R.~J. Asaro},
\newblock \bibinfo{title}{Crystal plasticity},
\newblock \bibinfo{journal}{Journal of Applied Mechanics} \bibinfo{volume}{50}
  (\bibinfo{year}{1983}) \bibinfo{pages}{921--934}.
  \DOIprefix\doi{10.1115/1.3167205}.
\bibitem[{Roters et~al.(2010)Roters, Eisenlohr, Hantcherli, Tjahjanto, Bieler,
  and Raabe}]{roters2010overview}
\bibinfo{author}{F.~Roters}, \bibinfo{author}{P.~Eisenlohr},
  \bibinfo{author}{L.~Hantcherli}, \bibinfo{author}{D.~D. Tjahjanto},
  \bibinfo{author}{T.~R. Bieler}, \bibinfo{author}{D.~Raabe},
\newblock \bibinfo{title}{Overview of constitutive laws, kinematics,
  homogenization and multiscale methods in crystal plasticity finite-element
  modeling: Theory, experiments, applications},
\newblock \bibinfo{journal}{Acta Materialia} \bibinfo{volume}{58}
  (\bibinfo{year}{2010}) \bibinfo{pages}{1152--1211}.
  \DOIprefix\doi{10.1016/j.actamat.2009.10.058}.
\bibitem[{Roters et~al.(2019)Roters, Diehl, Shanthraj, Eisenlohr, Reuber, Wong,
  Maiti, Ebrahimi, Hochrainer, Fabritius, Nikolov, Friak, Fujita, Grilli,
  Janssens, Jia, Kok, Ma, Meier, Werner, Stricker, Weygand, and
  Raabe}]{roters2019damask}
\bibinfo{author}{F.~Roters}, \bibinfo{author}{M.~Diehl},
  \bibinfo{author}{P.~Shanthraj}, \bibinfo{author}{P.~Eisenlohr},
  \bibinfo{author}{C.~Reuber}, \bibinfo{author}{S.~L. Wong},
  \bibinfo{author}{T.~Maiti}, \bibinfo{author}{A.~Ebrahimi},
  \bibinfo{author}{T.~Hochrainer}, \bibinfo{author}{H.-O. Fabritius},
  \bibinfo{author}{S.~Nikolov}, \bibinfo{author}{M.~Friak},
  \bibinfo{author}{N.~Fujita}, \bibinfo{author}{N.~Grilli},
  \bibinfo{author}{K.~G.~F. Janssens}, \bibinfo{author}{N.~Jia},
  \bibinfo{author}{P.~J.~J. Kok}, \bibinfo{author}{D.~Ma},
  \bibinfo{author}{F.~Meier}, \bibinfo{author}{E.~Werner},
  \bibinfo{author}{M.~Stricker}, \bibinfo{author}{D.~Weygand},
  \bibinfo{author}{D.~Raabe},
\newblock \bibinfo{title}{{DAMASK} -- the {D}{\"u}sseldorf advanced material
  simulation kit for modelling multi-physics crystal plasticity, damage, and
  thermal phenomena from the single crystal up to the component scale},
\newblock \bibinfo{journal}{Computational Materials Science}
  \bibinfo{volume}{158} (\bibinfo{year}{2019}) \bibinfo{pages}{420--478}.
  \DOIprefix\doi{10.1016/j.commatsci.2018.04.030}.
\bibitem[{Bhattacharya et~al.(2026)Bhattacharya, Cao, and
  Stuart}]{bhattacharya2026optimal}
\bibinfo{author}{K.~Bhattacharya}, \bibinfo{author}{L.~Cao},
  \bibinfo{author}{A.~Stuart},
\newblock \bibinfo{title}{Optimal experimental design for reliable learning of
  history-dependent constitutive laws},
\newblock \bibinfo{journal}{Computer Methods in Applied Mechanics and
  Engineering} \bibinfo{volume}{457} (\bibinfo{year}{2026})
  \bibinfo{pages}{119022}. \DOIprefix\doi{10.1016/j.cma.2026.119022}.
\bibitem[{Zhang and Bhattacharya(2024)}]{zhang2024iterated}
\bibinfo{author}{Y.~Zhang}, \bibinfo{author}{K.~Bhattacharya},
\newblock \bibinfo{title}{Iterated learning and multiscale modeling of
  history-dependent architectured metamaterials},
\newblock \bibinfo{journal}{Mechanics of Materials} \bibinfo{volume}{197}
  (\bibinfo{year}{2024}) \bibinfo{pages}{105090}.
  \DOIprefix\doi{10.1016/j.mechmat.2024.105090}.
\bibitem[{Kammler(1976)}]{kammler1976chebyshev}
\bibinfo{author}{D.~W. Kammler},
\newblock \bibinfo{title}{Chebyshev approximation of completely monotonic
  functions by sums of exponentials},
\newblock \bibinfo{journal}{SIAM Journal on Numerical Analysis}
  \bibinfo{volume}{13} (\bibinfo{year}{1976}) \bibinfo{pages}{761--774}.
  \DOIprefix\doi{10.1137/0713062}.
\bibitem[{Braess and Hackbusch(2005)}]{braess2005approximation}
\bibinfo{author}{D.~Braess}, \bibinfo{author}{W.~Hackbusch},
\newblock \bibinfo{title}{Approximation of {$1/x$} by exponential sums in
  {$[1,\infty)$}},
\newblock \bibinfo{journal}{IMA Journal of Numerical Analysis}
  \bibinfo{volume}{25} (\bibinfo{year}{2005}) \bibinfo{pages}{685--697}.
  \DOIprefix\doi{10.1093/imanum/dri015}.
\bibitem[{Beylkin and Monz{\'o}n(2010)}]{beylkin2010revisited}
\bibinfo{author}{G.~Beylkin}, \bibinfo{author}{L.~Monz{\'o}n},
\newblock \bibinfo{title}{Approximation by exponential sums revisited},
\newblock \bibinfo{journal}{Applied and Computational Harmonic Analysis}
  \bibinfo{volume}{28} (\bibinfo{year}{2010}) \bibinfo{pages}{131--149}.
  \DOIprefix\doi{10.1016/j.acha.2009.08.011}.
\bibitem[{Koyama(2023)}]{koyama2023exponential}
\bibinfo{author}{Y.~M. Koyama}, \bibinfo{title}{Exponential sum approximations
  of finite completely monotonic functions}, \bibinfo{year}{2023}. \URLprefix
  \url{https://arxiv.org/abs/2301.08931}.
  \DOIprefix\doi{10.48550/arXiv.2301.08931}.
  \href{http://arxiv.org/abs/2301.08931}{{\tt arXiv:2301.08931}}.
\bibitem[{Darrow and Stepaniants(2026)}]{darrow2026spectral}
\bibinfo{author}{D.~Darrow}, \bibinfo{author}{G.~Stepaniants},
\newblock \bibinfo{title}{A spectral theory of scalar {Volterra} equations},
\newblock \bibinfo{journal}{Communications of the American Mathematical
  Society} \bibinfo{volume}{6} (\bibinfo{year}{2026})
  \bibinfo{pages}{634--725}. \DOIprefix\doi{10.1090/cams/64}.
\bibitem[{Romero and Ortiz(2026)}]{romero2026identification}
\bibinfo{author}{I.~Romero}, \bibinfo{author}{M.~Ortiz},
\newblock \bibinfo{title}{Identification of optimal history variables and
  corresponding hereditary laws in linear viscoelasticity},
\newblock \bibinfo{journal}{Computer Methods in Applied Mechanics and
  Engineering} \bibinfo{volume}{461} (\bibinfo{year}{2026})
  \bibinfo{pages}{119122}. \DOIprefix\doi{10.1016/j.cma.2026.119122}.
\bibitem[{Lanthaler and Stuart(2026)}]{lanthaler2026parametric}
\bibinfo{author}{S.~Lanthaler}, \bibinfo{author}{A.~M. Stuart},
\newblock \bibinfo{title}{The parametric complexity of operator learning},
\newblock \bibinfo{journal}{IMA Journal of Numerical Analysis}
  \bibinfo{volume}{46} (\bibinfo{year}{2026}) \bibinfo{pages}{647--712}.
  \DOIprefix\doi{10.1093/imanum/draf028}.
\bibitem[{Kovachki et~al.(2024)Kovachki, Lanthaler, and
  Mhaskar}]{kovachki2024data}
\bibinfo{author}{N.~B. Kovachki}, \bibinfo{author}{S.~Lanthaler},
  \bibinfo{author}{H.~Mhaskar}, \bibinfo{title}{Data complexity estimates for
  operator learning}, \bibinfo{year}{2024}. \URLprefix
  \url{https://arxiv.org/abs/2405.15992}.
  \DOIprefix\doi{10.48550/arXiv.2405.15992}.
  \href{http://arxiv.org/abs/2405.15992}{{\tt arXiv:2405.15992}}.
\bibitem[{Masi and Stefanou(2023)}]{masi2023evolution}
\bibinfo{author}{F.~Masi}, \bibinfo{author}{I.~Stefanou},
\newblock \bibinfo{title}{Evolution {TANN} and the identification of internal
  variables and evolution equations in solid mechanics},
\newblock \bibinfo{journal}{Journal of the Mechanics and Physics of Solids}
  \bibinfo{volume}{174} (\bibinfo{year}{2023}) \bibinfo{pages}{105245}.
  \DOIprefix\doi{10.1016/j.jmps.2023.105245}.
\bibitem[{Ortiz and Pedregal(2026)}]{ortiz2026encoding}
\bibinfo{author}{M.~Ortiz}, \bibinfo{author}{P.~Pedregal},
  \bibinfo{title}{Optimal history encoding for elastic-plastic hereditary laws:
  Sharp input and constitutive approximation}, \bibinfo{year}{2026}. \URLprefix
  \url{https://arxiv.org/abs/2607.09974}.
  \DOIprefix\doi{10.48550/arXiv.2607.09974}.
  \href{http://arxiv.org/abs/2607.09974}{{\tt arXiv:2607.09974}},
  \bibinfo{note}{preprint}.
\bibitem[{Huang et~al.(2022)Huang, He, Chem, and Reina}]{huang2022variational}
\bibinfo{author}{S.~Huang}, \bibinfo{author}{Z.~He}, \bibinfo{author}{B.~Chem},
  \bibinfo{author}{C.~Reina},
\newblock \bibinfo{title}{Variational {Onsager} neural networks ({VONNs}): A
  thermodynamics-based variational learning strategy for non-equilibrium
  {PDEs}},
\newblock \bibinfo{journal}{Journal of the Mechanics and Physics of Solids}
  \bibinfo{volume}{163} (\bibinfo{year}{2022}) \bibinfo{pages}{104856}.
  \DOIprefix\doi{10.1016/j.jmps.2022.104856}.
\bibitem[{Logg et~al.(2012)Logg, Mardal, and Wells}]{logg2012automated}
\bibinfo{editor}{A.~Logg}, \bibinfo{editor}{K.-A. Mardal},
  \bibinfo{editor}{G.~Wells} (Eds.), \bibinfo{title}{Automated Solution of
  Differential Equations by the Finite Element Method},
  volume~\bibinfo{volume}{84} of \textit{\bibinfo{series}{Lecture Notes in
  Computational Science and Engineering}}, \bibinfo{publisher}{Springer},
  \bibinfo{address}{Berlin, Heidelberg}, \bibinfo{year}{2012}.
  \DOIprefix\doi{10.1007/978-3-642-23099-8}.
\bibitem[{Aln{\ae}s et~al.(2015)Aln{\ae}s, Blechta, Hake, Johansson, Kehlet,
  Logg, Richardson, Ring, Rognes, and Wells}]{alnaes2015fenics}
\bibinfo{author}{M.~S. Aln{\ae}s}, \bibinfo{author}{J.~Blechta},
  \bibinfo{author}{J.~Hake}, \bibinfo{author}{A.~Johansson},
  \bibinfo{author}{B.~Kehlet}, \bibinfo{author}{A.~Logg},
  \bibinfo{author}{C.~Richardson}, \bibinfo{author}{J.~Ring},
  \bibinfo{author}{M.~E. Rognes}, \bibinfo{author}{G.~N. Wells},
\newblock \bibinfo{title}{The {FEniCS} project version 1.5},
\newblock \bibinfo{journal}{Archive of Numerical Software} \bibinfo{volume}{3}
  (\bibinfo{year}{2015}) \bibinfo{pages}{9--23}.
  \DOIprefix\doi{10.11588/ans.2015.100.20553}.
\bibitem[{Paszke et~al.(2019)Paszke, Gross, Massa, Lerer, Bradbury, Chanan,
  Killeen, Lin, Gimelshein, Antiga, Desmaison, K{\"o}pf, Yang, DeVito, Raison,
  Tejani, Chilamkurthy, Steiner, Fang, Bai, and Chintala}]{paszke2019pytorch}
\bibinfo{author}{A.~Paszke}, \bibinfo{author}{S.~Gross},
  \bibinfo{author}{F.~Massa}, \bibinfo{author}{A.~Lerer},
  \bibinfo{author}{J.~Bradbury}, \bibinfo{author}{G.~Chanan},
  \bibinfo{author}{T.~Killeen}, \bibinfo{author}{Z.~Lin},
  \bibinfo{author}{N.~Gimelshein}, \bibinfo{author}{L.~Antiga},
  \bibinfo{author}{A.~Desmaison}, \bibinfo{author}{A.~K{\"o}pf},
  \bibinfo{author}{E.~Yang}, \bibinfo{author}{Z.~DeVito},
  \bibinfo{author}{M.~Raison}, \bibinfo{author}{A.~Tejani},
  \bibinfo{author}{S.~Chilamkurthy}, \bibinfo{author}{B.~Steiner},
  \bibinfo{author}{L.~Fang}, \bibinfo{author}{J.~Bai},
  \bibinfo{author}{S.~Chintala},
\newblock \bibinfo{title}{{PyTorch}: An imperative style, high-performance deep
  learning library},
\newblock in: \bibinfo{booktitle}{Advances in Neural Information Processing
  Systems}, volume~\bibinfo{volume}{32}, \bibinfo{publisher}{Curran Associates,
  Inc.}, \bibinfo{year}{2019}, pp. \bibinfo{pages}{8024--8035}. \URLprefix
  \url{https://proceedings.neurips.cc/paper/2019/hash/bdbca288fee7f92f2bfa9f7012727740-Abstract.html}.

\end{thebibliography}
\end{document}